\input{epsf}

\newcounter{letter}
\newenvironment{alphalist}{\begin{list}
{{\normalshape(\alph{letter})}}{\usecounter{letter}}}{\end{list}}

\newcommand{\two}{I\!I}
\newcommand{\T}{\cal T}
\newcommand{\C}{\cal C}
\newcommand{\D}{\cal D}
\newcommand{\A}{\cal A}
\newcommand{\B}{\cal B}
\newcommand{\tensor}{\otimes}
\newcommand{\ten}{\otimes}
\newcommand{\intersect}{\cap}
\newcommand{\iso}{\cong}
\newcommand{\union}{\cup}

\newcommand{\btensor}{\bigotimes}
\newcommand{\tto}{\Rightarrow}
\newcommand{\maps}{\colon}
\newcommand{\from}{\colon}
\newcommand{\R}{{\Bbb R}}
\newcommand{\N}{{\Bbb N}}

\newtheorem{thm}{Theorem} 

\newtheorem{lem}[thm]{Lemma}

\newtheorem{defn}[thm]{Definition}

\newcommand{\et}{\hspace{-0.08in}{\bf .}\hspace{0.1in}}
\newcommand{\BOX}{\hbox {$\sqcap$ \kern -1em $\sqcup$}}
\newcommand{\qed}{\hskip 3em \hbox{\BOX} \vskip 2ex}

\newcommand{\source} {{\rm source}}
\newcommand{\target}{{\rm target}}
\renewcommand{\ten}{\otimes}
\newcommand{\rarr}{\rightarrow}
\newcommand{\bull}{\bullet}

        \newcommand{\be}{\begin{equation}}
        \newcommand{\ee}{\end{equation}}
        \newcommand{\ba}{\begin{eqnarray}}
        \newcommand{\ea}{\end{eqnarray}}
        \newcommand{\ban}{\begin{eqnarray*}}
        \newcommand{\ean}{\end{eqnarray*}}
        \newcommand{\barr}{\begin{array}}
        \newcommand{\earr}{\end{array}}

\documentstyle[12pt,diagram,amsfonts]{article}
\begin{document}

      \begin{center}
      {\bf Higher-Dimensional Algebra IV: \\
        2-Tangles \\}
      \vspace{0.5cm}
      {\em John C.\ Baez \\}
      \vspace{0.3cm}
      {\small Department of Mathematics,  University of California\\
      Riverside, California 92521 \\
      USA\\ }
      \vspace{0.3cm}
      {\em Laurel Langford \\}
      \vspace{0.3cm}
      {\small Department of Mathematics, University of Wisconsin ---
       River Falls\\
      410 South Third St.\ \\
      River Falls, Wisconsin 54022 \\
      USA \\}
      \vspace{0.3cm}
      {\small email: baez@math.ucr.edu, laurel.langford@uwrf.edu\\}
      \vspace{0.3cm}
      {\small November 23, 1998 \\ }
      \end{center}

\begin{abstract}
Just as knots and links can be algebraically described as certain
morphisms in the category of tangles in 3 dimensions, compact surfaces
smoothly embedded in $\R^4$ can be described as certain 2-morphisms in
the 2-category of `2-tangles in 4 dimensions'.   Using the work of
Carter, Rieger and Saito, we prove that this 2-category is the `free
semistrict braided monoidal 2-category with duals on one unframed
self-dual object'.  By this universal property, any unframed self-dual
object in a braided monoidal 2-category with duals determines an
invariant of 2-tangles in 4 dimensions.
\end{abstract}

\section{Introduction}

One of the most exciting aspects of higher-dimensional algebra is how
weak $n$-categories seem to provide precisely the right mathematics for
algebraic topology.  From one point of view, weak $n$-categories are
purely algebraic structures consisting of objects, 1-morphisms between
objects, 2-morphisms between 1-morphisms, and so on up to $n$-morphisms,
together with various composition operations, satisfying laws that arise
naturally from algebraic considerations \cite{B2}.   But time and time
again, the mathematics of weak $n$-categories has turned out to be
perfectly suited to $n$-dimensional topology.

Until the late 1980's, the most striking instances of this phenomenon
came from homotopy theory.  By now there is a large body of evidence
supporting a conjecture that would completely explain the relation
between $n$-categories and homotopy theory \cite{BD3}.  In rough terms,
this conjecture states that spaces with vanishing homotopy groups above
dimension $n$ are equivalent to a certain class of weak $n$-categories,
the `weak $n$-groupoids'.   An weak $n$-groupoid is a weak $n$-category
where every $j$-morphism has a `weak inverse'.  For $j = n$, a weak
inverse for the $j$-morphism $f \maps x \to y$ is just an inverse in the
usual sense, while for $j < n$, a weak inverse for $f$ is recursively
defined to be a $j$-morphism $g \maps x \to y$ such that $fg$ and $gf$
are the identity {\it up to a weakly invertible $(j+1)$-morphism}.

Starting with the discovery of the Jones polynomial and a family of
related knot invariants, a new branch of algebraic topology has emerged
in the last decade.  It is often called `quantum topology' because of
its close ties to quantum field theory.  Its relation to more
traditional forms of algebraic topology based on homotopy theory was
initially very mysterious.  Now it appears that quantum topology goes
beyond homotopy theory precisely by exploiting a larger class of
$n$-categories, the `$n$-categories with duals'.   These are poorly
understood except in some low-dimensional cases, but some of their
essential features are already clear.  Most importantly, while every
$j$-morphism $f \maps x \to y$ has a `dual' $f^\ast \maps y \to x$, this
dual need not be a weak inverse of $f$.  One important example of a
category with duals is the category of Hilbert spaces, where the dual of
a linear operator is its Hilbert space adjoint.   Another is the
category of tangles in 3-dimensional space, where the dual of a tangle
is obtained by reflecting it to switch its source and target.

The category of tangles in 3 dimensions is especially important, because
it has a beautiful algebraic characterization in terms of a universal
property.  This was initially developed by Turaev \cite{T},
Freyd--Yetter \cite{FY,Y}, and Joyal--Street \cite{JS}, and it reached a
highly polished form in the work of Shum \cite{Shum}.  In our language
\cite{B,BD}, her result is that isotopy classes of framed oriented
tangles in 3 dimensions are the morphisms of the `free braided monoidal
category with duals on one object'.   Using this universal property, we
can easily obtain functors from this category to other braided monoidal
categories with duals, such as categories of representations of quantum
groups.  Any such functor gives an invariant of tangles, and in
particular, a knot invariant.   This is the easiest way to understand
the Jones polynomial and its relatives \cite{RT}.

The `tangle hypothesis' \cite{BD} suggests a vast generalization of this
result, applicable to $n$-dimensional surfaces embedded in
$(n+k)$-dimensional space for all $n$ and $k$.  This generalization
involves the notion of a `$k$-tuply monoidal $n$-category'.  A $k$-tuply
monoidal $n$-category is an $(n+k)$-category that has only trivial
$j$-morphisms for $j < k$.  By reindexing we can think of this as an
$n$-category with extra structure and properties.  Some low-dimensional
special cases are shown in Table 1 below.

Briefly put, the tangle hypothesis says that framed oriented $n$-tangles
in $(n+k)$ dimensions are the $n$-morphisms of the `free weak $k$-tuply
monoidal $n$-category with duals on one object'.  In this $n$-category,
the objects correspond to collections of points embedded in $[0,1]^k$.
The 1-morphisms correspond to compact 1-manifolds with boundary embedded
in $[0,1]^{k+1}$ going from one such object to another.  Similarly, the
2-morphisms correspond to compact 2-manifolds with corners embedded in
$[0,1]^{k+2}$ going from one 1-morphism $f \maps x \to y$ to another
1-morphism $g \maps x \to y$, and so on.  Finally, the $n$-morphisms
correspond to isotopy classes of $n$-manifolds with corners embedded in
$[0,1]^{n+k}$.  We call these `$n$-tangles in $n+k$ dimensions'.

\begin{center}
{\small
\begin{tabular}{|c|c|c|c|}  \hline
         & $n = 0$   & $n = 1$    & $n = 2$          \\     \hline
$k = 0$  & sets      & categories & 2-categories     \\     \hline
$k = 1$  & monoids   & monoidal   & monoidal         \\
         &           & categories & 2-categories     \\     \hline
$k = 2$  &commutative& braided    & braided          \\
         & monoids   & monoidal   & monoidal         \\
         &           & categories & 2-categories     \\     \hline
$k = 3$  &`'         & symmetric  & weakly involutory \\
         &           & monoidal   & monoidal         \\
         &           & categories & 2-categories     \\     \hline
$k = 4$  &`'         & `'         &strongly involutory\\
         &           &            & monoidal         \\
         &           &            & 2-categories     \\     \hline
$k = 5$  &`'         &`'          & `'               \\
         &           &            &                  \\
         &           &            &                  \\      \hline
\end{tabular}} \vskip 1em
Table 1.  $k$-tuply monoidal $n$-categories
\end{center}
\vskip 0.5em

Unfortunately the tangle hypothesis involves concepts from topology and
$n$-category theory that presently have only been made precise in certain
low-dimensional cases.  As a kind of warmup, we wish to prove a version
of this hypothesis in the case $n = k = 2$.  So far we have only
completed work on the unframed, unoriented case, which allows us to take
maximal advantage of the recent work of Carter, Rieger and Saito
\cite{CRS}.  Since the theory of $k$-tuply monoidal weak $n$-categories
is not yet well developed for $n = k = 2$, we use the better-understood
`semistrict' ones as a kind of stopgap.  These are also known as
`semistrict braided monoidal 2-categories'.  Our result is thus that the
2-category of unframed unoriented 2-tangles in 4 dimensions is the `free
semistrict braided monoidal 2-category with duals on one unframed
self-dual object'.

This result is closely related to the fact that the category of unframed
unoriented tangles in 3 dimensions is the free braided monoidal category
with duals on one unframed self-dual object.   In particular, the
Reidemeister moves, which arise as {\it equations} between morphisms in the
category of tangles in 3 dimensions, arise as {\it 2-isomorphisms} in our
context.  For this reason we say that our result is a `categorification'
of the 3-dimensional one.  For more on categorification and how it
relates to the tangle hypothesis, see our previous papers \cite{BD,BD3}.

The study of duality in $n$-categories is only beginning, so an
important part of this paper consists of finding an appropriate
definition of a braided monoidal 2-category `with duals'.  Given this,
we simply define a `self-dual' object $x$ to be one for which $x = x^*$.
On the other hand, the notion of an `unframed' object really takes
advantage of categorification.  In the category of tangles, a twist in
the framing corresponds to a morphism called the `balancing'.  In our
situation, an `unframed object' is not one for which the balancing {\it
equals} the identity, but one for which the balancing is {\it
isomorphic} to the identity via a certain 2-isomorphism.  This
2-isomorphism, corresponding to the Reidemeister I move, satisfies a
highly nontrivial equation of its own.

The study of universal properties for $n$-categories is also just
beginning, so we must clarify what is meant by the `free' braided
monoidal 2-category with duals on one unframed self-dual object.
Finally, since there is presently no general construction of the
$k$-tuply monoidal $n$-category of $k$-tangles in $n$ dimensions, we
must construct this `by hand' in the case $n = k = 2$ before proving our
result.  To obtain a semistrict braided monoidal 2-category, we cannot
let the objects be simply collections of points embedded in $[0,1]^2$.
Instead, we must introduce an equivalence relation on such collections,
and take objects to be equivalence classes.  Similarly, the morphisms in
our 2-category are certain equivalence classes of tangles.

We must choose these equivalence relations carefully, in order to avoid
the errors present in Fischer's attempt \cite{Fischer} to define a
2-category of 2-tangles.   Kharlamov and Turaev \cite{KT} have shown
that composition of 2-morphisms is not well-defined if, as Fischer did,
we take isotopy classes of tangles as our 1-morphisms.    Kharlamov and
Turaev showed how to avoid this problem by introducing a `height
function' on $[0,1]^3$ and saying that two tangles define the same
1-morphism only if they differ by an isotopy that preserves the order of
the heights of local extrema of this function.    Our work is based on
Carter, Rieger and Saito's recent combinatorial description of 2-tangles
\cite{CRS}, which places a somewhat stronger restriction on the
isotopies: they must preserve the order of heights of local maxima,
local minima, and crossings relative to a specified projection.

The plan of the paper is as follows.   In Section 2 we give a
topological description of a $2$-category $\T$ of unframed unoriented
2-tangles in 4 dimensions.  We define duality for monoidal and braided
monoidal 2-categories, and show that $\T$ is a braided monoidal
$2$-category with duals.   We also define the notion of an `unframed
self-dual object' in a braided monoidal 2-category with duals, and we
show that $\T$ has an unframed self-dual object $Z$ corresponding to a
single point in the unit square.  In Section 3 we give an alternate,
purely combinatorial description of a 2-category of $2$-tangles, which
we denote by $\C$.  Using the work of Carter, Rieger and Saito, we then
show that $\T$ and $\C$ are isomorphic.  In Section 4 we use this
isomorphism to show that $\T$ is generated, as a braided monoidal
$2$-category with duals, by the unframed self-dual object $Z$.   Given a
strict monoidal 2-functor $F \maps \T \to \B$, we define what it means
for $F$ to `preserve braiding and duals semistrictly on the generator'.
Finally, we show that for any braided monoidal 2-category with duals
$\B$ containing an unframed self-dual object $B$, there is a unique
strict monoidal 2-functor $F \maps \T \to \B$ with $F(Z) = B$ that
preserves braiding and duals semistrictly on the generator.  This is the
precise sense in which the 2-category of 2-tangles is the free braided
monoidal $2$-category with duals on an unframed self-dual object.

This paper is based upon the second author's Ph.D.\ thesis.  A summary
of the results here can be found in a previous paper of ours \cite{BL},
and also in the magnificently illustrated book by Carter and Saito
\cite{CS}.  In the Errata section at the end of this paper we correct
some errors in our previous one.  Also, as promised, we treat
the universal property of $\T$ more carefully here, allowing us to omit
the conditions $\tilde R_{(A|A,A)} = 1$ and $\tilde R_{(A|A,A)} = 1$
which previously appeared in the definition of an unframed self-dual
object $A$.

We refer to the paper in which the tangle hypothesis was first stated as
HDA0 \cite{BD}, and refer to the earlier papers in this series as HDA1
\cite{BN}, HDA2 \cite{B}, and HDA3 \cite{BD2}.

\section{A Topological Description of 2-Tangles}

In this section we describe the 2-category $\T$ of 2-tangles using the
language of differential topology, and prove that $\T$ is a braided
monoidal 2-category with duals.  First we carefully describe the
objects, 1-morphisms, and 2-morphisms of $\T$, and show that $\T$ has the
structure of a 2-category.  Then we show that $\T$ actually has the
structure of a braided monoidal 2-category with duals.

In what follows, all manifolds are assumed to be compact and smooth, but
possibly with corners shaped like the subset $\{(x_1,\dots,x_n) \colon
\; x_1,\dots,x_i \ge 0\}$.  All maps between them are assumed to be
smooth in the obvious sense, but not necessarily mapping corners to
corners.  A diffeomorphism, defined as a smooth map with smooth inverse,
automatically maps corners to corners.  An {\it ambient isotopy} of a
manifold $M$ is defined to be a map $H \maps M \times [0,1] \to M$ such
that $H(\cdot,s)$ is a diffeomorphism for each $s \in [0,1]$.  We define
a {\it fiber isotopy} of a bundle $p \maps E \to B$ to be a ambient
isotopy $H$ of the total space such that $H(\cdot,s)$ maps fibers to
(possibly different) fibers for each $s \in [0,1]$.

\subsection{The 2-Category $\T$}

We shall think of 2-tangles as lying in the 4-cube $I_1\times I_2\times
I_3\times I_4$, where $I_i = [0,1]$.  We take the coordinates of this
4-cube to be $x,y,z,$ and $t$, respectively.  We refer to points as
being {\it behind} or {\it in front} if they have greater or smaller $x$
values, to the {\it left} or {\it right} if they have smaller or greater
$y$ values, {\it above} or {\it below} for smaller or greater $z$
values, and {\it before} or {\it after} if they have smaller or greater
$t$ values.  We refer to the function $z$ as the {\it height}.
Sometimes a 2-tangle will be illustrated by what Carter, Rieger and
Saito \cite{CRS} call a `movie': a finite sequence of tangles that are
cross-sections of the 2-tangle at successive values of $t$.  We draw
tangles in the cube, with axes $x,y$ and $z$ as shown on the left of
Fig.\ 1.  Another way we illustrate a 2-tangle is to draw a generic
projection of it into the cube with axes $y,z,$ and $t$ as shown on
the right of Fig.\ 1.

\begin{figure}[h]
\centerline{\epsfysize=1.5in\epsfbox{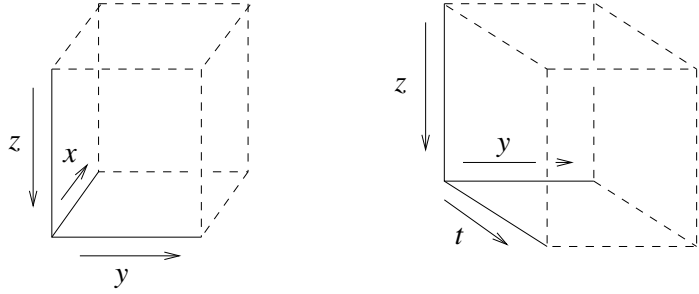}} \medskip
\caption{Orientation of illustrations}
\end{figure}

\subsubsection{Objects}

The objects of $\T$ are in one-to-one correspondence with the natural
numbers $0,1,2,\dots .$  Technically, we define objects to be equivalence
classes of generic finite subsets of the square.  We say a finite subset
of the square $I_1\times I_2$ is {\it generic} if it lies in the
interior of the square and no two points in the set have the same $y$
coordinate.  Two such sets $A$ and $B$ are said to be {\it equivalent}
if there is a level-preserving ambient isotopy $H$ of the square such
that $H(A,0)=A$ and $H(A,1)=B$, where we say $H$ is {\sl level
preserving} if it is a fiber isotopy of the bundle $p\colon I_1\times
I_2\to I_2$.  We call such an isotopy an {\it equivalence isotopy}
between generic finite subsets of the square.  This defines an
equivalence relation under which any two such sets with the same number
of elements are equivalent.

\subsubsection{1-Morphisms}

The 1-morphisms of $\T$ are equivalence classes of generic tangles.
More precisely:

\begin{defn}\et A {\rm tangle} is a 1-dimensional manifold $T$ with
boundary embedded in ${\rm int}(I_1 \times I_2) \times I_3$ such that:
\begin{enumerate}
\item The boundary points of $T$ lie in ${\rm int}
(I_1 \times I_2) \times \{0,1\}$.
\item $T$ has a {\rm product structure} near the top and bottom;
that is, there exists $\epsilon > 0$ such
that $(x,y,z)\in T$ if and only if $(x,y,z_0)\in T$
when $z$ is within $\epsilon$ of $z_0 = 0$ or $z_0 = 1$.
\end{enumerate}
We say a tangle $T$ is {\rm generic} if its projection to $I_2 \times
I_3$ is an embedding except for finitely many crossings (transverse
double points), the critical points of the height function on $T$ are
all nondegenerate local extrema, and all the crossings and
critical points occur at different heights.  \end{defn}

\noindent Note that any tangle is ambient isotopic to a generic tangle.
This allows us to restrict attention to generic tangles without loss of
generality.

\begin{defn} \label{equivalent.tangles} \et Two generic tangles are {\rm
equivalent} if there is an ambient isotopy $H$ carrying one to the other
such that:

\begin{enumerate}

\item $H$ is {\rm level preserving}, meaning that it is a fiber
isotopy for the bundle $p\colon I_1 \times I_2 \times I_3 \to
I_2 \times I_3$ and also for the bundle $\pi \colon I_1 \times I_2
\times I_3 \to I_3$.

\item $H$ has a {\rm product structure} in a neighborhood of $I_1 \times
I_2 \times \partial I_3$, meaning that for some $\epsilon > 0$,
$H(x,y,z,s)$ is of the form
\[   (X(x,y,z_0,s),Y(x,y,z_0,s),z) \]
if $z$ is within $\epsilon$ of $z_0 = 0$ or $z_0 = 1$.

\end{enumerate}
We call such an isotopy an {\rm equivalence isotopy} between generic
tangles.
\end{defn}

The level-preserving properties of this equivalence relation imply that
generic tangles whose projections onto the square differ by Reidemeister
moves or by changing the order of heights of crossings or local
extrema are not equivalent.  We may thus represent 1-morphisms of $\T$ by
planar diagrams of tangles for which the crossings and local extrema
of $z$ are ordered by their height.

Given a 1-morphism $f$ in $\T$ represented by a generic tangle $T$, we
define its {\it source} to be the object represented by to the set $T
\intersect (I_1 \times I_2 \times \{0\})$.  Similarly, we define its
{\it target} to be the object corresponding to $T \intersect (I_1 \times
I_2 \times \{1\})$.  The restriction of an equivalence isotopy between
generic tangles to $I_1\times I_2\times \{0\}$ or $I_1\times I_2\times
\{1\}$ is an equivalence isotopy between generic finite subsets of the
square, so the source and target of a 1-morphism are well-defined
objects.

\subsubsection{2-Morphisms}

Similarly, the 2-morphisms of $\T$ are equivalence classes of
generic 2-tangles.

\begin{defn}\et  A {\rm 2-tangle} is a 2-manifold $S$ with corners
embedded in $I_1 \times I_2 \times I_3 \times I_4$ such that:

\begin{enumerate}

\item The boundary of $S$ is embedded in $I_1 \times I_2 \times \partial
(I_3\times I_4)$, the intersection $S \intersect (I_1 \times I_2 \times
I_3 \times \partial I_4)$ is a pair of tangles, and $S \cap (I_1 \times
I_2 \times \partial I_3 \times I_4)$ consists of finitely many straight
lines of the form $(x,y,z)\times I_4$.

\item $S$ has a {\rm product structure} near the boundary.  That
is, there exists $\epsilon > 0$ such that if $(x,y,z,t) \in S$
then $(x,y,z',t) \in S$ if both $z$ and $z'$ are within $\epsilon$ of
either 0 or 1, and $(x,y,z,t') \in S$ if both $t$ and $'$ are within
$\epsilon$ of either 0 or 1.
\end{enumerate}
We say a 2-tangle $S$ is {\rm generic} if its intersection with the
hyperplanes of constant $t$ are generic tangles except for finitely many
values of $t$ at which one of the `full set of elementary string
interactions' occurs.  Briefly, these are: the three
Reidemeister moves, the birth or death of an unknotted circle, a saddle
point of the function $t$ on $S$, a cusp on a fold line, a double point
arc crossing a fold line, and moves that interchange the relative height
of two crossings and/or local extrema.   \end{defn}

\noindent Just as for tangles, any 2-tangle is ambient isotopic to a
generic one.  For a proof of this and a more detailed description
of the full set of elementary string interactions, see Carter, Saito
and Rieger \cite{CRS}.

\begin{defn} \label{equivalent.2-tangles} \et Two generic 2-tangles are {\rm
equivalent} if there is an ambient isotopy
\[H\maps I_1 \times I_2 \times I_3 \times I_4 \times [0,1]
\to I_1 \times I_2 \times I_3 \times I_4\]
carrying one to the other such that:

\begin{enumerate}
\item The restriction of $H$ to $I_1 \times I_2 \times I_3 \times
\{t_0\}$ for $t_0 = 0$ or $1$ is an equivalence of generic tangles.

\item The restriction of $H$ to $I_1 \times I_2\times \{z_0\} \times
I_4$ for $z_0 = 0$ or $1$ is {\rm level preserving} in the sense that it
is a fiber isotopy for the bundle $p' \from I_1\times I_2\times I_4\to
I_2\times I_4$, and also for the bundle $\pi'\from I_1\times I_2\times
I_4\to I_4$.

\item $H$ has a {\rm product structure} near $I_1 \times I_2\times \partial
(I_3 \times I_4)$; specifically, there is an $\epsilon >0$ such that
$H$ is of the form
\[H(x,y,z,t,s) = (X(x,y,z,t_0,s),Y(x,y,z,t_0,s),Z(x,y,z,t_0,s),t)\]
if $t$ is within $\epsilon$ of $t_0 = 0$ or $t_0 = 1$, and of the form
\[H(x,y,z,t,s) = (X(x,y,z_0,t,s),Y(x,y,z_0,t,s),z,T(x,y,z_0,t,s))\]
if $z$ is within $\epsilon$ of $z_0 = 0$ or $z_0 = 1$.
\end{enumerate}
We call such an isotopy an {\rm equivalence isotopy} between generic
2-tangles.
\end{defn}

Given a 2-morphism of $\T$ represented by a generic 2-tangle $S$, we
define its {\it source} to be the 1-morphism represented by the generic tangle
$S\cap (I_1 \times I_2 \times I_3 \times \{0\})$.  We define its {\it target}
to be the 1-morphism represented by $S\cap (I_1 \times I_2 \times I_3
\times \{1\})$.  The restriction of an equivalence isotopy between
generic 2-tangles to $I_1\times I_2\times I_3\times \{0\}$ or $I_1\times
I_2\times I_3\times \{1\}$ gives an equivalence isotopy between generic
tangles, so the source and target of a 2-morphism are well-defined
1-morphisms.

\begin{figure}[ht]
\centerline{\epsfysize=1.5in\epsfbox{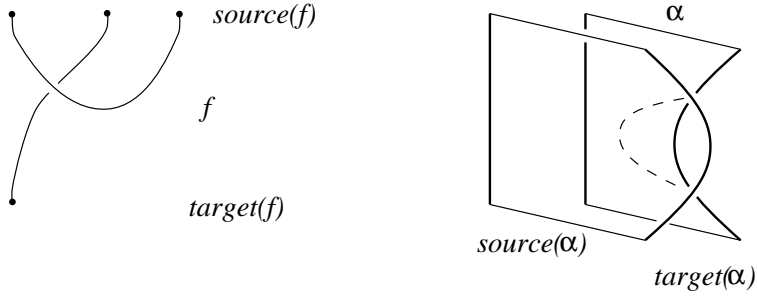}} \medskip
\caption{Morphisms and 2-morphisms}
\end{figure}

Note that generic 2-tangles without boundary are equivalent if and
only if they are ambient isotopic.  It follows that 2-morphisms in
$\T$ having the 1-morphism corresponding to the empty tangle as both
source and target are the same thing as ambient isotopy classes of
closed 2-dimensional submanifolds of $[0,1]^4$, or equivalently, $\R^4$.
This is the precise sense in which our algebraic characterization of
$\T$ gives an algebraic description of knotted surfaces in 4-dimensional
space.

In what follows we often use the same notation for an object (resp.\
morphism, 2-morphism) of $\T$ and a subset of the square (resp.\ generic
tangle, generic 2-tangle) representing it.  The difference should be
clear from the context.

\subsubsection{Composition of 1-morphisms and 2-morphisms}

We now describe how to compose 1-morphisms and 2-morphisms in $\T$.
We use the same conventions and notation regarding 2-categories as in
HDA1.  In particular, composition of 1-morphisms, horizontal
composition of a 1-morphism and a 2-morphism in either order, and
horizontal composition of 2-morphisms is denoted by
juxtaposition.  Vertical composition of 2-morphisms is denoted by
$\cdot$.  We use the ordering in which, for example, the composite of
$f \maps A \to B$ and $g \maps B \to C$ is written as $fg$.  We treat 
composition of 1-morphisms, horizontal composition of 2-morphisms, and 
vertical composition of 2-morphisms as fundamental, and define horizontal
composition of a 1-morphism and a 2-morphism in terms of these in
the usual way:
\[   f\alpha := 1_f \alpha , \qquad
     \alpha f := \alpha 1_f.  \]

Composition of 1-morphisms corresponds to gluing tangles along their
source and target sets.  Since we require that tangles be straight near
their source and target, the resulting tangle is indeed a smooth
submanifold.  Specifically, let $f\colon A\to B$ and $g\colon B\to C$ be
1-morphisms in $\T$.  Then the composite $fg\colon A\to C$ is defined as
follows.  Choose generic tangles representing $f$ and $g$, which by
abuse of language we also call $f$ and $g$.  Assume that the set
representing the target of $f$ equals the set representing the source of
$g$.  By abuse of language we call this set $B$.  Also assume that $f$
has a product structure in the $z$ direction for $z \in [1/4,1]$ ---
i.e., $f \cap I_1 \times I_2 \times \{z\} =B$ for such $z$ --- and that
$g$ has a product structure in the $z$ direction for $z \in [0,3/4]$.
Then $fg$ is the 1-morphism given by the tangle
\[(f \cap I_1 \times I_2 \times [0,1/2])\cup
(g \cap I_1 \times I_2 \times [1/2,1]).\]
The tangles $f$ and $g$ agree on $I_1 \times I_2 \times [1/4,3/4]$, so
this indeed gives a tangle, and in fact a generic tangle, as shown in Fig.\ 3.

\begin{figure}[ht]
\centerline{\epsfxsize=6in\epsfbox{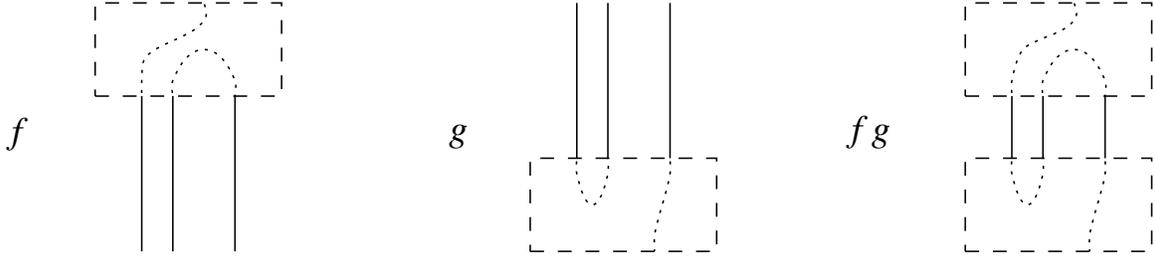}} \medskip
\caption{ Composition of 1-morphisms}
\end{figure}

Horizontal composition of 2-morphisms corresponds to gluing generic
2-tangles along the strands at their top and bottom.  Since the top and
bottom of a 2-tangle consists of finitely many straight strands, there
is a unique way to glue them together, and since we require 2-tangles to
have a product structure near their top and bottom, this gluing indeed
results in a 2-tangle, as shown in Fig.\ 4.

More precisely, suppose we have 2-morphisms $\alpha\maps f\tto f'$ and
$\beta\maps g\tto g'$ in $\T$ such that $f,f' \maps A \to B$ and $g,g'
\maps B \to C$.  Then we define their horizontal composite $\alpha
\beta$ as follows.  We choose generic 2-tangles representing
these 2-morphisms, which by abuse of language we also call $\alpha$ and
$\beta$, such that
\[ \alpha \intersect I_1 \times I_2 \times \{1\} \times I_4=
\beta \intersect I_1 \times I_2 \times \{0\} \times I_4 ,\]
$\alpha$ has a product structure in the $z$ direction for $z \in
[1/4,1]$, and $\beta$ has a product structure in the $z$ direction for
$z \in [0,3/4]$.  Then $\alpha \beta$ is the 2-morphism
represented by
\[(\alpha \intersect I_1 \times I_2 \times [0,1/2]\times I_4)\union
(\beta \intersect I_1 \times I_2 \times [1/2,1]\times I_4) .  \]

\begin{figure}[ht]
\centerline{\epsfysize=2in\epsfbox{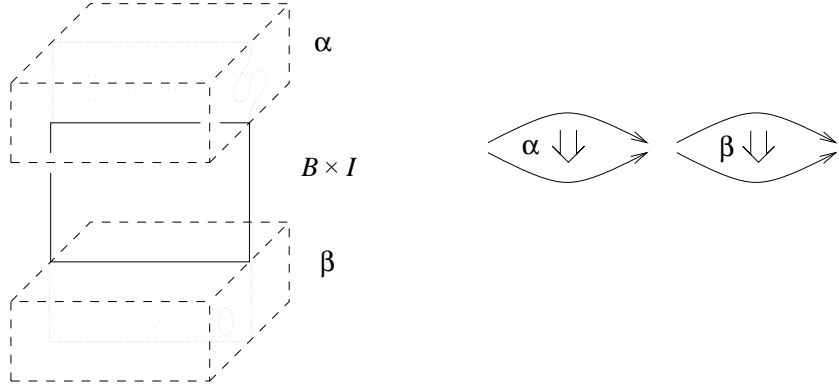}} \medskip
\caption{ Horizontal composition of 2-morphisms}
\end{figure}

Vertical composition of 2-morphisms $\alpha \maps f \tto g$ and $\beta
\maps g \tto h$ corresponds to gluing together 2-tangles representing
$\alpha$ and $\beta$ along a tangle representing $g$.  For vertical
composition to be well-defined, it is crucial that tangles with
different height orderings of their local extrema cannot represent the
same 1-morphism.  This is built into Definition
\ref{equivalent.tangles}.  Since we require 2-tangles to have a product
structure near their source and target, this gluing indeed gives a
smooth 2-tangle, as shown in Fig.\ 5.

More precisely, given 2-morphisms $\alpha \colon f\tto g$ and $\beta
\colon g\tto h$ in $\T$, we define the vertical composite $\alpha \cdot
\beta$ as follows.  We choose representatives of these 2-morphisms,
which we also call $\alpha$ and $\beta$, such that
\[ \alpha \intersect I_1 \times I_2 \times I_3 \times \{1\}=
\beta \intersect I_1 \times I_2 \times I_3\times \{0\},\]
$\alpha$ has a product structure in the $t$
direction for $t \in [1/4,1]$, and $\beta$ has a product structure in
the $t$ direction for $t \in [0,3/4]$.  Then $\alpha\cdot\beta$ is the
2-morphism represented by
\[(\alpha \intersect I_1 \times I_2 \times I_3 \times [0,1/2])\union
(\beta \intersect I_1\times I_2\times I_3\times [1/2,1]) .\]

\begin{figure}[ht]
\centerline{\epsfysize=2in\epsfbox{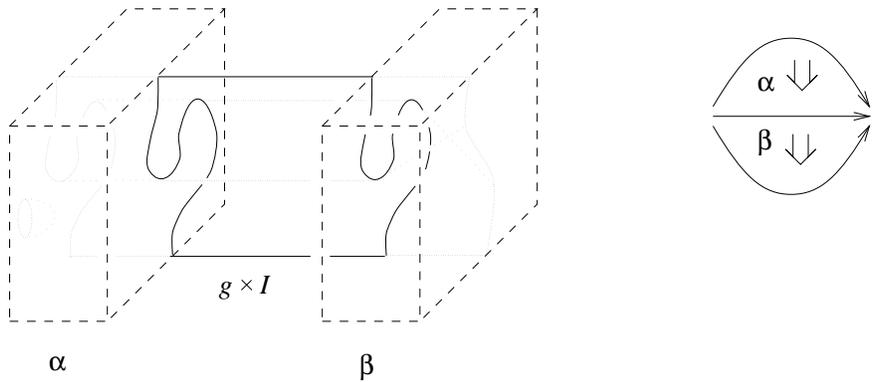}} \medskip
\caption{ Vertical composition of 2-morphisms}
\end{figure}

\begin{lem} \label{composition} \et Composition of 1-morphisms and
horizontal and vertical composition of 2-morphisms are well
defined. \end{lem}

Proof - Suppose we have a composable pair of 1-morphisms.  To check that
their composite is well defined, choose a pair of generic tangles $f,g$
representing these 1-morphisms and satisfying the conditions in the
definition of composition, and also choose another pair, say $f'$ and
$g'$.  Let $H(x,y,z,s)$ be an equivalence isotopy between $f$ and $f'$
with a product structure in the $z$ direction for $z \in [1/4,1]$
--- by which we mean that $H(x,y,z,s) = (X(x,y,s),Y(x,y,s),z)$ for $z\in
[1/4,1]$.  Similarly, let $J(x,y,z,s)$ be an equivalence isotopy between
$g$ and $g'$ with a product structure in the $z$ direction for $z
\in [0,3/4]$.  Glue these together as follows:
\[K(x,y,z,s) = (1-\phi(z)) H(x,y,z,s) +\phi(z) J(x,y,z,s) \]
where $\phi \maps [0,1]\to [0,1]$ is a smooth monotone function that
equals $0$ on $[0,1/4+\epsilon]$ and $1$ on $[3/4-\epsilon,1]$.  We
claim that $K$ is an equivalence isotopy between the tangles
\[(f \cap I_1 \times I_2 \times [0,1/2] )\cup
(g\cap I_1 \times I_2 \times [1/2,1])\]
and
\[(f'\cap I_1 \times I_2 \times [0,1/2] )\cup
(g'\cap I_1 \times I_2 \times [1/2,1]).\]

Clearly $K$ is a homotopy carrying the first tangle to the
second.  The fact that $K$ has a product structure in a neighborhood of
$I_1\times I_2\times \partial(I_3)$ follows from the same properties for
$H$ and $J$.  Similarly, one can show that $K$ is level preserving using
the fact that restricted to any given value of $z$, $K$ is a convex
linear combination of the level-preserving maps $H$ and $J$.

Finally we check that $K$ is really an ambient isotopy.  Note that for
each fixed $s$, the $z$ coordinate of $H(x,y,z,s)$ depends only on $z$,
defining a monotone increasing diffeomorphism of the interval.
Similarly, for each fixed $z$ and $s$, the $y$ coordinate of
$H(x,y,z,s)$ depends only on $y$, defining a monotone increasing
diffeomorphism of the interval.  Also, for each fixed $y,z$ and $s$, the
$x$ coordinate of $H(x,y,z,s)$ depends on $x$, defining a monotone
increasing diffeomorphism of the interval.  The same properties also
hold for $J$.  Since $K$ is a convex linear combination of $H$ and $J$
for each fixed $z$, it follows that $K$ has these properties as well, so
$K(\cdot,s)$ is a diffeomorphism for each $s$.

We can construct similar equivalence isotopies to show that horizontal
and vertical composition of 2-morphisms are well defined.  Since this
is where the height function for tangles is used, we give a
fairly complete sketch of the proof for the more interesting case of
vertical composition.

Suppose we have a vertically composable pair of 2-morphisms.  To show
that their composite is well defined, choose a pair of generic 2-tangles
$\alpha,\beta$ representing them and satisfying the conditions in the
definition of composition, and also choose another such pair, $\alpha'$
and $\beta'$.  In particular, $\alpha$ and $\alpha'$ must have a product
structure in the $t$ direction for $t\in [1/4, 1]$, while $\beta$ and
$\beta'$ have a product structure in the $t$ direction for $t \in
[0,3/4]$.  Let $H$ be an equivalence isotopy between $\alpha$ and
$\alpha'$ that has a product structure in the $t$ direction for $t\in
[1/4, 1]$, meaning that $H(x,y,z,t,s) = (X(x,y,z,s), Y(x,y,z,s),
Y(x,y,z,s),t)$ for $t$ in this interval.  Similarly, let $J$ be an
equivalence isotopy between $\beta$ and $\beta'$ with a product
structure in the $t$ direction for $t \in [0,3/4]$.  Then we glue $H$
and $J$ to get a homotopy
\[ K(x,y,z,t,s) = (1-\phi(t)) H(x,y,z,s) +\phi(t) J(x,y,z,s) \]
where $\phi \maps [0,1]\to [0,1]$ is a smooth monotone function that
equals $0$ on $[0,1/4 +\epsilon]$ and $1$ on $[3/4-\epsilon,1]$.

We claim that $K$ is an equivalence isotopy between the 2-tangle formed
by gluing $\alpha$ and $\beta$ together and that formed by gluing
$\alpha'$ and $\beta'$ together.  Clearly $K$ is a homotopy carrying the
first 2-tangle to the second.  The conditions for $K$ to be an
equivalence isotopy are trivially satisfied for $t$ outside $[1/4,
3/4]$.  For $t$ in this interval, both $H$ and $J$ have a product
structure in the $t$ direction, so we may define a homotopy $K_t$ to be
the restriction of $K$ to a specific value of $t$ in this interval.
Since $\phi$ depends only $t$, this homotopy $K_t$ is a convex linear
combination of the corresponding isotopies $H_t$ and $J_t$.  Since $H$
and $J$ have a product structure for $t$ in this  interval, and $H_1$
and $J_0$ are equivalences of generic tangles, the isotopies $H_t$ and
$J_t$ are also equivalences of generic tangles, and in particular they
are level preserving.  Being a convex linear combination of
level-preserving ambient isotopies, $K_t$ is itself a level-preserving
ambient isotopy.  It follows that $K$ is an ambient isotopy and that the
restriction of $K$ to $I_{1}\times I_{2}\times \{z_0\} \times I_{4}$ is
level preserving for $z_0 = 0$ and $1$.   Similarly, $K$ has a product
structure near $I_1 \times I_2 \times \partial(I_3 \times I_4)$ because
$H$ and $J$ do.

To show that $\alpha \beta$ is well defined, we define a similar
isotopy, pasting together isotopies of representatives of $\alpha$ and
$\beta$ along $1/4 \le z \le 3/4$.  The product structure and
level-preserving properties of these isotopies imply that the resulting
map is an equivalence isotopy, using an argument similar to the previous
ones.  \qed

\subsubsection{Verifying the conditions}

We conclude by checking that the structures described above make $\T$
into a 2-category.

\begin{lem}\label{2cat} \et $\T$ is a 2-category. \end{lem}

Proof - By Lemma \ref{composition}, composition of 1-morphisms and
horizontal and vertical composition of 2-morphisms are well
defined.  One can easily check that the composites $fg$,
$\alpha\cdot\beta$ and $\alpha\beta$ have the desired source and
target.  In addition, the property that a generic 2-tangle intersected
with $I_1 \times I_2
\times \partial I_3 \times I_4$ consists of finitely many straight lines
of the form of the form $(x,y,z)\times I_4$ implies that
\[\target (\source (\alpha)) = \target (\target (\alpha))\]
and
\[\source (\target (\alpha)) = \source (\source (\alpha))\]
for any 2-morphism $\alpha$.

\begin{figure}[p]
\centerline{\epsfysize=1.5in\epsfbox{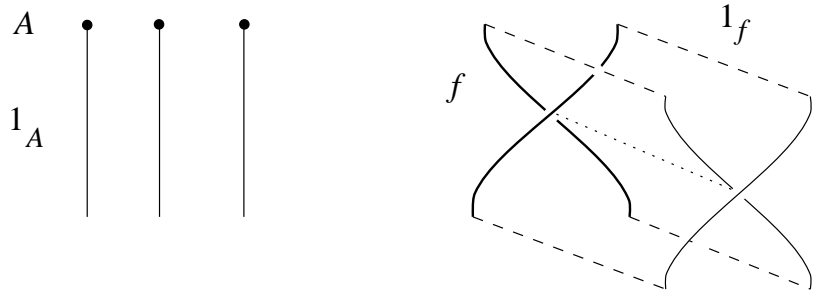}} \medskip
\caption{ Identity 1-morphisms and 2-morphisms}
\end{figure}

\begin{figure}[p]
\centerline{\epsfysize=1.25in\epsfbox{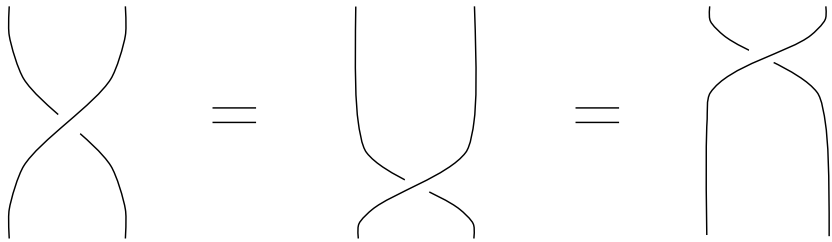}}
\medskip \caption{ $f = 1_A f = f 1_A$}
\end{figure}

\begin{figure}[p]
\centerline{\epsfysize=2in\epsfbox{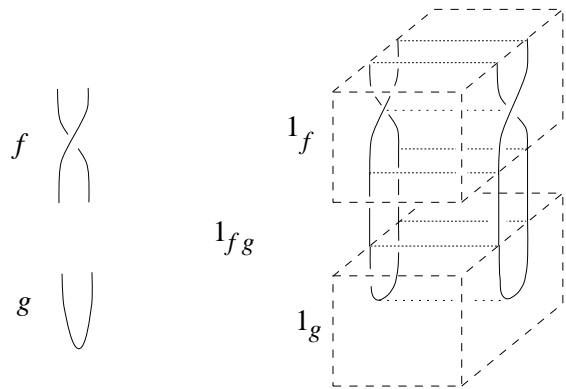}} \medskip
\caption{ $1_{fg} = 1_f 1_g$}
\end{figure}

\begin{figure}[p]
\centerline{\epsfysize=1.5in\epsfbox{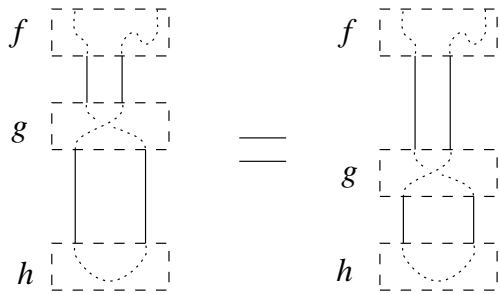}} \medskip
\caption{ $(fg)h = f(gh)$}
\end{figure}

\begin{figure}[p]
\centerline{\epsfysize=2in\epsfbox{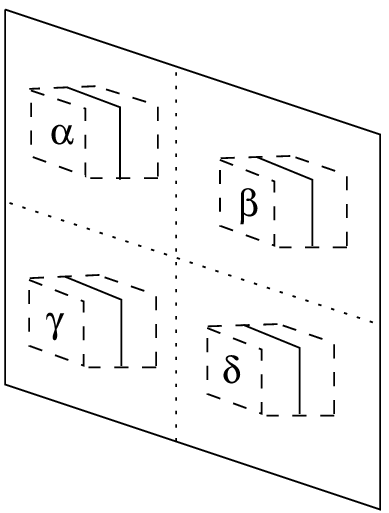}} \medskip
\caption{$(\alpha\cdot\beta)(\gamma\cdot\delta) =
(\alpha\gamma)\cdot(\beta\delta)$ }
\end{figure}

We define identity 1-morphisms and 2-morphisms as shown in Fig.\ 6.  Given
an object represented by the generic finite subset $A$ of the square, we
let the identity of that object be the 1-morphism represented by the
generic tangle $1_A = A \times I_3$.  Similarly, given a morphism
represented by the generic tangle $f$, we let the identity of that
morphism be the 2-morphism represented by the generic 2-tangle $1_f = f
\times I_4$.  Given an equivalence isotopy between two generic finite
subsets $A$ and $A'$ of the square, we can take its product with $I_3$
to get an equivalence isotopy between $1_A$ and $1_{A'}$, and similarly
an equivalence isotopy between generic tangles $f$ and $f'$ gives an
equivalence isotopy between $1_f$ and $1_{f'}$, so the identity of an
object or 1-morphism is well defined.

For any $f\colon A\to B$, the identity 1-morphisms satisfy $1_Af =f1_B =
f$, as shown in Fig.\ 7.  More precisely, notice that we can find
representatives of $f$ that are straight on specified intervals of the
form $[0,\epsilon]$ or
$[\epsilon,1]$ for any $\epsilon \in (0,1)$; a representative that is
straight on $[0,3/4]$ would also clearly be a representative of $1_A f$,
and a representative that is straight on $[1/4,1]$ is, likewise, clearly
a representative of $f1_B$.  Similarly, for any $f,g\colon A\to B$ and
$\alpha \colon f\tto g$, the identity 2-morphisms satisfy $1_f\cdot
\alpha = \alpha \cdot 1_g = \alpha$ and $1_{1_A} \alpha = \alpha
1_{1_B} = \alpha$.

We also have $1_{fg} = 1_f 1_g$ whenever the composite 1-morphism $fg$
is well-defined.  To see this, recall that in the definition of
composition of 1-morphisms, we take representatives of $f$ and $g$ that
have a product structure on certain intervals in the $z$ direction;
likewise, in the definition of horizontal composition of 2-morphisms we
take representatives of $1_f$ and $1_g$ that have a product structure on
the same intervals in the $z$ direction.  If we use representatives of
$f$ and $g$ to obtain representatives of $1_f$ and $1_g$ as above, then
the 2-morphisms $1_{fg}$ and $1_f 1_g$ are represented by the same
surface $(fg) \times I_4$, as shown in Fig.\ 8.

The operations of composition of 1-morphisms and horizontal and vertical
composition of 2-morphisms are associative.  This is true in each case
because the equivalence relations for 1-morphisms and 2-morphisms allow
for smooth changes in the coordinate directions, if small neighborhoods
of the source and target are fixed. This is sufficient to transform
a standard ``follow the definitions'' representative of one composition to
a standard representative of the other, as shown in Fig.\ 9.

Finally we need to check the exchange identity.
Let $\alpha, \beta, \gamma$ and $\delta$ be 2-morphisms such that the
compostions below are defined.  To see that
\[(\alpha\cdot\beta)(\gamma\cdot\delta) =
(\alpha \gamma)\cdot(\beta \delta),\]
consider representatives of the 2-morphism
$(\alpha\cdot\beta)(\gamma\cdot\delta)$ that has a product
structure in the $t$ direction outside small regions representing
$\alpha, \beta, \gamma$ and $\delta$, as in Fig.\ 10.  Splitting this
horizontally, it is clear that this is a representative of
$(\alpha\cdot\beta)(\gamma\cdot\delta)$, but splitting it
vertically, it is also clear that this is a representative of
$(\alpha\gamma)\cdot(\beta\delta)$.  \qed

\subsection{$\T$ is a Monoidal 2-Category}

By a {\it monoidal 2-category}, we mean a semistrict monoidal
2-category as defined by Kapranov and Voevodsky \cite{KV}.  The more
compact formulation later given in HDA1 can be unpacked to give
precisely their definition.  In the following sections we first
introduce some extra structures on the 2-category $\T$, and then prove
they make it into a monoidal 2-category.  These structures are:
\begin{enumerate}
\item An object $I$, called the {\it unit object}.
\item For any objects $A$ and $B$, an object $A \tensor B$.
For any object $A$ and 1-morphism $f\maps B \to C$, 1-morphisms
\[   A \tensor f \maps A \tensor B \to A \tensor C, \qquad
     f \tensor A \maps B \tensor A \to C \tensor A .\]
For any object $A$ and 2-morphism $\alpha \maps f \tto g$, 2-morphisms
\[   A \tensor \alpha \maps A \tensor f \tto A \tensor g , \qquad
      \alpha \tensor A \maps f \tensor A \tto g \tensor A .\]
\item For any 1-morphisms $f\colon A\to A'$ and $g\colon B\to B'$,
a 2-morphism
\[\btensor{}_{f,g}\colon (A\tensor g)(f\tensor B')\tto
(f\tensor B)(A'\tensor g), \]
called the {\it tensorator}.
\end{enumerate}

\subsubsection{The unit object}

We define $I$ to be the object represented by the empty set.

\subsubsection{Tensoring by an object}

Given objects $A,B \in \T$, choose representatives --- which by abuse of
language we also call $A$ and $B$ --- such that the $y$ coordinate of
every point in $A$ is less than the $y$ coordinate of every point in
$B$.  Then we define $A\otimes B \in \T$ to be the object represented
by $A \cup B$.  To simplify notation, we sometimes suppress the
$\tensor$ in writing the tensor product of objects; that is, $A\tensor
B$ may be written as $AB$.

Given an object $A$ and a 1-morphism $f\colon B\to C$ in $\T$, the
1-morphism $A\otimes f\from A\ten B\to A\ten C$ is represented by a
tangle with $1_A$ to the left of $f$.  More precisely, choose
representatives such that the $y$ coordinates of every point in $A$ is
less than the $y$ coordinate of every point in $f$.  Then we define
$A\otimes f$ to be the 1-morphism represented by $(A \times I_3)\cup
f$.  Clearly the source and target of $A \tensor f$ are $A \tensor B$
and $A \tensor C$, respectively.  The product $f\otimes A\from B\ten A
\to C\ten A$ is defined similarly using a tangle with $1_A$ to the right
of $f$.

Given an object $A$ and a 2-morphism $\alpha\from f\tto g$ in $\T$,
the 2-morphism $A\otimes \alpha\from A\ten f\to A\ten g$ is represented
by a 2-tangle with $1_{1_A}$ to the left of $\alpha$.
More precisely, choose representatives such that the $y$ coordinate
of every point of $A$ is less than the $y$ coordinate of every point
of $\alpha$.  Then we define $A\otimes \alpha$ to be the 2-morphism
represented by $(A \times I_3\times I_4)\cup \alpha$.   Clearly
the source and target of $A \tensor \alpha$ are $A \tensor f$ and
$A \tensor g$, respectively.  The product $\alpha\ten A \maps
f \tensor A \tto g \tensor B$ is defined similarly using a 2-tangle
with $1_{1_{A}}$ to the right of $\alpha$.

\begin{figure}[ht]
\centerline{\epsfysize=1.5in\epsfbox{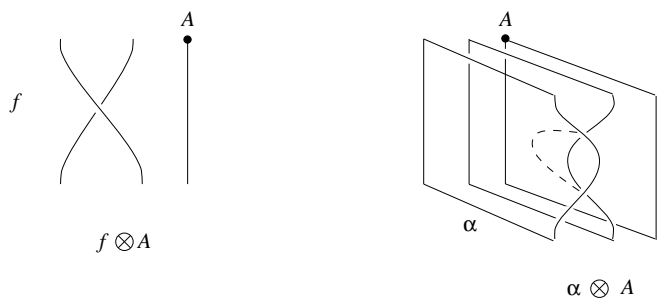}} \medskip
\caption{ Tensor product of an object $A$ with a 1-morphism $f$ and a
2-morphism $\alpha$}
\end{figure}

\begin{lem}\et  The tensor product of an object, 1-morphism or 2-morphism
of $\T$ by an object of $\T$ is well defined. \end{lem}

Proof - We outline the proof for the tensor product $A \tensor \alpha$
of an object $A$ and a 2-morphism $\alpha$; the other results follow
from similar arguments.  Let $A,A'$ be representatives of the
object $A$ and let $\alpha,\alpha'$ be representatives of the
2-morphism $\alpha$, such that the $y$ coordinate of every point
in $A$ is less than some number $y_0$, the $y$ coordinate of every
point in $\alpha$ is greater than $y_0$, and similarly for $A'$ and
$\alpha'$ for some number $y_0'$.  Without loss of generality we
assume that $y_0 \le y_0'$.  
Let $J$ be an equivalence isotopy between $\alpha$ and $\alpha'$ such
that $J(x,y,z,t,s) = (x,y,z,t)$ for all $y < y_0$, and note that $J$
is also an equivalence isotopy between $A \times I_3 \times I_4$ and
itself.  Let $H$ be an equivalence isotopy between
$A$ and $A'$ such that $H(x,y,s) = (x,y)$ for all $y>y_0'$, and let
$K(x,y,z,t,s) = (H(x,y,s),z,t)$.  Note that $K$ is an equivalence
isotopy between $A \times I_3 \times I_4$ and $A' \times I_3 \times
I_4$, and also an equivalence isotopy between $\alpha'$ and itself.
Following $J$ by $K$ we get an equivalence isotopy between
$(A \times I_3\times I_4)\union \alpha$ and $(A'\times I_3\times
I_4)\union \alpha'$.  Thus $A\tensor \alpha$ is independent of the
choice of representatives used to define it.  \qed

\subsubsection{The tensorator}

For any 1-morphisms $f\colon A\to A'$ and $g\colon B\to B'$ we define
the tensorator
\[{\btensor}_{f,g}\colon (A\tensor g)(f\tensor B')\tto
(f\tensor B)(A'\tensor g)\]
as follows.  Take a representative of the 1-morphism $(A\tensor
g)(f\tensor B')$ that consists of straight vertical lines outside
small regions containing representatives of $f$ and $g$.  Let $H\colon
I_1\times I_2\times I_3\times [0,1] \to I_1\times I_2\times I_3$ be an
isotopy that slides the region containing $f$ up, and the region
containing $g$ down, until $f$ is above $g$, as in Fig.\ 12.  We may
choose $H$ to be independent of $s$ near $s = 0,1$, and also choose
it so that
\[ S = \{(H(x,y,z,s),s)\colon \;(x,y,z)\in (A\tensor g)(f\tensor B')\} \]
has generic projections $p(S)$ and $\pi(p(S))$.  Then $S$ represents a
2-morphism, which we define to be $\btensor_{f,g}$.

\begin{figure}[ht]
\centerline{\epsfysize=1.5in\epsfbox{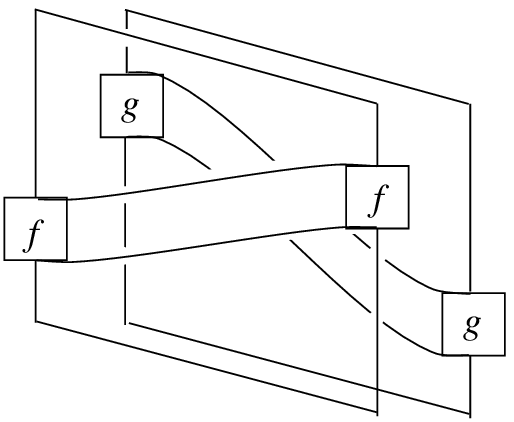}} \medskip
\caption{ The tensorator $\btensor_{f,g}$}
\end{figure}

For a fixed choice of representatives of $f$ and $g$, any two
isotopies of the above type will determine the same 2-morphism
$\btensor_{f,g}$, since they act only by shifting the heights of $f$
and $g$.  Also, any two representatives of $f$ (or $g$) that satisfy
the above conditions will be isotopic by an equivalence isotopy that
is constant outside the regions specified for $f$ and $g$, and this
isotopy can be extended in a natural way to an equivalence isotopy
between the two corresponding 2-tangles.  Hence $\btensor_{f,g}$ is
well defined.  Since $H$ is an isotopy, $\btensor_{f,g}$ is a
2-isomorphism, with ${\btensor}_{f,g}^{-1}$ being the 2-morphism
represented by
\[\{(H(x,y,z,s),1-s)\colon \; (x,y,z)\in (A\tensor
g)(f\tensor B')  \}.\]

\subsubsection{Verifying the conditions}

We conclude by checking that the structures defined above make
$\T$ into a monoidal 2-category.

\begin{lem} \label{m2cat} \et  $\T$ is a monoidal 2-category. \end{lem}

Proof - By Lemma \ref{2cat}, $\T$ is a 2-category.  To prove that the
structures defined above make $\T$ into a monoidal 2-category, we
check that it satisfies the following conditions listed in Lemma 4 of
HDA1.

{\it (i) For any object $A$, the maps $A\tensor -\colon
\T \to \T$ and $-\tensor A\colon \T\to \T$ are 2-functors.}
We only consider the case of tensoring on the left by $A$, as the
tensoring on the right is similar.  First we check that tensoring by
$A$ preserves identities.  Choosing appropriate representatives of $A$
and $B$, the 1-morphism $A \tensor 1_B$ is represented by $(A \times
I_3)\union (B \times I_3)$.  Since this equals $(A \union B)
\times I_3$, which represents $1_{A \tensor B}$, we have
$A \tensor 1_B = 1_{A \tensor B}$.  One can similarly check that
$A\tensor 1_f = 1_{A\tensor f}$ for any 1-morphism $f$.

Next we should check that tensoring with $A$ preserves all three forms
of composition:
\[   A\tensor fg = (A\tensor f)(A\tensor g),  \]
\[   A\tensor (\alpha \beta) = (A\tensor \alpha)(A\tensor \beta) ,\]
\[   A\tensor (\alpha \cdot \beta) = (A\tensor \alpha)\cdot
(A\tensor \beta) \]
Since the arguments for all three cases are similar, we consider only
the third.   Choose representatives of $\alpha$ and $\beta$
that satisfy the conditions in the definition of vertical composition
(the target tangle of $\alpha$ equals the source tangle of
$\beta$, and $\alpha, \beta$ have a product structure in the $t$
direction for $t\in [1/4,1]$ and $t \in [0,3/4]$, respectively)
and lie to the right of a representative of $A$.  Then the
2-tangle
\[ (A\times I_3\times I_4)\union
(\alpha \intersect I_1\times I_2\times I_3\times [0,1/2]) \union
(\beta\intersect I_1\times I_2\times I_3\times [1/2,1])\]
representing $A\tensor (\alpha \cdot \beta)$ and the 2-tangle
\[((A\times I_3\times I_4\union \alpha)\intersect
I_1\times I_2\times I_3\times [0,1/2])\union ((A\times I_3\times
I_4\union\beta)\intersect I_1\times I_2\times I_3\times [1/2,1]))\]
representing $(A\tensor \alpha)\cdot (A\tensor \beta)$ are equal,
as desired.

{\it (ii) For x any object, morphism or 2-morphism,
$x\tensor I = I\tensor x = x$.}  Since $I$ is represented by the empty
set, its product with the intervals $I_3$ and $I_4$ is also empty,
hence $x\tensor I = I\tensor x = x$ for any object, morphism or
2-morphism $x$.

{\it (iii) For $x$ any object, morphism or 2-morphism, and
for any objects $A$ and $B$, we have $A\tensor (B\tensor x) =
(A\tensor B)\tensor x, A\tensor (x\tensor B) = (A\tensor x)\tensor B$
and $x\tensor (A\tensor B) = (x\tensor A)\tensor B$.}  This follows
from the property that equivalence isotopies of objects, 1-morphisms and
2-morphisms all allow shifts in the $y$ direction.

{\it (iv) For any 1-morphisms $f\colon A\to A'$, $g\colon B\to B'$ and
$h\colon C\to C'$, we have ${\btensor}_{A\tensor g,h} =
A\tensor {\btensor}_{g,h}$, ${\btensor}_{f\tensor B,h} =
{\btensor}_{f,B\tensor h}$ and ${\btensor}_{f,g\tensor C} =
{\btensor}_{f,g}\tensor C$.}  The first follows from the fact that we
may choose a representative of ${\btensor}_{A\tensor g,h}$ for which
the tangle representing $1_A$ is straight, and the component
containing it is flat and disjoint from the surface containing the
representative of $g$.  Arguments for the other cases are similar.

{\it (v) For any objects $A$ and $B$ we have $1_A\tensor B =
A\tensor 1_B = 1_{A\tensor B}$, and for any 1-morphisms $f\colon A\to
A'$, $g\colon B\to B'$, we have $\btensor_{1_A,g} =
1_{A\tensor g}$ and $\btensor_{f,1_B} = 1_{f\tensor B}$.}  The
properties for objects are clear from the definitions. The properties
for 1-morphisms follow from the definition of the tensorator:
sliding a straight vertical segment up or down has no effect, so the
only effect of the isotopy defining $\btensor_{1_A,g}$ or
$\btensor_{f,1_B}$ is to slide $g$ or $f$ up or down.  This sliding is
an equivalence isotopy, so the 2-tangle representing
$\btensor_{1_A,g}$ or $\btensor_{f,1_B}$ is equivalent to one with a
product structure in the $t$ direction.  It thus represents an
identity 2-morphism.

{\it (vi) For any 1-morphisms $f\colon A\to A'$, $g,h\colon B\to B'$
and any 2-morphism $\beta\colon g\tto h$,
\[((A\tensor \beta)(1_f\tensor B'))\cdot {\btensor}_{f,h} =
{\btensor}_{f,g}\cdot ((f\tensor B) (A'\tensor \beta)). \]}

{\it (vii) For any 1-morphisms $f,g\colon A\to A'$, $h\colon B\to B'$ and
any 2-morphism $\alpha \colon f\tto g$,
\[((A\tensor h)(\alpha\tensor B'))\cdot {\btensor}_{g,h}
= {\btensor}_{f,h}\cdot ((\alpha\tensor B)(A'\tensor h)) . \]}
The arguments for (vi) and (vii) are similar, so we consider only the
latter.  A picture of the proof is shown in Fig.\ 13 below.  We
slide a piece of the surface that represents $\alpha$ in
the standard representative of $((A\tensor h)(\alpha\tensor B'))\cdot
\btensor_{g,h}$ along the path specified by $\btensor_{g,h}$,
resulting in a surface in which $\alpha$ follows $\btensor_{f,h}$, and
which represents $\btensor_{f,h}\cdot ((\alpha\tensor B)(A'\tensor
h))$.  Clearly, this sliding can be accomplished by an equivalence
isotopy.

\begin{figure}[ht]
\centerline{\epsfysize=1.5in\epsfbox{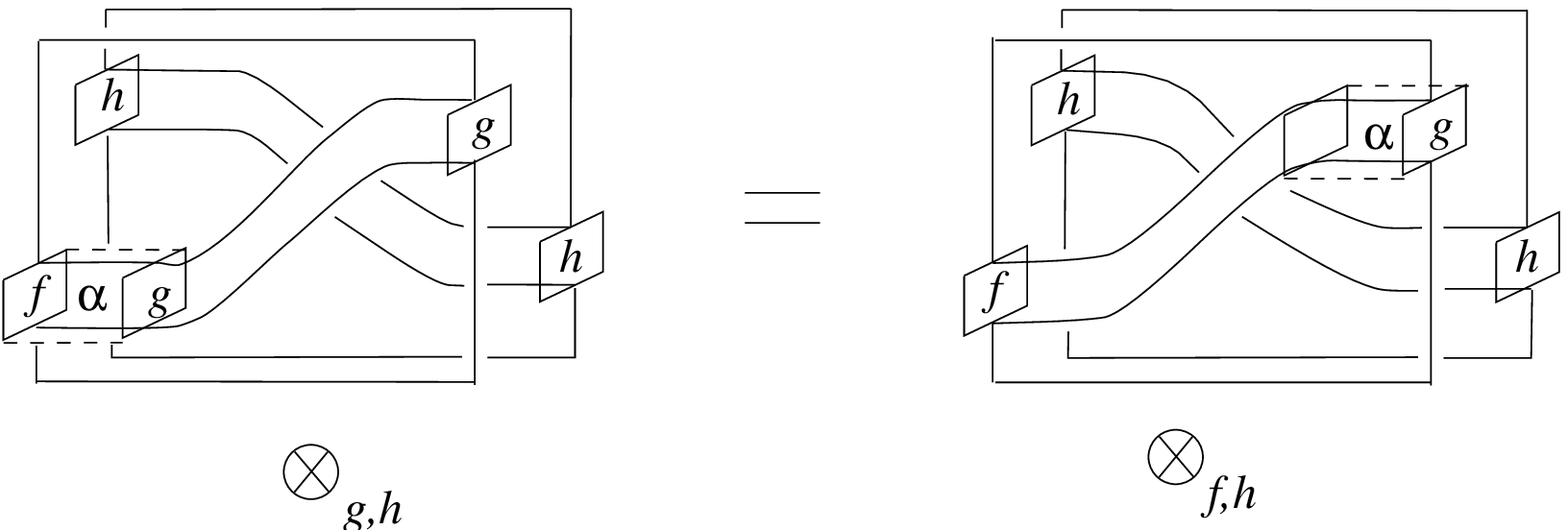}} \medskip
\caption{ $((A\tensor h)(\alpha\tensor B'))\cdot {\btensor}_{g,h} =
{\btensor}_{f,h}\cdot ((\alpha\tensor B)(A'\tensor h))$ }
\end{figure}

\begin{figure}[ht]
\centerline{\epsfxsize=5.5in\epsfbox{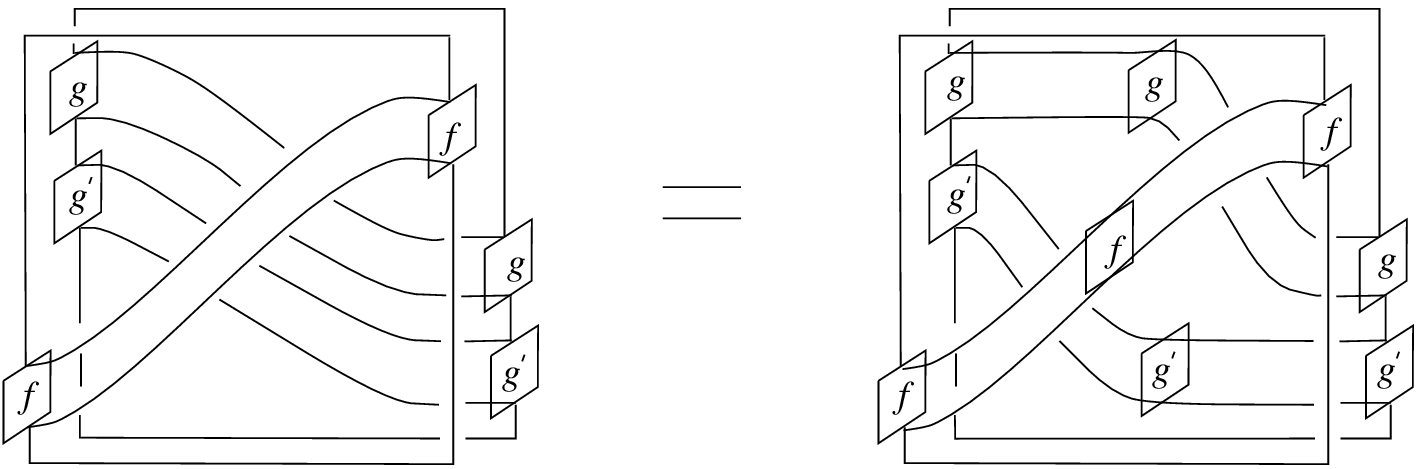}} \medskip
\caption{ $\btensor_{f,gg'} =
((A\tensor g) \btensor_{f,g'})\cdot (\btensor_{f,g} (A'\tensor g'))$}
\end{figure}

{\it (viii) For any 1-morphisms $f\colon A\to A'$, $g\colon B\to B'$
and $g'\colon B'\to B''$,
\[\btensor{}_{f,gg'} =
((A\tensor g) {\btensor}_{f,g'})\cdot ({\btensor}_{f,g} (A'\tensor g')), \]
and for any 1-morphisms $f\colon
A\to A'$, $f'\colon A'\to A''$ and $g\colon B\to B'$,}
\[\btensor{}_{ff',g}= (\btensor{}_{f,g}(f'\ten B'))\cdot ((f\ten B)
\btensor{}_{f',g}).\]
We only the check the first of these conditions, as the second is
similar.  Pictures of $((A\tensor g) \btensor_{f,g'})\cdot
(\btensor_{f,g} (A'\tensor g'))$ and $\btensor_{f,gg'}$ are shown
in Fig.\ 14, for which the representatives of $f$, $g$ and $g'$ are
straight outside small regions.
Clearly there is an equivalence isotopy between these surfaces.  \qed

\subsection{$\T$ is a Braided Monoidal 2-Category}

The definition of `braided monoidal 2-category' has a somewhat complex
history.  The first definition was given by Kapranov and Voevodsky
\cite{KV}.  This definition was modified somewhat in HDA1.
These modifications are necessary for the proper treatment of 2-tangles,
and especially for an unambiguous statement of the Zamolodchikov
tetrahedron equation, as had been noted by Breen \cite{Breen}.  Day and
Street \cite{DS} later gave a more terse formulation of the definition
in HDA1.  Then Crans \cite{C} noted an error in the proof of Theorem 18
of HDA1, and explained how to fix it by adding some conditions
concerning the unit object to the definition of a braided monoidal
2-category.  In what follows, by a `braided monoidal 2-category' we mean
a semistrict braided monoidal 2-category as defined by Crans.

In the following sections we first introduce some extra structures on
the monoidal 2-category $\T$, and then prove they make it into a
braided monoidal 2-category.  These structures are:
\begin{enumerate}
\item For any objects $A,B$ an equivalence
$R_{A,B}\maps A \otimes B \to B \otimes A$, called the
{\it braiding} of $A$ and $B$.
\item For any 1-morphism $f\maps A \to A'$ and object $B$, a
2-isomorphism
\[   R_{f,B} \maps (f \tensor B)R_{A',B} \tto R_{A,B}(B \tensor f) \]
called the {\it braiding} of $f$ and $B$, and
for any object $A$ and 1-morphism $g \maps B \to B'$, a 2-isomorphism
\[   R_{A,g} \maps (A\ten g)R_{A,B'}\tto R_{A,B}(g\ten A) \]
called the {\it braiding} of $A$ and $g$.
\item For any objects $A,B,$ and $C$, 2-isomorphisms
\[\tilde{R}_{(A|B,C)}\maps
(R_{A,B}\ten C)(B\ten R_{A,C})\tto R_{A,B\ten C}\]
and
\[  \tilde{R}_{(A,B|C)} \maps (A \tensor R_{B,C})(R_{A,C} \tensor B)
\tto R_{A \ten B,C},\]
called {\it braiding coherence 2-morphisms}.
\end{enumerate}

\subsubsection{Braiding for objects} \label{braiding.objects}

Given objects $A,B \in \T$, we define the braiding $R_{A,B}\from A\ten
B \to B\ten A$ to be the 1-morphism represented by a tangle consisting
of only positive crossings, such that each strand beginning at a point
$p$ of $A$ crosses all strands beginning at points of $B$ before any
strand starting at a point of $A$ to the left of $p$ crosses any
strands of $B$, as in Fig.\ 15.  We define the 1-morphism
$R^*_{A,B}$ by the reflection in $z$ of a tangle representing
$R_{A,B}$.  Clearly, there are surfaces (which can be defined in terms
of repeated Reidemeister II moves) that represent 2-isomorphisms
between $R_{A,B} R^*_{A,B}$ and $1_{A\ten B}$, and between
$R^*_{A,B} R_{A,B}$ and $1_{B\ten A}$.  Thus $R_{A,B}$ is an
equivalence.

\subsubsection{Braiding for an object and a 1-morphism}

Given a 1-morphism $f\colon A\to A'$ and an object $B$ in $\T$, we
define the 2-morphism
\[  R_{f,B}\from (f\ten B)R_{A',B}\tto R_{A,B}(B\ten f) \]
as follows.  Consider a representative of the source
$(f\ten B)R_{A',B}$ for which $f$ lies behind (has greater $x$ values
than) the strands beginning at $B$.  Then there is an isotopy that
moves $f$ past the strands of $B$, as in Fig.\ 15.  If $H(x,y,z,t)$ is
such an isotopy, the generic 2-tangle representing $R_{f,B}$ is given by
\[      \{(H(x,y,z,t),t)\colon \;(x,y,z)\in (f\ten B)R_{A',B}\} .\]
Since this 2-tangle is traced out by an isotopy in this manner,
$R_{A,f}$ has an inverse represented by the 2-tangle traced out
by the reverse isotopy:
\[    \{(H(x,y,z,1-t),t)\colon \;(x,y,z)\in (f\ten B)R_{A',B}\} .\]
Given an object $A$ and 1-morphism $g \maps B \to B'$, we define the
2-isomorphism
\[    R_{A,g} \maps (A\ten g)R_{A,B'}\tto R_{A,B}(g\ten A)    \]
in a similar way.

\begin{figure}[ht]
\centerline{\epsfxsize=6in\epsfbox{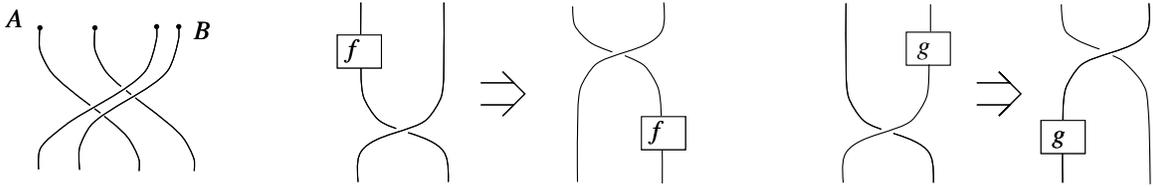}}
\vspace*{-.5in}
\caption{ $R_{A,B}$, $R_{f,B}$ and $R_{A,g}$}
\end{figure}

\subsubsection {Braiding coherence 2-morphisms}

Given objects $A,B,C\in \T$, we define the 2-isomorphism
\[\tilde{R}_{(A|B,C)}\maps
(R_{A,B}\ten C)(B\ten R_{A,C})\tto R_{A,B\ten C}\]
to be represented by a surface built from an isotopy that changes only
the order of distant crossings (crossings with no strands in common)
in a representative of $(R_{A,B}\ten C)(B\ten R_{A,C})$ to give a
representative of $R_{A,B\ten C}$.  Since this 2-morphism is defined
in terms of an isotopy, we can construct an inverse for it using the
isotopy with $t$ reversed, so it is a 2-isomorphism.  Note that if $Z$
is the object represented by a single point, $\tilde R_{(Z|B,C)}$ is
an identity 2-morphism.  We illustrate a nontrivial case of $\tilde
R_{(A|B,C)}$ in Fig.\ 16.

Given objects $A,B,C \in \T$, we define $\tilde{R}_{(A,B|C)}$ to be
the identity.  This makes sense because its source and target are
equal: $(A\ten R_{B,C})(R_{A,C}\ten B) = R_{A\ten B,C}$.

\begin{figure}[h]
\centerline{\epsfysize=1in\epsfbox{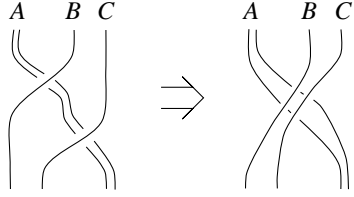}} \medskip
\caption{ $\tilde{R}_{(A|B,C)}$}
\end{figure}

\subsubsection{Verifying the conditions}

We conclude by checking that the structures defined above make
$\T$ into a braided monoidal 2-category.

\begin{lem}\label{bm2cat}\et $\T$ is a braided monoidal 2-category.
\end{lem}

Proof - By Lemma \ref{m2cat}, $\T$ is a monoidal 2-category.  To
prove that the structures defined above make $\T$ into a braided
monoidal 2-category, we verify the conditions listed in Lemma 7
of HDA1, together with the conditions added by Crans.  As in HDA1, we
list some of these conditions using the `hieroglyphic' notation of
Kapranov and Voevodsky.

{\it $(\to \tensor \to)$
For any 1-morphisms $f \maps A \to A'$ and $g \maps B \to B'$,
we have
\[((f\tensor B) R_{A',g})\cdot (R_{f,B} (g\tensor A'))\cdot
(R_{A,B}\; {\btensor}_{g,f})\]
\[ = ({\btensor}_{f,g}^{-1}\; R_{A',B'})\cdot
((A\tensor g) R_{f,B'})\cdot (R_{A,g}(B'\tensor f)) .\] }
This equation corresponds to the equivalence of the 2-tangles shown in
Fig.\ 17.  Since $f$ lies behind $g$ (i.e., the value of the $x$
coordinate on points in the tangle representing $f$ is greater than on
points of that representing $g$), it is clear that the 2-tangles shown
are equivalent.

\begin{figure}[ht]
\centerline{\epsfysize=2.5in\epsfbox{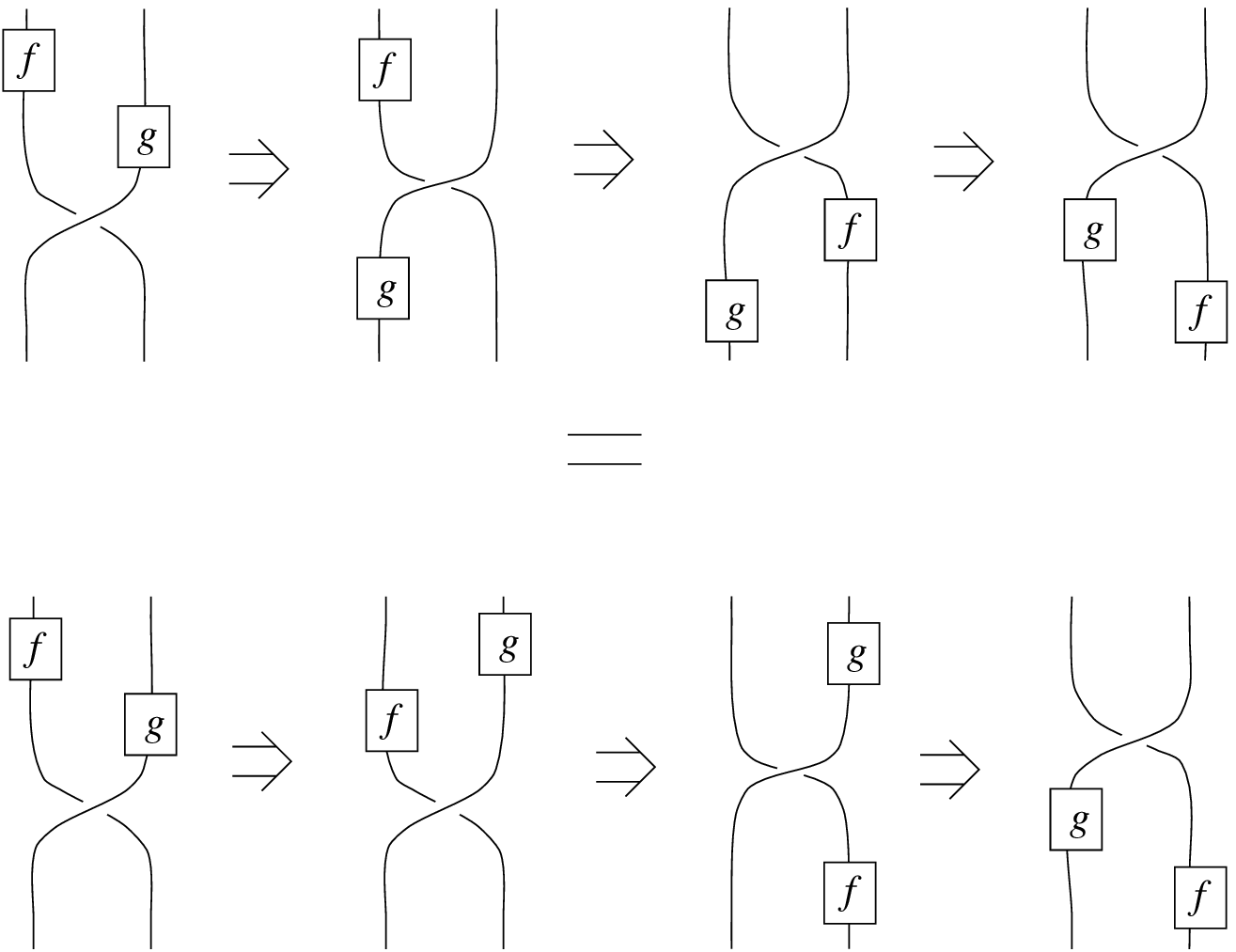}} \medskip
\caption{ $(\to \tensor \to)$}
\end{figure}

{\it $(\bullet \tensor \Downarrow)$
For any 1-morphisms
$f,f' \maps A \to A'$, 2-morphism $\alpha \maps f \tto f'$, and
object $B$, we have
\[    R_{f,B} \cdot (R_{A,B}(B \tensor \alpha)) =
      ((\alpha \tensor B)R_{A',B}) \cdot R_{f',B} . \] }
This equation corresponds to the equivalence of the 2-tangles shown in
Fig.\ 18.

\begin{figure}[ht]
\centerline{\epsfysize=1.5in\epsfbox{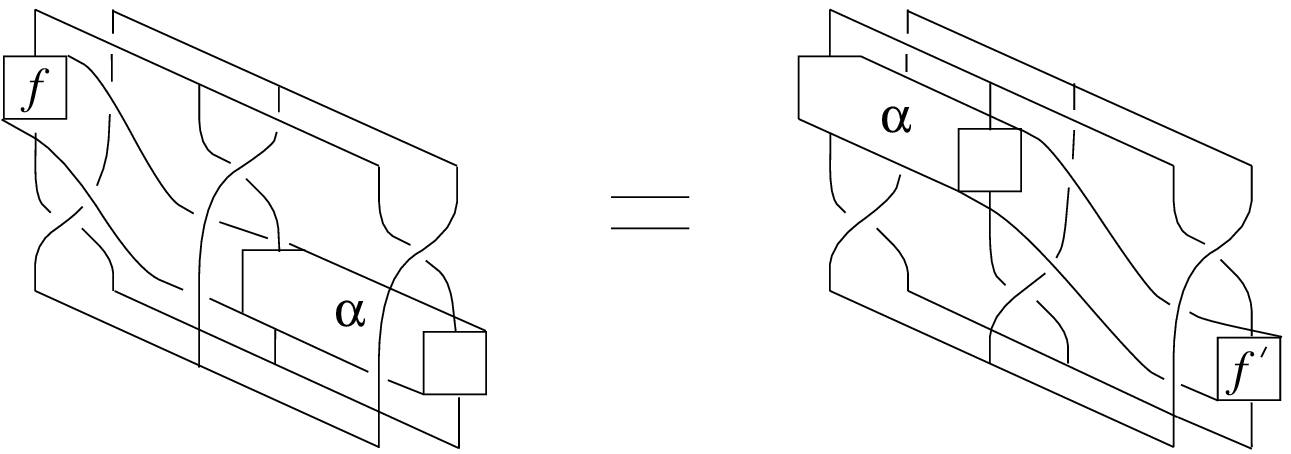}} \medskip
\caption{ $(\bullet \tensor \Downarrow)$}
\end{figure}

$(\Downarrow \tensor \bullet)$ This is similar to $(\bullet \tensor
\Downarrow)$, and follows from an analogous argument.

{\it $(\to\to\ten\bull)$
For any pair of 1-morphisms $f\maps A \to A'$,
$f' \maps A' \to A''$ and object $B$, we have
\[((f\tensor B) R_{f',B})\cdot (R_{f,B}(B\tensor f')) = R_{ff',B}. \] }
This equation corresponds to the equivalence of the 2-tangles
illustrated in Fig.\ 19.  Since the tangles representing $f$ and $f'$
lie behind those representing the braidings $R_{A,B}$, $R_{A',B}$ and
$R_{A'',B}$, it is clear that there is an isotopy between the 2-tangles
shown.

\begin{figure}[ht]
\centerline{\epsfysize=1.5in\epsfbox{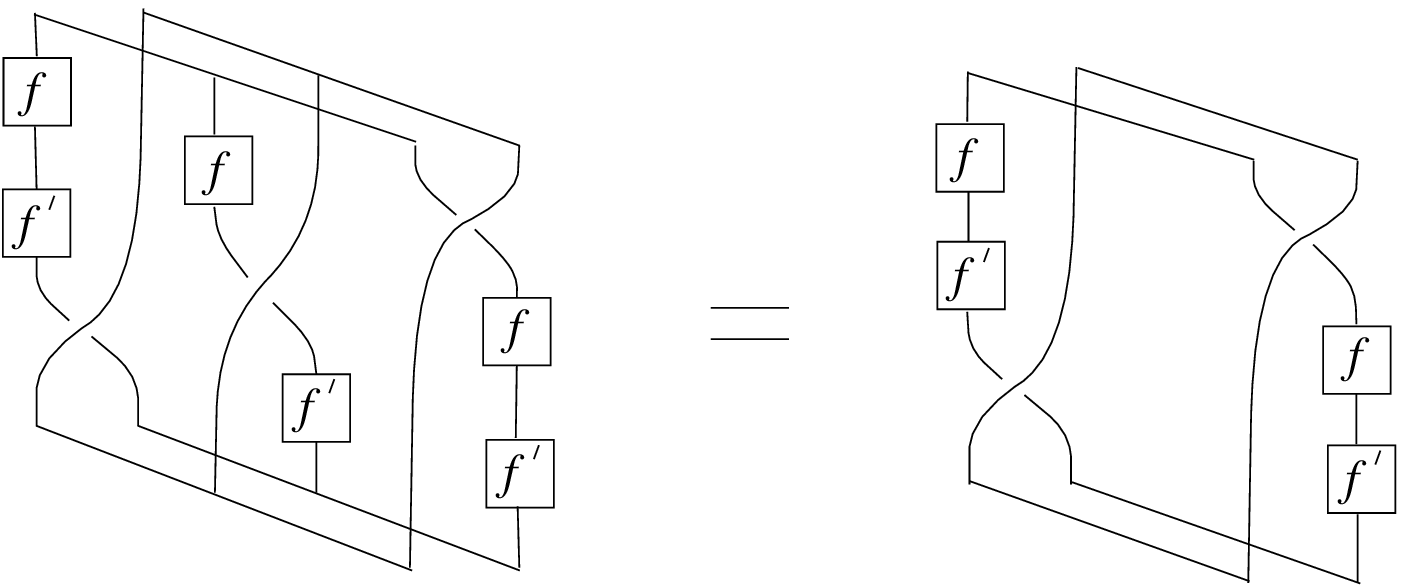}} \medskip
\caption{ $(\to\to\ten\bull)$}
\end{figure}

$(\bull\ten\to\to)$  This is similar to $(\to\to\ten\bull)$.

{\it $(\bullet \tensor (\bullet \tensor \to))$
For any objects $A,B,C$ and 1-morphism $f \maps C \to C'$, we have
\[((AB\tensor f) \tilde{R}_{(A|B,C')})\cdot R_{A,B\tensor f} =\]
\[ (({\btensor}_{R_{A,B},f}(B\tensor R_{A,C'}))\cdot
   ((R_{A,B}\tensor C)(B\tensor R_{A,f}))\cdot
   (\tilde{R}_{(A|B,C)}(B\tensor f\tensor A)). \] }
This equation corresponds to the equivalence of the 2-tangles shown in
Fig.\ 20.

{\it $((\bullet \tensor \bullet) \tensor \to)$
For any objects $A,B,C$ and 1-morphism $f \maps C \to C'$, we have
\[   (A \tensor B \tensor f)\tilde R_{(A,B|C')} \cdot R_{A\tensor B,f} =
((A \tensor R_{B,f})(R_{A,C'} \tensor B)) \cdot
((A \tensor R_{B,C})(R_{A,f} \tensor B)) \cdot
\tilde R_{(A,B|C)}.  \] }
This is similar to $(\bullet \tensor (\bullet \tensor \to))$, but simpler,
because $R_{(A,B|C)}$ and $R_{(A,B|C')}$ are identity 2-morphisms.

$(\to \otimes (\bullet \otimes \bullet))$,
$((\to \tensor \bullet) \tensor \bullet)$,
$((\bullet \otimes \to) \otimes \bullet)$,
$(\bullet \otimes (\to \otimes \bullet))$
These conditions are similar to the relations
$(\bullet \tensor (\bullet \tensor \to))$ and $((\bullet \tensor
\bullet) \tensor \to)$, and 
and are proved analogously.

\begin{figure}[ht]
\centerline{\epsfysize=2.5in\epsfbox{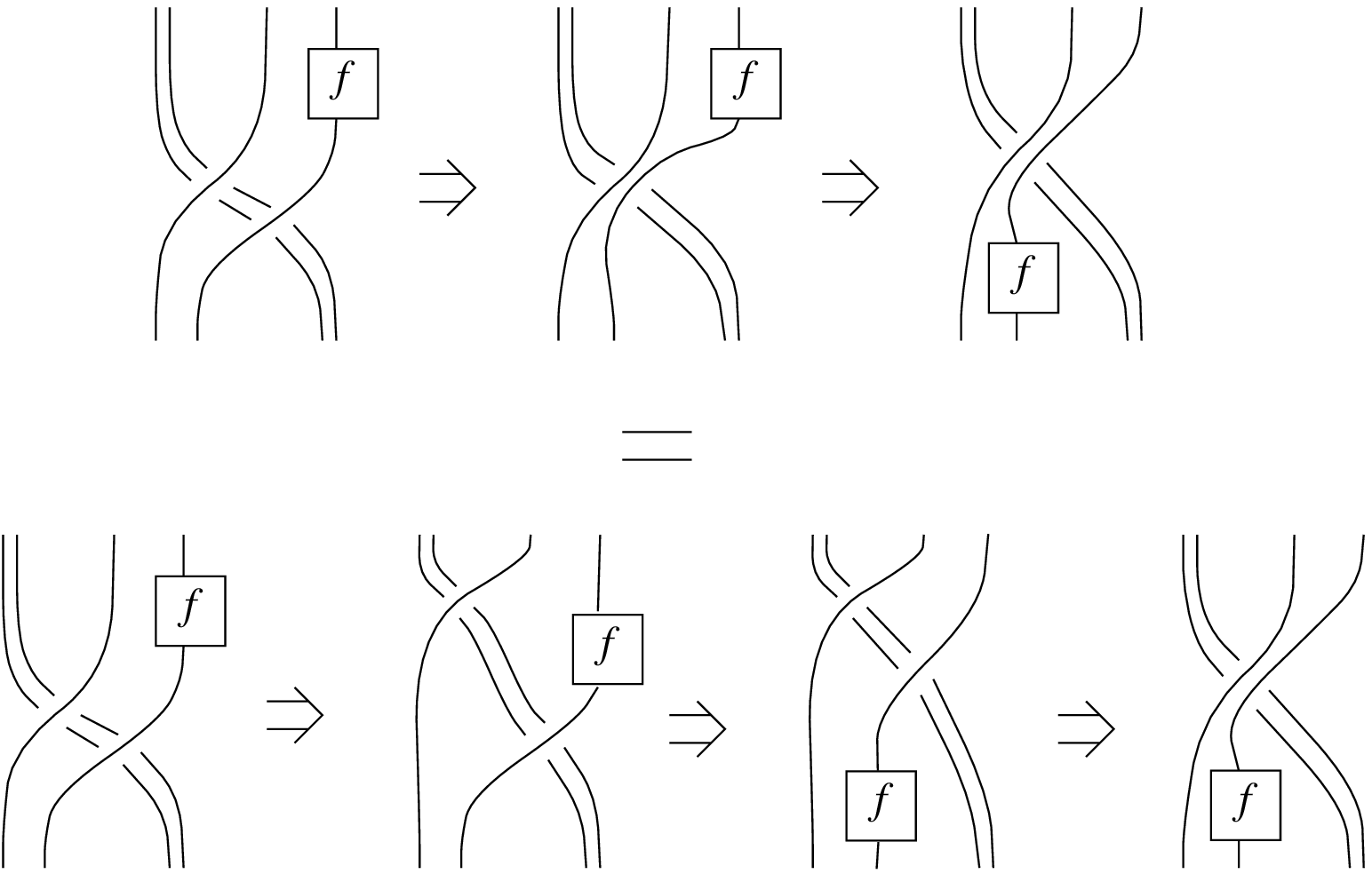}} \medskip
\caption{ $(\bullet \tensor (\bullet \tensor \to))$}
\end{figure}

$((\bullet\tensor\bullet\tensor \bullet)\tensor\bullet)$
This condition holds trivially, since all of the 2-morphisms involved
are of the form $\tilde{R}_{(\cdot,\cdot|\cdot)}$, which are identity
morphisms.

{\it $(\bullet\tensor (\bullet\tensor\bullet\tensor\bullet))$
For any objects $A,B,C,D$, we have
\[((R_{A,B} \tensor C \tensor D)(B\tensor \tilde{R}_{(A|C,D)})\cdot
\tilde{R}_{(A|B,C\tensor D)} =
((\tilde{R}_{(A|B,C)}\tensor D)(B \tensor C \tensor R_{A,D}))\cdot
\tilde{R}_{(A|B\tensor C,D)} . \] }
This equation corresponds to the equivalence of 2-tangles shown in
Fig.\ 21.

\begin{figure}[h]
\centerline{\epsfysize=2.5 in\epsfbox{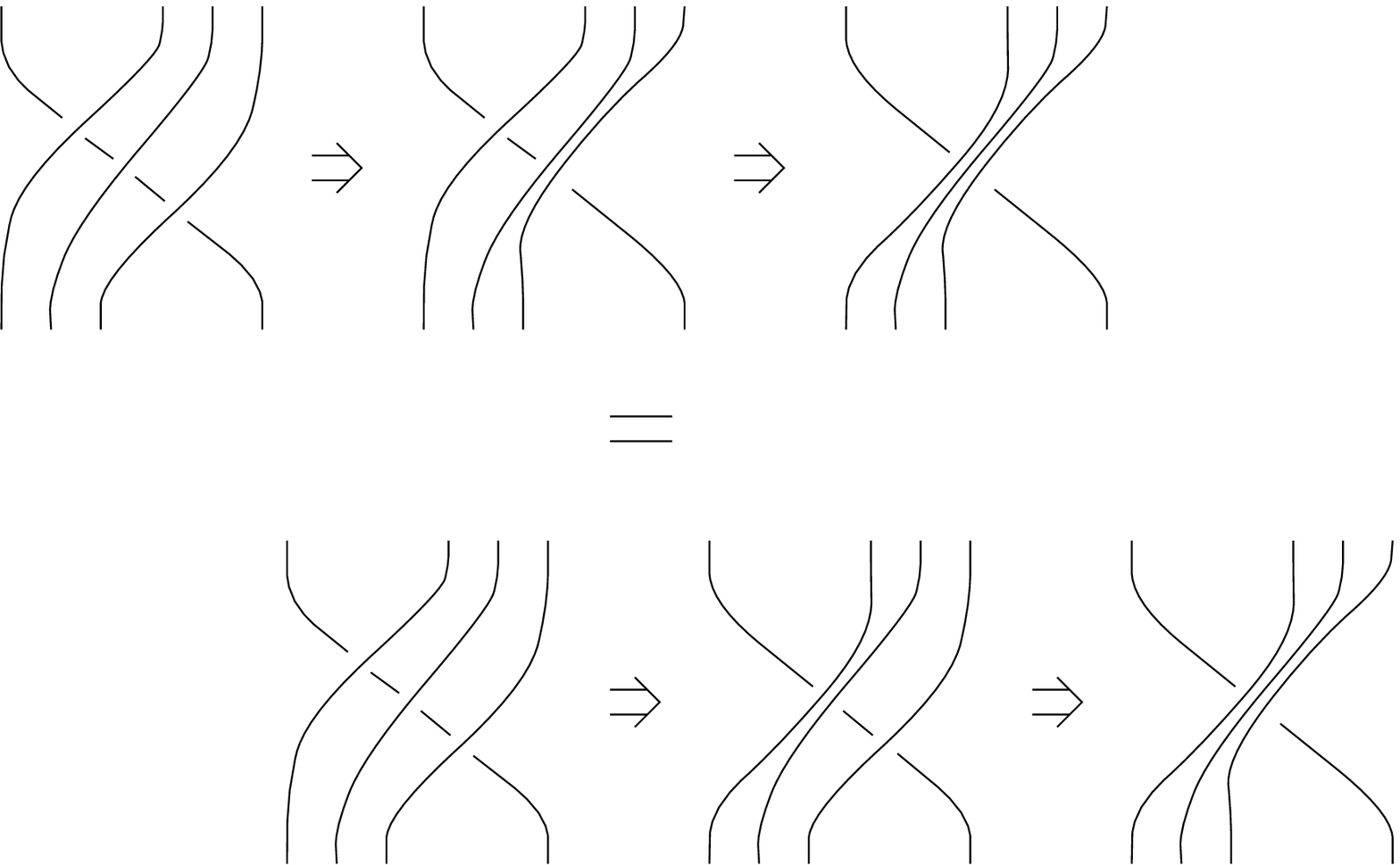}} \medskip
\caption{ $(\bullet\tensor (\bullet\tensor\bullet\tensor\bullet))$ }
\end{figure}

{\it $((\bullet\tensor \bullet )\tensor(\bull\tensor \bullet))$
For any objects $A,B,C,D$, we have
\[((A\tensor R_{B,C}\tensor D)
{\btensor}_{R_{A,C},R_{B,D}} (C\tensor R_{A,D}\tensor B))
\cdot ((\tilde{R}_{(A,B|C)} \tensor D) (C\tensor
\tilde{R}_{(A,B|D)}))\cdot \tilde{R}_{(A \tensor B|C,D)}\]
\[= ((A\tensor \tilde{R}_{(B|C,D)})(\tilde{R}_{(A|C,D)}\tensor B))
\cdot \tilde{R}_{(A,B|C\tensor D)}.\]  }
\noindent Using the fact that 2-morphisms of the form
$\tilde{R}_{(\cdot,\cdot|\cdot)}$ are identities, this equation
simplifies to
\[((A\tensor R_{B,C}\tensor D) {\btensor}_{R_{A,C},R_{B,D}}
(C\tensor R_{A,D}\tensor B)) \cdot \tilde{R}_{(A \tensor B|C,D)}\]
\[=(A\tensor \tilde{R}_{(B|C,D)})(\tilde{R}_{(A|C,D)}\tensor B).\]
When the objects $A, B$ are represented by a single point, the right
side of this equation is an identity 2-morphism because
$\tilde{R}_{(B|C,D)}$ and $\tilde{R}_{(A|C,D)}$ are identities.  The
left side is also an identity 2-morphism, because $\tilde{R}_{(A
\tensor B|C,D)}$ is the inverse of the other factor.  Thus the
equation is true in this case.  The result for other objects $A,B$
follows similarly, and indeed follows inductively from the case we
just considered.

{\it $S^+_{A,B,C} = S^-_{A,B,C}$.
For any objects $A,B,C$, we have
\[   (\tilde R^{-1}_{(A|B,C)} (R_{B,C} \tensor A)) \cdot
     R^{-1}_{A,R_{B,C}} \cdot
     ((A \tensor R_{B,C})\tilde R_{(A|B,C)}) = \]
\[   ((R_{A,B} \tensor C) \tilde R_{(A,B|C)}) \cdot
     R_{R_{A,B},C} \cdot
     (\tilde R^{-1}_{(A,B|C)} (C \tensor R_{A,B})).  \] }
This condition says that the isotopies corresponding to two different
factorizations of the Reidemeister III move give rise to equivalent
2-tangles.   Using the fact that 2-morphisms of the form
$\tilde{R}_{(\cdot,\cdot|\cdot)}$ are identities, it corresponds to the
equivalence of the 2-tangles shown in Fig.\ 22.

\begin{figure}[ht]
\centerline{\epsfysize=1.75in \epsfbox{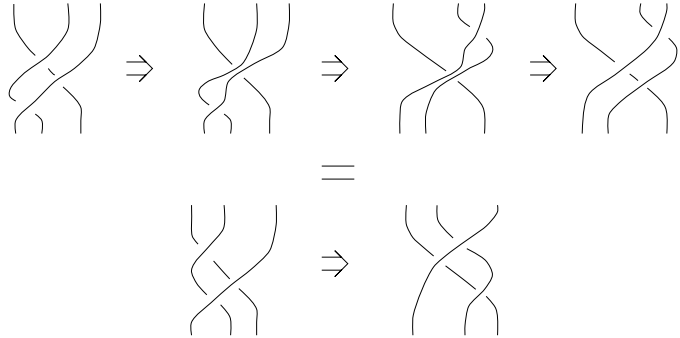}} \medskip
\caption{ $S^+_{A,B,C} = S^-_{A,B,C}$}
\end{figure}

Lastly, we need to check the extra conditions introduced by Crans.
These say that $R_{\cdot,\cdot}$, $\tilde R_{(\cdot|\cdot,\cdot)}$,
and $\tilde R_{(\cdot,\cdot|\cdot)}$ are the identity whenever one of
the arguments is the unit object $I$.  These follow from the fact that
$I$ is represented by the empty set.            \qed

\subsection{$\T$ is a Monoidal 2-Category with Duals}

Now we introduce still more structure on $\T$ and show that this makes
$\T$ into a `monoidal 2-category with duals'.  This is a categorification
of the concept of `monoidal category with duals', which was discussed in
HDA0 and HDA2.  Before giving a precise definition of a monoidal 2-category
with duals, let us sketch the key points, using $\T$ as an example.

There are three levels of duality in $\T$.  First, we can form the
dual of an object in $\T$ by reflecting a set representing it in the
$y$ direction.  Second, we can form the dual of a 1-morphism in $\T$
by reflecting a generic tangle representing it in the $z$ direction.
Third, we can form the dual of a 2-morphism in $\T$ by reflecting a
generic 2-tangle representing it in the $t$ direction.

In addition, the dual of any object $A$ comes equipped with `unit' and
`counit' 1-morphisms:
\[     i_A\maps I\to A\tensor A^*,\qquad  e_A\maps A^*\tensor A\to I , \]
familiar from other contexts, such as monoidal categories with duals.
If $Z$ is the object corresponding to a single point in the unit square,
the unit $i_Z$ corresponds to a tangle with one strand and a single
maximum, while the counit $e_Z$ corresponds to a tangle with one strand
and a single minimum, as shown in Fig.\ 23.   In singularity theory,
the singularity occurring at such a maximum or minimum is called a `fold'.

\begin{figure}[ht]
\centerline{\epsfysize=1.5in\epsfbox{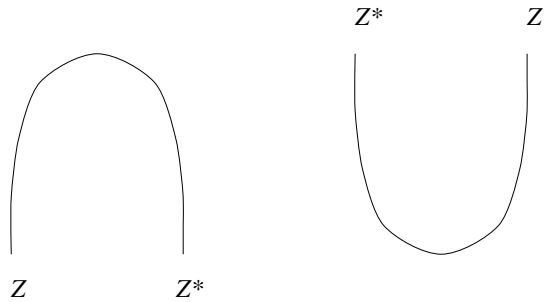}} \medskip
\caption{ Folds corresponding to the unit $i_Z$ and counit $e_Z$}
\end{figure}

In a monoidal category with duals, the unit and counit satisfy
`triangle identities' saying the following diagrams commute:
\[
\begin{diagram} [A \tensor A^\ast \tensor A]
\node{A} \arrow[2]{e,t}{1_A} \arrow{se,b}{i_A \tensor A}
\node[2]{A}   \\
\node[2]{A \tensor A^\ast \tensor A} \arrow{ne,r}{A \tensor e_A}
\end{diagram}
\]
\[
\begin{diagram} [A^\ast \tensor A \tensor A^\ast]
\node{A^\ast} \arrow[2]{e,t}{1_A} \arrow{se,b}{A^\ast \tensor i_A}
\node[2]{A^\ast}   \\
\node[2]{A^\ast \tensor A \tensor A^\ast} \arrow{ne,r}{e_A \tensor A^\ast}
\end{diagram}
\]
In a monoidal 2-category with duals these diagrams no longer commute
on the nose.  Instead, they commute up to specified 2-isomorphisms.
Using duality it turns out to be sufficient to consider only the
first, so we demand the existence of 2-isomorphism
\[   T_A\maps (i_A\tensor A)(A\tensor e_A)\tto 1_A \]
called the `triangulator'.   When $Z$ is the object corresponding
to a single point, $T_Z$ corresponds to the 2-tangle shown in Fig.\ 24.
This 2-tangle describes the process of cancellation of two folds,
a maximum and a minimum.  The singularity occurring at the moment of
cancellation is known as a `cusp'.  In this context, Carter, Rieger
and Saito call it a `cusp on a fold line'.

\begin{figure}[ht]
\centerline{\epsfysize=1.5in\epsfbox{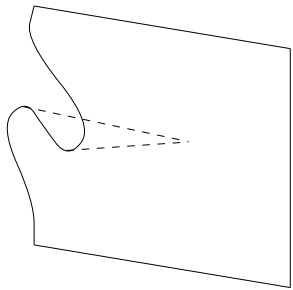}} \medskip
\caption{ Cusp corresponding to the triangulator $T_Z$}
\end{figure}

As usual, when we categorify and replace equations by specified
isomorphisms, the isomorphisms should satisfy new equations of their
own, called coherence laws.  Fig.\ 25 depicts the most interesting
coherence law satisfied by the triangulator in the case $A  = Z$.
The left side of the equation is a 2-tangle with two cusps, while
the right side has no cusps.  The equation arises from an isotopy
between these two 2-tangles corresponding to a cancellation of
cusps.  The singularity occurring when the cusps cancel is called
a `swallowtail'.

\begin{figure}[ht]
\centerline{\epsfysize=1.75in\epsfbox{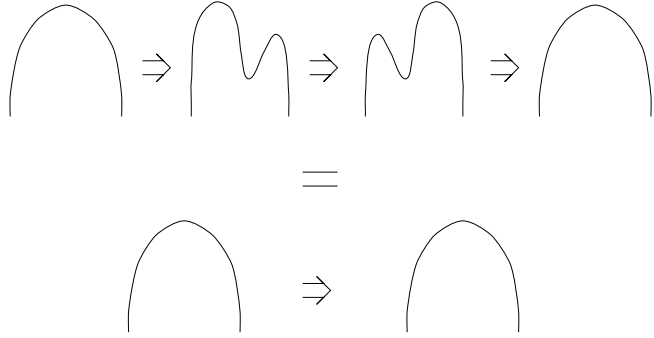}} \medskip
\caption{ Swallowtail coherence law for the triangulator $T_Z$}
\end{figure}

There is a fascinating recursive pattern here.  The fold $e_Z$
describes the process of cancellation of the object $Z$ and its dual
$Z^\ast$.  The cusp $T_Z$ describes the process of cancellation of two
folds.  Similarly, the swallowtail coherence law for $T_Z$ comes from
the process of cancellation of two cusps.  It will be interesting to
see if this pattern continues in higher-dimensional situations.

There is one more piece of structure in a monoidal 2-category with
duals.  Namely, the dual of a 1-morphism $f \maps A \to B$ comes
equipped with `unit' and `counit' 2-morphisms:
\[       i_f\maps 1_A\tto f f^*, \qquad e_f\maps f^* f\tto 1_B .\]
These units and counits satisfy the triangle identities strictly, as
equations.  When $f$ is the braiding $R_{Z,Z}$, the unit $i_f$ and
counit $e_f$ correspond to two forms of the Reidemeister II move, as
shown in Fig.\ 26.   When $f$ is the unit $i_Z$, the unit $i_f$ and
counit $e_f$ correspond to surfaces with a minimum and a saddle point,
respectively, as shown in Fig.\ 27.  The minimum is also called `the
birth of a circle'.  Finally, when $f$ is the counit $e_Z$, the counit
$e_f$ and unit $i_f$ correspond to surfaces with a maximum and a
different sort of saddle point, as shown in Fig.\ 28.  The maximum is
also called `the death of a circle'.  The triangle identities satisfied
by $i_{i_Z}$, $e_{i_Z}$, $i_{e_Z}$ and $e_{e_Z}$ correspond to
cancellation of critical points as in Morse theory, with the $t$
coordinate serving as the Morse function.

\begin{figure}[p]
\centerline{\epsfysize=1.25in\epsfbox{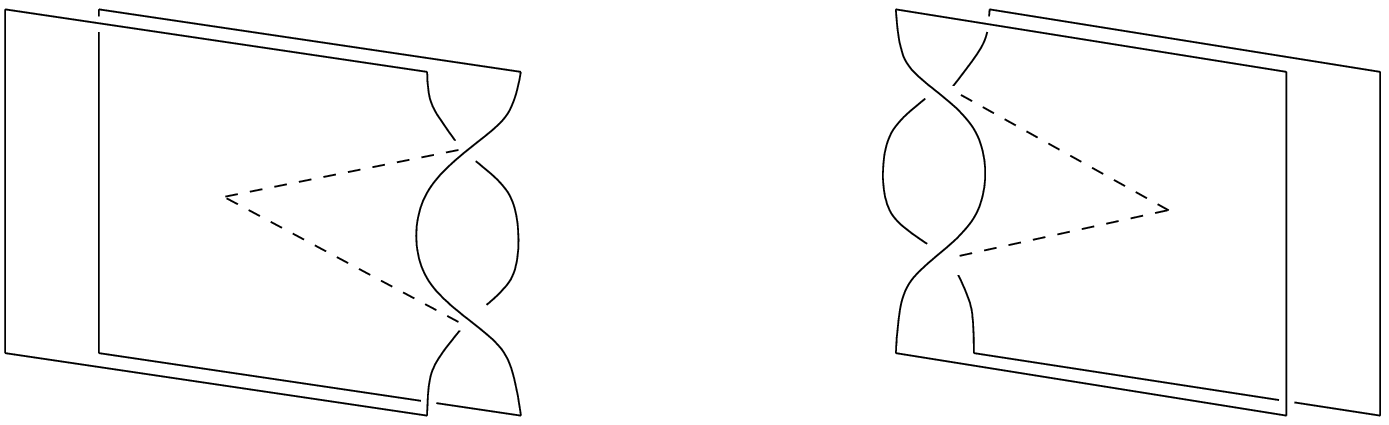}} \medskip
\caption{ 2-Tangles corresponding to $i_{R_{Z,Z}}$ and
$e_{R_{Z,Z}}$}
\end{figure}

\begin{figure}[p]
\centerline{\epsfysize=1.0in\epsfbox{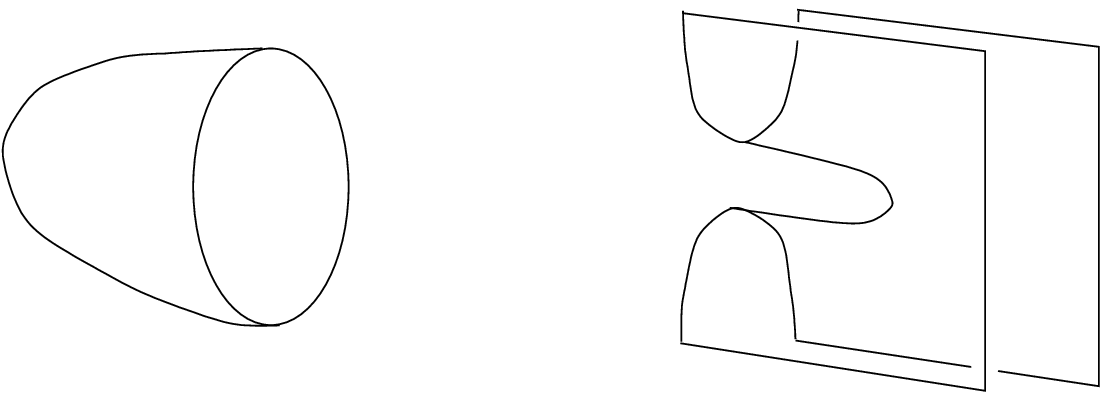}} \medskip
\caption{ 2-Tangles corresponding to $i_{i_Z}$ and $e_{i_Z}$}
\end{figure}

\begin{figure}[p]
\centerline{\epsfysize=1.0in\epsfbox{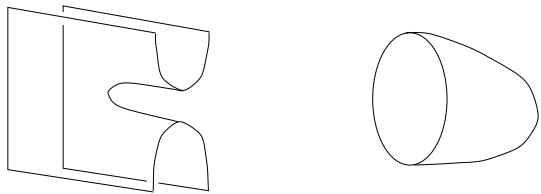}} \medskip
\caption{ 2-Tangles corresponding to $i_{e_Z}$ and $e_{e_Z}$}
\end{figure}

There is no unit or counit associated to the dual of a 2-morphism,
since these would have to be 3-morphisms.  In general, we expect
duality to be `truncated' like this in a monoidal $n$-category with
duals: there should be duals of objects, 1-morphisms, and so on up to
$n$-morphisms, with units and counits for these duals except at the
$n$-morphism level.

Finally, there are various coherence laws that need to be satisfied.
Apart from the swallowtail, these are of four forms.  First, there are
formulas for the duals of the 1-morphisms $i_A,e_A$ and the 2-morphisms
$i_f,e_f,$ and $T_A$. Second, there are compatibility relations between
duality and the various forms of composition: tensoring, composition
of 1-morphisms, and vertical and horizontal composition of 2-morphisms.
Third, the triangulator of the unit object is an identity 2-morphism.
Fourth, there is a compatibility relation between duality for
1-morphisms and 2-morphisms, which is item 12 in the definition below.

The precise definition is as follows:

\vfill \eject
\begin{defn} \label{m2cat.w.duals.defn} \et A {\rm monoidal 2-category
with duals} is a monoidal 2-category equipped with the following structures:

\begin{enumerate}
{\rm

\item For every 2-morphism $\alpha \maps f \tto g$ there is a
2-morphism $\alpha^*\maps g \tto f$ called the {\it dual} of $\alpha$.

\item For every morphism $f \from A\to B$ there is a morphism
$f^*\from B\to A$ called the {\it dual} of $f$, and 2-morphisms
$i_f\maps 1_A\tto f f^*$ and $e_f\maps f^* f\tto 1_B$, called the {\it
unit} and {\it counit} of $f$, respectively.

\item For any object $A$, there is a object $A^*$ called the {\it
dual} of $A$, 1-morphisms $i_A\maps I\to A\tensor A^*$ and $e_A\maps
A^*\tensor A\to I$ called the {\it unit} and {\it counit} of $A$,
respectively, and a 2-morphism $T_A\maps (i_A\tensor A)(A\tensor
e_A)\tto 1_A$ called the {\it triangulator} of $A$.

}
\end{enumerate}

\noindent We say that a 2-morphism $\alpha$ is {\rm unitary} if it is
invertible and $\alpha^{-1} = \alpha^\ast$. Given a 2-morphism $\alpha
\maps f \tto g$, we define the {\rm adjoint} $\alpha^\dagger \maps
g^* \tto f^*$ by
\[ \alpha^\dagger= (g^*i_f) \cdot (g^*\alpha f^*) \cdot (e_g f^*)
.\]

\noindent In addition, the structures above are required to
satisfy the following conditions:

{\rm
\begin{enumerate}

\item  $X^{**} = X$ for any object, morphism or 2-morphism $X$.

\item $1_X^* = 1_X$ for any object or morphism $X$.

\item  For all objects $A,B$, 1-morphisms $f,g$, and 2-morphisms
$\alpha,\beta$ for which both sides of the following equations are
well-defined, we have
\[   (\alpha \cdot \beta)^* = \beta^*\cdot \alpha^*,\]
\[  (\alpha \beta)^* = \alpha^* \beta^*,\]
\[  (fg)^* = g^*f^*,\]
\[     (A\tensor \alpha)^* = A\tensor \alpha^*,  \qquad
(\alpha\tensor A)^* = \alpha^* \tensor A,  \]
\[    (A\tensor f)^* = A\tensor f^*, \qquad
 (f\tensor A)^* =  f^* \tensor A,  \]
and
\[   (A\tensor B)^* = B^*\tensor A^* .\]

\item For all 1-morphisms $f$ and $g$, the 2-morphism $\btensor_{f,g}$ is
unitary.

\item For any object or 1-morphism $X$ we have
$i_{X^*} = e^*_X$ and $e_{X^*} = i^*_X$.

\item For any object $A$, the 2-morphism $T_A$ is unitary.

\item If $I$ is the unit object, $T_I = 1_{1_I}$.

\item For any objects $A$ and $B$ we have
\[   i_{A\tensor B} = i_A (A\tensor i_B \tensor A^*), \]
\[   e_{A\tensor B} = (B^*\tensor e_A \tensor B) e_B, \]
and
\[   T_{A\tensor B} =
[(i_A\tensor A\tensor B)
(A\tensor {\btensor}^{-1}_{i_B,e_A} \tensor B)
(A \tensor B\tensor e_B)]
\cdot [(T_A \tensor B)(A\tensor T_B)] .  \]

\item For any object $A$ and morphism $f$ we have
\[      i_{A \tensor f} = A \tensor i_f , \qquad
     i_{f \tensor A} = i_f \tensor A ,\]
\[      e_{A \tensor f} = A \tensor e_f , \qquad
     e_{f \tensor A} = e_f \tensor A .\]

\item For any 1-morphisms $f$ and $g$, $i_{fg} = i_f\cdot (f i_g f^*)$ and
$e_{fg} = (g^* e_f g) \cdot e_g$.

\item For any 1-morphism $f$, $i_f f\cdot fe_f = 1_f$ and $f^* i_f\cdot
e_f f^* = 1_{f^*}$.

\item For any 2-morphism $\alpha$, $\alpha^{\dagger *} = \alpha^{* \dagger}$.

\item For any object $A$ we have
\[ [i_A(A\tensor T^{\dagger }_{A^*})] \cdot
[{\btensor}_{i_A,i_A}^{-1}(A\tensor e_A\tensor A^*)]\cdot [i_A(T_A
\tensor A^*)] = 1_{i_A}.\]

\end{enumerate}
}
\end{defn}

In what follows we first introduce the structures on $\T$ that make it
into a monoidal 2-category with duals, and then verify that they
satisfy the conditions in the above definition.  In fact, some of these
conditions are redundant.  All the equational laws involving
counits but not units can be derived from those involving units
but not counits by taking duals.  For example, starting from the
first equation in condition 8,
\[  i_{A \tensor B} = i_A (A \tensor i_B \tensor A^*), \]
and taking duals, we obtain
\[  e_{B^* \tensor A^*} = (A \tensor e_{B^*} \tensor A^{*}) e_{A^*}, \]
which, since it holds for all objects $A$ and $B$, implies
the second equation in condition 8:
\[  e_{A \tensor B} = (B^* \tensor e_A \tensor B) e_B .\]
Conversely, of course, the equational laws involving units but not
counits can be derived from those involving counits but not units.
Also, the two equations $i_{X^*} = e_X^*$ and $e_{X^*} = i_X^*$
imply each other.

\subsubsection{Duality for objects}

Given an object $A \in \T$, we define its dual $A^\ast$ to equal $A$.
Given a specific representative of an object in $A \in \T$, we use its
image under the reflection $y \mapsto 1-y$ as a standard
representative of $A^\ast$.

To obtain a representative of the unit $i_A\maps I\to A\tensor A^*$, we
first rotate a representative of $A$ lying in $I_1 \times [0,1/2] \times
\{1\}$ through an angle of 180 degrees around the axis $\{y=1/2, z =
1\}$; the submanifold of the cube traced out by this rotation is a
disjoint union of semicircles with endpoints in the target plane,
$\{z=1\}$, with the right endpoints representing $A^* = A$, and the left
endpoints representing $A$.  We then straighten this submanifold using a
small isotopy so that it has a product structure in the $z$ direction
near $z = 1$; if the isotopy is sufficiently close to the identity,
we obtain a generic tangle, which we take as a representative of $i_A$.
The 1-morphism $i_A$ is independent of our choice of a representative
for $A$, since any two such representatives are equivalent by an
equivalence isotopy that is the identity except for $y < 1/2$, and
such an equivalence isotopy can be extended to the tangle representing
$i_A$.  Moreover, $i_A$ is independent of the isotopy used for
straightening near $z = 1$, provided this isotopy is sufficiently
close to the identity.

Similarly, we obtain a representative of the counit $e_A\maps
A^*\tensor A\to I$ by rotating a representative of $A$ lying in
$I_1 \times [1/2,1] \times \{0\}$ through 180 degrees around the axis
$\{y=1/2, z = 0\}$.  The submanifold traced out by this rotation is a
disjoint union of semicircles with endpoints in the source plane,
$\{z=0\}$, with the right endpoints representing $A$, and the left
endpoints representing $A^*$.  Straightening it near $z = 0$, we
obtain a generic tangle which we take as a representative of $e_A$.

To define the triangulator $T_A$ we proceed inductively, using the
fact that any object $A \in \T$ is a tensor product of copies of $Z$,
where $Z$ is the object represented by the one-point set.  We define
$T_I$ to be $1_I$.  We define $T_Z$ by first choosing a representative
of $Z$, and then choosing an isotopy between the standard
representative of $1_Z$ and the standard representative of
$(i_Z\tensor Z)(Z\tensor e_Z)$.  We require that this isotopy have the
property that the surface $S$ traced out by this isotopy is a generic
2-tangle, and then let $T_Z$ be the 2-morphism represented by $S$.
One can check that $T_Z$ is independent of the choices made.  Since
any other object $A \in \T$ is of the form $Z\ten A'$ for some object
$A'$, we define $T_A$ inductively by the relation in item 7 of
Definition \ref{m2cat.w.duals.defn}:
\[ T_A = [(i_Z\tensor Z\tensor A') (Z\tensor
{\btensor}^{-1}_{i_{A'},e_Z} \tensor A')(Z \tensor A'\tensor e_{A'})]
\cdot [(T_Z \tensor A') (Z\tensor T_{A'})] . \]

\subsubsection{Duality for 1-morphisms}

Given a 1-morphism $f \maps A \to B$ in $\T$, we define the dual
1-morphism $f^\ast \maps B \to A$ as follows: given any representative
of $f$, we take its image under the reflection $z \mapsto 1-z$ as a
representative for $f^\ast$.  Any equivalence between representatives
of $f$ can be similarly reflected, so $f^*$ is well defined.

To define $i_f$, we first rotate a representative of $f$ lying in $I_1
\times I_2 \times [0,1/2] \times \{1\}$ through an angle of 180 degrees
around the plane $\{z=1/2, t = 1\}$; the submanifold of the 4-cube
traced out by this rotation intersects the target hyperplane $\{t = 1\}$
in a generic tangle representing $ff^\ast$.  We then straighten this
submanifold using a small isotopy to make it have a product structure in
the $t$ direction near $t=1$; if the isotopy is sufficiently small,
we obtain a generic 2-tangle, which we take as a representative of
$i_f\maps 1_A \tto ff^\ast$.  As in the previous section, one can
show that this 2-morphism $i_f$ is independent of the choices made.

We define $e_f$ in a similar way by rotating an representative of
$f^\ast$ lying in $I_1 \times I_2 \times [1/2,1] \times \{0\}$ around
$\{z=1/2, t = 0\}$, and then straightening the result.  We obtain
a generic 2-tangle which we take as a representative of $e_f \maps
f^\ast f \tto 1_B$.

\subsubsection{Duality for 2-morphisms}

Given a 2-morphism $\alpha \maps f \tto g$ in $\T$, we define the dual
2-morphism $\alpha^\ast \maps g \tto f$ by taking the image of any
representative of $\alpha$ under the reflection $t \mapsto 1-t$.  Any
equivalence isotopy between representatives of $\alpha$ can be
likewise reflected, so $\alpha^*$ is well defined.

\subsubsection{Verifying the conditions}

We conclude by checking that the structures defined above make
$\T$ into a monoidal 2-category with duals.

\begin{lem}\label{m2cat.w.duals} \et $\T$ is a monoidal 2-category
with duals.  \end{lem}

Proof - We check that $\T$ equipped with the structures given in
the previous sections satisfies the conditions listed in Definition
\ref{m2cat.w.duals.defn}.

1.  Duals of objects, 1-morphisms or 2-morphisms are defined by
reflecting representatives.  Since reflecting twice is an identity
map, $X=X^{**}$ for any object, 1-morphism or 2-morphism $X$.

2.  A standard representative of $1_A$ is a collection of vertical line
segments.  This is unaffected by the reflection $z \mapsto 1-z$, so $1_A
= 1_A^*$.  Similarly, $1_f = 1_f^\ast$, since a standard representative
for $1_f$ is a product of a representative of $f$ with the unit
interval in the $t$ direction, and this 2-tangle is unaffected by the
reflection $t \mapsto 1 - t$.

3.  In a standard representative of the vertical composite $\alpha
\cdot\beta$, $\beta$ follows $\alpha$ in the $t$ direction, so when
we apply the reflection $t \mapsto 1-t$, not only are the $\beta$ and
$\alpha$ components reflected, but also the reflected $\beta$ now
precedes the reflected $\alpha$, hence $(\alpha \cdot \beta)^* =
\beta^*\cdot \alpha^*$.

In a standard representative of the horizontal composite
$\alpha\beta$, $\beta$ is below $\alpha$.  This $z$ ordering is
not changed by the reflection $t \mapsto 1-t$, so $(\alpha
\beta)^* = \alpha^* \beta^*$.

In a standard representative of $fg$, $g$ is below $f$.  Since
$(fg)^\ast$ is obtained by the reflection $z \mapsto 1-z$, the order
of the reflected $f$ and $g$ is reversed, so $(fg)^* = g^*f^*$.

In a standard representative of $A \tensor B$, $B$ is to the right of
$A$.  Since $(A \tensor B)^\ast$ is obtained by the reflection $y
\mapsto 1 - y$, the order of the reflected $A$ and $B$ is reversed, so
$(A \tensor B)^\ast = B^\ast \tensor A^\ast$.

The tensor product of an object $A$ with a 1-morphism (resp.\
2-morphism) $X$ is represented by a disjoint union of $X$ with $1_A$
(resp.\ $1_{1_A}$) on the right or left.  Since the reflections
defining duals for 1-morphisms and 2-morphisms do not change the $y$
coordinate, the order of the tensor product remains the same after
reflection.  Since $1_X^* = 1_X$, we have $(A\tensor X)^* = A
\tensor X^*$ and $(X\ten A)^* = X^*\ten A$.

4.  A standard representative of $\btensor_{f,g}$ is generated by an
isotopy that moves $f$ and $g$ past each other in the $z$ direction.
A representative of $\btensor_{f,g}^{-1}$ is obtained using the same
isotopy with $t$ reversed.  Reversing $t$ in the isotopy amounts to
reflecting the 2-tangle representing $\btensor_{f,g}$ in the $t$
direction, so $\btensor_{f,g}^{-1} = \btensor_{f,g}^*$.

5.  A standard representative of $e_A$ is obtained by rotating a
representative of $A$ lying in the right half of $I_1 \times I_2  \times
\{0\}$ about the axis $\{y = 1/2, z = 0\}$, and then straightening the
resulting submanifold near $z = 0$.  Thus a representative of $e^*_A$
can be obtained by rotating a representative of $A$ lying in the right
half of $I_1 \times I_2 \times \{1\}$ about the axis $\{y = 1/2, z =
1\}$, and straightening the resulting submanifold near $z = 1$.
However, the resulting tangle also serves as a representative of
$i_{A^*}$, since we may also obtain this tangle by rotating a representative 
of $A^*$ lying in the left half of $I_1 \times I_2 \times \{1\}$ about the
axis $\{y = 1/2, z = 1\}$, and straightening the resulting submanifold
near $z = 1$.  Thus we have $i_{A^*} = e_A^*$. Similarly, one can show
that a representative for $i_{f^*}$ is also a representative for
$e_f^*$, so these 2-morphisms are equal.  The conditions $e_{A^*} =
i_A^*$ and $e_{f^*} = i_f^*$ follow by taking duals.

6.  A standard representative of $T_Z$ is obtained from an isotopy
between its source and target tangles, so the 2-morphism $T_Z$ is
unitary (it has an inverse that is given by reversing $t$ in the
isotopy).  The fact that $T_A$ is unitary in general follows from the
fact that $T_Z$ and $\btensor_{\cdot,\cdot}$ are unitary, together with
the relation that defines $T_A$.

7.  We have $T_I = 1_I$ by definition.

8.  A standard representative of $i_{A\tensor B}$ is defined by rotating
a representative of $A \tensor B$ lying in the left half of $I_1 \times
I_2 \times \{1\}$ about the axis $\{y = 1/2, z = 1\}$, and straightening
the resulting submanifold near $z = 1$.  This produces a generic tangle
consisting of nested semicircles, and the maxima of the semicircles
with boundary $A\cup A^*$ are above the maxima of the semicircles with
boundary $B\cup B^*$.  If we take a small isotopy that straightens the
strands for $z$ in a small interval just above the maxima of the
semicircles with boundary $B\cup B^*$, we get an equivalent generic
tangle that represents $i_A \cdot (A\tensor i_B\tensor A^*)$.  It
follows that $i_{A\tensor B} = i_A \cdot (A\tensor i_B \tensor A^*)$.
The condition $e_{A\tensor B} = (B^*\tensor e_A \tensor B) \cdot e_B$
follows from this by taking duals.

The relation
\[T_{A\tensor B} = [(i_A\tensor A\tensor B) (A\tensor
{\btensor}^{-1}_{i_B,e_A} \tensor B)(A \tensor B\tensor e_B)] \cdot
[(T_A \tensor B)(A\tensor T_B)]\]
is a consequence of the relation defining $T_A$, and in fact is
identical to that relation for $A=Z$.

9.   A standard representative of $i_{A \ten f}$ is generated by
rotating a representative of $A \ten f$ lying in $I_1 \times I_2 \times
[0,1/2] \times \{1\}$ around the plane $\{z=1/2, t = 1\}$.   A standard
representative of $A \ten f$ consists of straight strands for the points
of $A$ to the left of a tangle representing $f$.  When rotated, the
strands associated with $A$ become flat planes to the left of a 2-tangle
representing $i_f$, so $i_{A \ten f} = A \ten i_f$.   A similar argument
shows that $i_{f \tensor A} = i_f \tensor A$.  The analogous conditions
for the counits, $e_{A \tensor f} = A \tensor e_f$ and $e_{f \tensor A}
= e_f \tensor A$, follow by taking duals.

10.  A representative of $i_{fg}$ is generated by rotating a
representative of $fg$ lying in $I_1 \times I_2 \times [0,1/2] \times
\{1\}$ around the plane $\{z=1/2, t = 1\}$.   We may choose a
representative of $fg$ so that the tangle representing $f$ lies in the
region where $0 \le z \le 1/4$, while that representing $g$ lies in the
region where $1/4 \le z \le 1/2$.  Then, after straightening the
resulting surface near $t
= 3/4$, the 2-tangle representing $i_{fg}$ is the vertical composite of
a 2-tangle in $I_1 \times I_2 \times I_3 \times [0,3/4]$ representing
$i_{f\tensor f^*}$ and a 2-tangle in $I_1 \times I_2 \times I_3 \times
[3/4,1]$ representing $f i_g f^*$.  We thus have $i_{fg}=i_f \cdot(f
i_g f^*)$.  The condition $e_{fg} = (g^* e_f g) \cdot e_g$ follows by
taking duals.

11.  We specify a representative of $i_f f\cdot fe_f$ as in the
left-hand side of Fig.\ 29.   More precisely, we start by  choosing a
representative of $f$ that consists of straight vertical lines outside
the region where $z < 1/3 - \epsilon$ for some $\epsilon > 0$.  By abuse
of language let us call this representative simply $f$.   Then we take
$(f \cap I_1 \times I_2 \times [0,1/3]) \times \{1/2\}$ and rotate it
180 degrees around  $\{z=1/3, t=1/2 \}$, so that it sweeps out a surface
lying in the region with $t \le 1/2$.  We straighten this to obtain a
surface $S$ that has a product structure in the $t$ direction near $t =
1/2$.   Next, we take $S \cap (I_1 \times I_2 \times [1/3,2/3] \times
\{1/2\})$ and rotate it 180 degrees around $\{z=2/3, t = 1/2\}$ so that
it sweeps out a surface lying in the region with $t \ge 1/2$.  We
straighten this surface to obtain a surface $S'$ with a product
structure in the $t$ direction near $t = 1/2$, for which the union $S
\cup S'$ is a smooth submanifold of the 4-cube.  Finally, we continue $S
\cup S'$ to the boundary of the 4-cube by the product structure in $t$,
obtaining a generic 2-tangle that represents $i_f f \cdot f e_f$.  This
2-tangle is isotopic to a generic 2-tangle representing $1_f$, as shown
in the right-hand side of Fig.\ 29. Hence we conclude that the triangle
identity $i_f f\cdot fe_f = 1_f$  holds.   A similar argument shows that
$fe_f\cdot i_ff= 1_f$.

\begin{figure}[ht]
\centerline{\epsfysize=1.5in\epsfbox{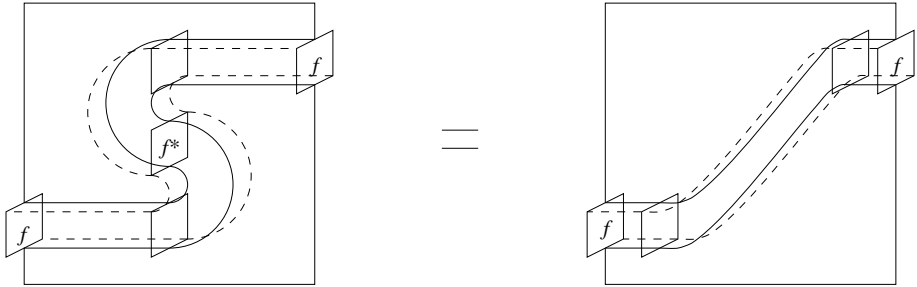}} \medskip
\caption{ The triangle identity for 1-morphisms}
\end{figure}

12.  We can find a representative of $\alpha^{*\dagger}$ as in the
left-hand side of Fig.\ 30.  More precisely, we can find a
representative with a product structure in the $t$ direction outside the
representatives  of $\alpha^*$ in $I_1 \times I_2 \times [1/3,2/3]
\times [1/3,2/3]$ $i_g$ in $I_1 \times I_2 \times [1/3,1] \times
[0,1/3]$, and $e_f$ in $I_1 \times I_2 \times [0,2/3] \times [2/3,1]$.
We may also assume that the surfaces representing $i_g$ and $e_f$ are
formed by rotating the appropriate tangles around $\{z=2/3, t=1/3\}$ and
$\{z=1/3, t=2/3\}$, and then straightening to obtain surface with a
product structure in the $t$ direction near $t = 1/3$ and $t = 2/3$. As
shown in Fig.\ 30, this representative of $\alpha^{*\dagger}$  is
equivalent to a representative of $f^*i_f\cdot e_f f^*\cdot \alpha^*_r$,
where $\alpha^*_r$ is obtained by applying the transformation
$(x,y,z,t) \mapsto (x,y,-z,-t)$ to a representative of $\alpha^*$. Since
$fe_f\cdot i_ff= 1_f$, $i_{f^*} = e^*_f$ and $e_{f^*}= i^*_f$, we have
$f^*i_f\cdot e_f f^*\cdot \alpha^*_r = \alpha^*_r$.  Similarly, we can
find an equivalence isotopy from a representative of $\alpha^{\dagger *}$
to $\alpha^*_r$.  Thus we conclude $\alpha^{\dagger *} = \alpha^{* \dagger}$.

\begin{figure}[ht]
\centerline{\epsfysize=1.5in\epsfbox{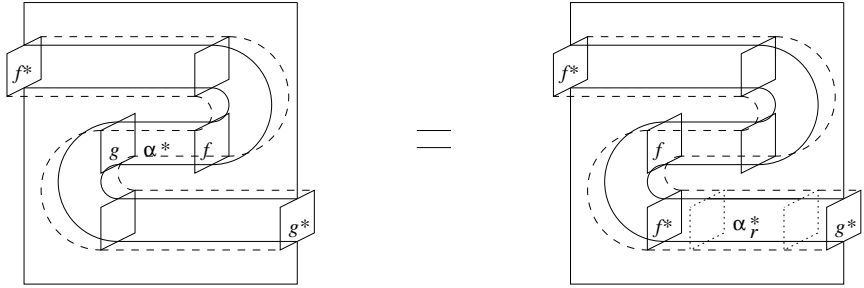}} \medskip
\caption{ Representatives of $\alpha^{* \dagger} $}
\end{figure}

13.  The swallowtail coherence law
\[ [i_A(A\tensor T^{\dagger }_{A^*})] \cdot
[{\btensor}_{i_A,i_A}^{-1}(A\tensor e_A\tensor A^*)]\cdot
[i_A(T_A \tensor A^*)] = 1_{i_A}\]
clearly holds when $A$ is the object $Z$ represented by a one-element
set; see Fig.\ 25.  For $A$ represented by an $n$-element set,
we can inductively straighten sheets beginning either at the outside
or inside, and get the result
\[ [i_A(A\tensor T^{\dagger }_{A^*})] \cdot
[{\btensor}_{i_A,i_A}^{-1}(A\tensor e_A\tensor A^*)] \cdot
[i_A(T_A \tensor A^*)] = 1_{i_A}\]
for any object $A$.
\qed

\subsection{$\T$ is a Braided Monoidal 2-Category with Duals}
\label{bm2cat.w.duals.section}

In general, when an $n$-category with duals is equipped with extra
structure, it is natural to demand that the structural $n$-isomorphisms be
unitary.   For example, in our definition of a monoidal 2-category with
duals, we demanded that the tensorator and triangulator be unitary.  We
incorporate this principle in our definition of a `braided
monoidal 2-category with duals' by requiring that the braiding
coherence 2-morphisms and the braiding for an object and a 1-morphism
2-isomorphisms be unitary.  On the other hand, since the braiding for
a pair of objects is a 1-morphism, we require only that it be unitary
{\it up to specified unitary 2-morphisms}.  In other words, given a
pair of objects $A,B$, we do not demand that $R_{A,B} R_{A,B}^* =
1_{A \tensor B}$ and $R_{A,B}^* R_{A,B} = 1_{B \tensor A}$.  Instead,
we demand only that the 2-morphisms
\[     i_{R_{A,B}} \maps 1_{A \ten B} \tto R_{A,B} R^*_{A,B} \]
and
\[     e_{R_{A,B}} \maps R^*_{A,B} R_{A,B} \tto 1_{B \ten A} \]
be unitary.   For more on this point, see Section \ref{braiding.objects}.

\begin{defn} \label{bm2cat.w.duals.defn} \et  A {\rm braided monoidal
2-category with duals} is a monoidal 2-category with duals that is also
a braided monoidal 2-category for which the braiding is unitary in the
sense that:
\begin{enumerate}
\item For any objects $A,B$, the 2-morphisms $i_{R_{A,B}}$ and
$e_{R_{A,B}}$ are unitary.
\item For any object $A$ and morphism $f$, the 2-morphisms $R_{A,f}$ and
$R_{f,A}$ are unitary.
\item For any objects $A,B,C$, the 2-morphisms $\tilde{R}_{(A,B|C)}$ and
$\tilde{R}_{(A|B,C)}$ are unitary.
\end{enumerate} \end{defn}

In HDA2 we showed that in a braided monoidal category with duals every
object $A$ has an automorphism $b_A \maps A \to A$ called the
`balancing', given by
\[  b_A = (e_A^* \tensor A)(A^* \tensor R_{A,A})(e_A \tensor A) .\]
In the study of tangles, the balancing is closely related to the subtle
issue of framings.   For example, the category of framed oriented
tangles is the free braided monoidal category with duals on an object
$Z$ corresponding to a single positively oriented point in the unit
square.  In this category the balancing $b_Z$ corresponds to a strand
with a single twist in its framing.   The category of unframed oriented
tangles has an extra relation saying that $b_Z = 1$.  In the language
of knot theory, this relation is called the Reidemeister I move.

Similar ideas apply to braided monoidal 2-categories with duals.  We may
define the balancing by the same formula as in a braided monoidal category
with duals.  We expect that the 2-category of framed oriented 2-tangles
is the free braided monoidal 2-category on the object $Z$ corresponding
to a single positively oriented point in the unit square.  However, the
2-category of unframed oriented 2-tangles should be obtained, not by
setting $b_Z$ equal to the identity, but instead by adjoining a unitary
2-morphism $V_Z \maps b_Z \tto 1_Z$ satisfying a new coherence law of
its own.  This 2-morphism corresponds to the {\it process} of performing
the Reidemeister I move.  In other words, it describes the process of
undoing a twist in the framing.

Actually, the connection to the work of Carter, Rieger and Saito
will be clearer if we work not with the balancing $b_Z$ but with the
closely related morphism $i_{Z^*}R_{Z^*,Z} \maps 1 \to Z \tensor Z^*$.
The process of inserting a right-handed twist in the framing then
corresponds to a unitary 2-morphism
\[           W_Z \maps i_Z \tto  i_{Z^*}R_{Z^*,Z} \]
which we call the `writhing', as shown in Fig.\ 31.  For any object
$A$ in a braided monoidal 2-category with duals, one can construct a
unitary 2-morphism $V_A \maps b_A \tto 1_A$ given a unitary 2-morphism
$W_A \maps i_A \tto i_{A^*}R_{A^*,A}$, and conversely.

\begin{figure}[ht]
\centerline{\epsfysize=0.75in\epsfbox{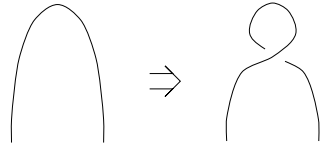}}\medskip
\caption{ The writhing $W_Z$}
\end{figure}

Since we are studying unoriented unframed 2-tangles, where the
generating object $Z$ satisfies $Z = Z^*$, we shall only give the
coherence law for the writhing in the case of a self-dual object.

\begin{defn} \et A {\rm self-dual object}
in a braided monoidal 2-category with duals is an object $A$ with
$A^* = A$. \end{defn}

To state the coherence law for the writhing, it is convenient to
introduce the following 2-morphisms for any pair of objects $X,Y$
in a braided monoidal 2-category with duals:
\[  H_{X,Y} \maps (i_X \ten Y)(X \ten R_{X^*,Y}) \tto
                  (Y \ten i_X)(R^*_{X,Y} \ten X^*)  \]
and
\[  \bar{H}_{X,Y} \maps (i_X \ten Y)(X \ten R^*_{Y,X^*}) \tto
                  (Y \ten i_X)(R_{Y,X} \ten X^*),  \]
defined as follows:
\ban
H_{X,Y} &=& [(i_X \ten Y)(X \ten R_{X^*,Y})(i_{R_{X,Y}} \ten X^*)] \cdot \\
&&          [(i_X \ten Y) \tilde{R}_{(X,X^*|Y)} (R^*_{X,Y} \ten X^*)] \cdot \\
&&          [R_{i_X,Y} (R^*_{X,Y} \ten X^*)]
\ean
\ban
\bar{H}_{X,Y} &=& [(i_X \ten Y)(X \ten R^*_{Y,X*})(i_{R^*_{Y,X}} \ten X^*)]
\cdot \\
&& [(i_X \ten Y) \tilde{R}^{\dagger *}_{(Y|X,X^*)} (R_{Y,X} \ten X^*)] \cdot \\
&&          [R^\dagger_{Y,i^*_X} (R_{Y,X} \ten X^*)]
\ean
When both $X$ and $Y$ are the object $Z \in \T$, Carter, Rieger and
Saito \cite{CRS} call these 2-morphisms `a double point arc crossing
a fold line'.  They can be visualized as in Fig.\ 32.

\begin{figure}[ht]
\centerline{\epsfysize=0.75in\epsfbox{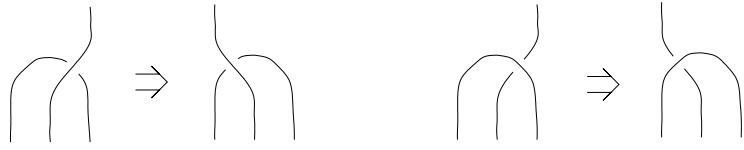}}\medskip
\caption{ $H_{X,Y}$ and $\bar{H}_{X,Y}$}
\end{figure}

\begin{defn}\et A self-dual object $A$ in a braided monoidal 2-category
with duals is {\rm unframed}  if it is equipped with a unitary
2-morphism
\[  W_A \maps i_A \tto i_A R_{A,A} \]
satisfying the equation
\[
T_A^\dagger \cdot
((A\tensor W_A)(e_A \tensor A)) \cdot
((A\tensor i_A)\bar{H}^{\dagger *}_{A,A}) = \]
\[ T_A^{-1} \cdot
((i_A\tensor A)(A\tensor i_{R_{A,A}} e_A))\cdot
((i_A\tensor A)(A\tensor R_{A,A} W_A^\dagger))\cdot
(H_{A,A}(A \tensor e_A))  \]
A unitary 2-morphism satisfying this equation is called a {\rm writhing}
for $A$. \end{defn}

Fig.\ 33 depicts the coherence law satisfied by the writhing
in the case where $A$ is the object $Z \in \T$ corresponding to
a single point in the unit square.

\begin{figure}[h]
\centerline{\epsfysize=2.5in\epsfbox{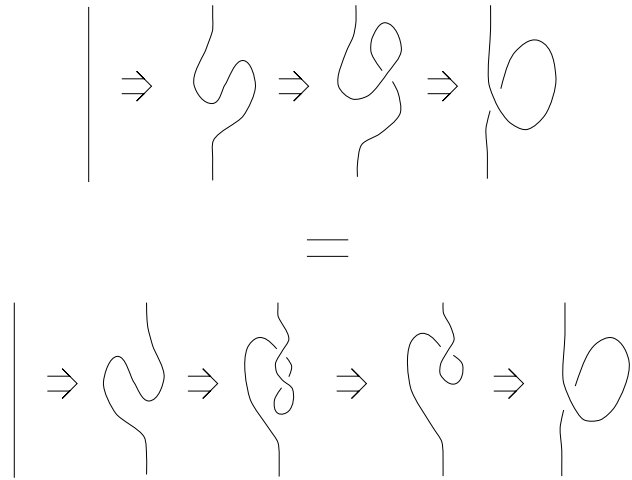}} \medskip
\caption{ The coherence law for the writhing}
\end{figure}

\begin{lem}\et $\T$ is a braided monoidal 2-category with duals,
and the object $Z \in \T$ corresponding to a single point is an
unframed self-dual object.
\end{lem}

Proof - For any objects $A,B \in \T$, one can straighten out a
representative of $i_{R_{A,B}}\cdot i^*_{R_{A,B}}$ using an equivalence
isotopy to obtain the standard representative of $1_{1_{A \ten B}}$.
One can similarly straighten out representatives of $i^*_{R_{A,B}}\cdot
i_{R_{A,B}}$, $e_{R_{A,B}}\cdot e^*_{R_{A,B}}$, and $e^*_{R_{A,B}}\cdot
e_{R_{A,B}}$, showing that $i_{R_{A,B}}$ and $e_{R_{A,B}}$ are unitary.

For any objects $A,B,C \in \T$ and any 1-morphism $f$ in $\T$, the
2-morphisms $R_{A,f}$, $R_{f,A}$, $\tilde{R}_{(A,B|C)}$ and
$\tilde{R}_{(A|B,C)}$ are represented by surfaces traced out by an
isotopy.  It follows that their inverses are defined by the isotopy with
time reversed, so these 2-morphisms are unitary.

Finally, we show that the object $Z$ represented by a single object is
an unframed self-dual object.  We have $Z = Z^*$ by definition, and we
define the writhing $W_Z$ to be the 2-morphism represented by the
generic 2-tangle traced out by the Reidemeister I move.  Since this move
is an isotopy, $W_Z$ is unitary.   Finally, the coherence law
\[
T_Z^\dagger \cdot
((Z\tensor W_Z)(e_Z \tensor Z)) \cdot
((Z\tensor i_Z)\bar{H}^{\dagger *}_{Z,Z}) = \]
\[ T_Z^{-1} \cdot
((i_Z\tensor Z)(Z\tensor i_{R_{Z,Z}} e_Z))\cdot
((i_Z\tensor Z)(Z\tensor R_{Z,Z} W_Z^\dagger))\cdot
(H_{Z,Z}(Z \tensor e_Z))  \]
corresponds the equation shown in Fig.\ 33.  Carter, Rieger and Saito
\cite{CRS} call this `a branch point passing through a cusp' and show
that it holds in $\T$. \qed

The coherence law for the writhing is somewhat mysterious from an
algebraic point of view, but we can offer a partial explanation for it
as follows.   Suppose that $A$ is a self-dual object in a braided
monoidal 2-category with duals, and that $A$ is equipped with a
2-morphism
\[    W_A \maps i_A \tto i_A R_{A,A}. \]
Then we can form a unitary 2-morphism from $i_A$ to $i_A R^*_{A,A}$ in
two different ways, and the coherence law for the writhing says that
these are equal.  When $A$ is the object $Z \in \T$, both these
2-morphisms insert a left-handed twist in the framing.

The first way to form a 2-morphism from $i_A$ to $i_A R^*_{A,A}$ uses
all three levels of duality.  Suppose we have a 2-morphism $\alpha \maps
f \tto g$ between 1-morphisms $f,g \maps X \to Y$.  Then each level of
duality gives us a different way to `reverse' $\alpha$.  As already
discussed, duality at the 2-morphism level gives us a 2-morphism
$\alpha^* \maps g \tto f$, while duality at the 1-morphism level gives
us the 2-morphism $\alpha^\dagger \maps g^* \tto f^*$ defined by
$\alpha^\dagger = (g^* i_f) \cdot (g^* \alpha f^*) \cdot (e_g f^*)$.  In
addition, duality at the object level lets us define 1-morphisms
\[        f^\dagger, g^\dagger \maps Y^* \to X^*  \]
and a 2-morphism
\[      \hat{\alpha} \maps g^\dagger \tto f^\dagger \]
as follows:
\[       f^\dagger = (Y^* \tensor i_X) (Y^* \tensor f \tensor X^*)
(e_Y \tensor X^*),   \]
\[       g^\dagger = (Y^* \tensor i_X) (Y^* \tensor g \tensor X^*)
(e_Y \tensor X^*),   \]
\[       \hat{\alpha} = (Y^* \tensor i_X)(Y^* \tensor \alpha \tensor X^*)
(e_Y \tensor X^*).   \]
A more detailed analysis of the relationships between these operations
can be found in the recent work of Mackaay \cite{M}; this is one place
where there may be room for improvement in our definition of `monoidal
2-category with duals'.

\begin{figure}[h]
\centerline{\epsfysize=1.0in\epsfbox{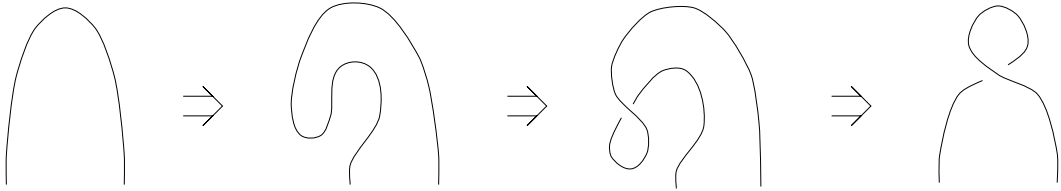}} \medskip
\caption{ The 2-morphism $\bar{W}_Z$ }
\end{figure}

Using these operations, we can form the 2-morphism
\[    \bar{W}_A \maps i_A \tto i_A R^*_{A,A} \]
as the vertical composite of the three 2-morphisms depicted in Fig.\ 34
in the case $A = Z$. The first and third 2-morphisms in this composite
are really just `padding'.   The meat of the sandwich is the second
2-morphism, which is formed by reversing $W_A$ in all three ways listed
above.   More precisely, we have
\[     \bar{W}_A = \gamma \cdot \widehat{W_A^{\dagger *}} \cdot \delta \]
where
\[    \gamma = i_A (T_A^\dagger \ten A)   \]
and
\ban \delta
&=&[i_A (H^*_{A,A} \ten A) (e_A \ten A \ten A)] \cdot
[i_A (i_A \ten A \ten A) (\bar{H}^{\dagger *}_{A,A} \ten A)] \cdot \\
&& [{\bigotimes}_{i_A R^*_{A,A}, i_A} (A \ten e_A \ten A)] \cdot
[i_A R^*_{A,A} (A \ten T_A^{\dagger *})] . \ean

The second way to form a 2-morphism from $i_A$ to $i_A R^*_{A,A}$
is to form the vertical composite of
\[  i_A i_{R_{A,A}} \maps i_A \tto i_A R_{A,A} R^*_{A,A}  \]
and
\[  W^*_A R^*_{A,A} \maps i_A R_{A,A} R^*_{A,A} \tto i_A R^*_{A,A}. \]
In the case $A = Z$, this amounts to doing a Reidemeister II move and
then removing a right-handed twist in the framing, as in Fig.\ 35.

\begin{figure}[h]
\centerline{\epsfysize=1.0in\epsfbox{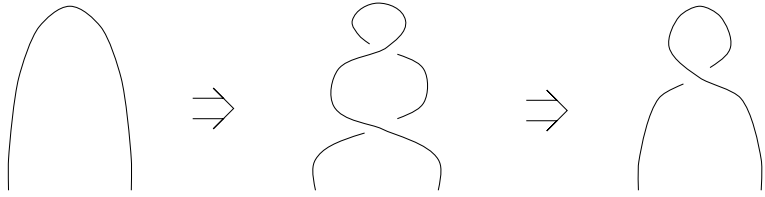}} \medskip
\caption{ The 2-morphism
$(i_Z i_{R_{Z,Z}}) \cdot (W^*_Z R^*_{Z,Z})$
}
\end{figure}

\begin{lem} \label{writhing.lemma} \et Let $A$ be an self-dual object
in a braided monoidal 2-category with duals.  Then a unitary 2-morphism
$W_A \maps i_A \tto i_A R_{A,A}$ is a writhing for $A$ if and only if it
satisfies the equation
\[   \bar{W}_A = (i_A i_{R_{A,A}}) \cdot (W^*_A R^*_{A,A}) .\]
\end{lem}

Proof - We omit the proof, since it is a straightforward although
lengthy calculation.  \qed

\section{A Combinatorial Description of 2-Tangles}

In this section we describe a 2-category $\C$ that is isomorphic to the
2-category $\T$ of unframed unoriented 2-tangles.  Our description is
purely combinatorial, in the sense that we list the objects of $\C$ and
describe its 1-morphisms and 2-morphisms using the method of
generators and relations.   Our list of generators and relations is
based on the work of Carter, Rieger and Saito, so we can use their
results to show that $\C$ and $\T$ are isomorphic \cite{CRS}.   This
isomorphism makes $\C$ into a braided monoidal 2-category with duals.

Our presentation of the 2-category $\C$ makes use of a `bar' operation
on morphisms and 2-morphisms that corresponds in $\T$ to taking the
image of a representative tangle or 2-tangle under the reflection $x
\mapsto 1 - x$.  Thus given a 1-morphism $f \maps A \to B$ we have
$\bar{f} \maps A \to B$, and given a 2-morphism $\alpha \maps f \tto g$
we have $\bar{\alpha} \maps \bar{f} \tto  \bar{g}$.  Some examples
have already been mentioned in Section 2.5.

In a deeper treament, this bar operation would be built into the
definition of `braided monoidal 2-category with duals'.  Since a braided
monoidal 2-category is a special sort of 4-category, a braided monoidal
2-category with duals should really have four duality operations.  In
the study of 2-tangles these operations correspond to reflection along
the $x,y,z$, and $t$ axes.  The last three of these correspond to
duality for objects, morphisms, and 2-morphisms.  The first one remains
obscure  in our approach, but it gives rise to the bar operation.

\subsection{The 2-category $\C$}

In what follows we list the objects of $\C$ and give presentations
for its 1-morphisms and 2-morphisms.  For a more systematic
treatment of 2-category presentations, see the work of Street
\cite{Street}.

\subsubsection{Objects}

Objects in $\C$ correspond to natural numbers, and are denoted by $A_n$,
where $n \in \N$.  These correspond to finite sets of points in the unit
square, where $n$ is the number of points.

\subsubsection{1-Morphisms}

Each 1-morphism in $\C$ represents a planar diagram of a generic tangle.
Using the notation of Carter, Rieger, and Saito, we describe these
1-morphisms using the technique of generators and relations.  The
generating 1-morphisms are
\[         X_{i,j} \maps A_{i+j+2} \to A_{i+j+2},  \]
\[   \bar{X}_{i,j} \maps A_{i+j+2} \to A_{i+j+2},  \]
\[      \cup_{i,j} \maps A_{i+j+2} \to A_{i+j},    \]
\[      \cap_{i,j} \maps A_{i+j}   \to A_{i+j+2},  \]
and the identity 1-morphisms
\[             1_n \maps A_n \to A_n .\]
Here $X$ denotes a right-handed crossing, $\bar{X}$ denotes a
left-handed crossing, $\cup$ denotes a minimum, $\cap$ denotes a
maximum, and the subscripts $i,j \ge 0$ denote the number of vertical
strands to the left and right respectively of this crossing, maximum
or minimum.  The identity 1-morphism $1_n$ corresponds to a tangle
consisting of nothing but $n$ vertical strands.

Every 1-morphism in $\C$ is a composite of these generators, with the
only relations being associativity and the left and right unit laws.
Note that these relations do not include equivalence under Reidemeister
moves, nor changing the height order of crossings and local extrema in
the corresponding tangles.  

We define duals of 1-morphisms in $\C$ as follows.  
For any morphism $f$, there is a 1-morphism $f^*$ defined recursively 
by the following equations: $X_{i,j}^*
= \bar{X}_{i,j}$, $\cap^*_{i,j} = \cup_{i,j}$, $1_n^* = 1_n$, 
$f^{**} = f$, and $(fg)^* = g^*f^*$.  

\subsection{2-Morphisms} \label{2-morphisms}

Our description of the 2-morphisms in $\C$ is compatible with that given
by Carter, Rieger and Saito, with the notation changed to better fit the
2-category structure: we describe surfaces explicitly in terms of
2-morphisms, rather than the sources and targets of these 2-morphisms.
Again, we use the technique of generators and relations.  Having already 
specified the 1-morphisms of $\C$, we now give a list of 2-morphism 
generators going between these 1-morphisms, together with a list of 
relations satisfied by these generators.  Every 2-morphism in $\C$ is 
formed from these generating 2-morphisms by taking vertical and 
horizontal composites.  Two such composites define the same 2-morphism 
in $\C$ if and only if one can get from one to the other using the relations 
in our list, together with the equational laws in the definition of a 
2-category.

The generating 2-morphisms correspond to what Carter, Rieger and Saito
call the `full set of elementary string interactions', together with
identity 2-morphisms for all the 1-morphisms in $\C$.    We list these
generating 2-morphisms below.  The subscripts denote the number of
strands to the left and right that are not affected.

\begin{enumerate}

\item The birth of a circle:

$I_{m,n}\maps 1_{m+n}\tto \cap_{m,n}\cup_{m,n}$.

\item A saddle point:

$E_{m,n}\maps \cup_{m,n}\cap_{m,n}\tto 1_{m+n+2}$.

\item The Reidemeister I move:

$W_{m,n}\maps \cap_{m,n}\tto \cap_{m,n}X_{m,n}$,

\item The Reidemeister II move:

$\two_{m,n}\maps 1_{m+n+2}\tto X_{m,n}\bar{X}_{m,n}$.

\item Three forms of the Reidemeister III move:

$S_{0;m,n}\maps X_{m,n+1}X_{m+1,n}X_{m,n+1}\tto
X_{m+1,n}X_{m,n+1}X_{m+1,n}$,

$S_{1;m,n}\maps X_{m,n+1}X_{m+1,n}\bar{X}_{m,n+1}\tto
\bar{X}_{m+1,n}X_{m,n+1}X_{m+1,n}$,

$S_{2;m,n}\maps X_{m,n+1}\bar{X}_{m+1,n}\bar{X}_{m,n+1}\tto
\bar{X}_{m+1,n}\bar{X}_{m,n+1}X_{m+1,n}$.

\item A double point arc crossing a fold line:

$H_{m,n} \maps \cap_{m,n+1}X_{m+1,n}\tto \cap_{m+1,n}\bar{X}_{m,n+1}$,

\item A cusp on a fold line:

$T_{m,n}\maps \cap_{m,n+1}\cup_{m+1,n}\tto 1_{m+n+1}$,

\item Shifting relative heights of distant crossings and local extrema:

$ N_{Y_{m,n}, Z_{i,j}}\from Z_{i',j}Y_{m,n'}\tto
Y_{m,n}Z_{i,j}$

where $Y$ and $Z$ stand for $X, \bar{X}, \cup$ or $\cap$,  where $i\geq
m+2$ if $Y\neq \cup$, and $i\geq m$ if $Y=\cup$, and $i,j,m,n,i',n'$ are
chosen so that the composite 1-morphisms $Z_{i',j}Y_{m,n'}$ and
$Y_{m,n}Z_{i,j}$ are well defined.

\item Identity 2-morphisms $1_f\maps f\tto f$ for any 1-morphism $f$.

\item For each of the above generating 2-morphisms $\alpha\maps f\tto
g$, a generating 2-morphism $\bar{\alpha}\maps \bar{f} \tto \bar{g}$,
where $\bar{f}$ is obtained from $f$ by replacing each occurence of
$X_{i,j}$ by $\bar{X}_{i,j}$ and vice versa in a product of 
generating 1-morphisms representing $f$.

\item For each of the above generating 2-morphisms $\alpha\maps f\tto g$,
a generating 2-morphism $\alpha^\dagger \maps g^* \tto f^*$.

\item For each of the above generating 2-morphisms $\alpha\maps f\tto
g$, a generating 2-morphism $\alpha^*\maps g\tto f$.

\end{enumerate}

Next we list the relations that these generating 2-morphisms satisfy.
With the help of the relations at the end of our list, relations 1-30
below are equivalent to the correspondingly numbered `movie moves' in
Theorem 3.5.5 of Carter, Rieger and Saito \cite{CRS}.  (We omit their
31st movie move, since it follows from the definition of a 2-category.)
These authors explain how their movie moves arise from singularity
theory, and illustrate many of them with beautiful figures.  To help the
reader find these figures in their paper, we include Carter, Rieger and
Saito's names for the relations below when possible.

\begin{enumerate}

\item
An elliptic confluence of branch points: $W_{m,n}W^*_{m,n} = 1$.

\item
A hyperbolic confluence of branch points: $W^*_{m,n} W_{m,n} = 1$.

\item
An elliptic confluence of double points: $\two_{m,n}\two^*_{m,n} = 1$.

\item
A hyperbolic confluence of double points: $\two^*_{m,n}\two_{m,n} = 1$.

\item Cancelling triple points: $S_{i;m,n}S_{i;m,n}^* = 1$ and
$S_{i;m,n}^* S_{i;m,n} = 1$ for $i = 0,1,2$.

\item A quadruple point in the isotopy (the Zamolodchikov tetrahedron
equation):

$(Z_{i;m,n+1} X_{m+2,n}X_{m+1,n+1}X_{m,n+2})
\cdot (C_{m+1,n+1}B_{m,n+2} Z_{A;m+1,n}X_{m,n+2})$

$\cdot (C_{m+1,n+1} N^*_{B_{m,n+2},X_{m+2,n}}
X_{m+1,n+1}A_{m+2,n}X_{m,n+2})$

$\cdot (C_{m+1,n+1}X_{m+2,n}B_{m,n+2}X_{m+1,n+1}
N_{X_{m,n+2},A_{m+2,n}})$

$\cdot (C_{m+1,n+1}X_{m+2,n}Z_{B;m,n+1} A_{m+2,n})
\cdot (Z_{C;m+1,n}  X_{m,n+2}B_{m+1,n+1}A_{m+2,n})$

$\cdot (X_{m+2,n}X_{m+1,n+1}  N_{X_{m,n+2},C_{m+2,n}}
B_{m+1,n+1}A_{m+2,n})$

$=
(A_{m,n+2}B_{m+1,n+1}  N^*_{C_{m,n+2}X_{m+2,n}} X_{m+1,n+1}X_{m,n+2})$

$\cdot (A_{m,n+2}B_{m+1,n+1}X_{m+2,n}  Z_{C;m,n+1})
\cdot (A_{m,n+2} Z_{B;m+1,n}  X_{m,n+2}C_{m+1,n+1})$

$\cdot (N^*_{A_{m,n+2},X_{m+2,n}}
X_{m+1,n+1}B_{m+2,n}X_{m,n+2}C_{m+1,n+1})$

$\cdot (X_{m+2,n}A_{m,n+2}X_{m+1,n+1} N_{X_{m,n+2},B_{m+2,n}}
C_{m+1,n+1})$

$\cdot (X_{m+2,n}  Z_{A;m,n+1}  B_{m+2,n} C_{m+1,n+1})
\cdot (X_{m+2,n} X_{m+1,n+1} X_{m,n+2} Z_{i;m+1,n})$

where $Z=S$ or $\bar{S}$, $i = 0,1,2$, $A,B,C$ equal
either $X$ or $\bar{X}$ in such a way that
$\source (Z_{i;m,n+1})=A_{m,n+2}B_{m+1,n+1}C_{m,n+2}$, and
$Z_{Y;j,k}$ equals $S_{0,j,k}$ if $Y=X$ or $\bar{S}_{2,j,k}$ if
$Y=\bar{X}$.

\item A branch point moving through a triple point:

$(\cap_{m+1,n}  \two_{m,n+1}  \bar{X}_{m+1,n}X_{m,n+1})
\cdot (\bar{H}^*_{m,n}  \bar{X}_{m,n+1}\bar{X}_{m+1,n}X_{m,n+1})$

$\cdot (\cap_{m,n+1}\bar{X}_{m+1,n}  \bar{S}_{1;m,n})
\cdot (\cap_{m,n+1}  \bar{\two}^*_{m+1,n}
\bar{X}_{m,n+1}\bar{X}_{m+1,n})
\cdot (\bar{W}^*_{m,n+1}  \bar{X}_{m+1,n})$

$ = (\bar{W}^*_{m+1,n}  X_{m,n+1}) \cdot \bar{H}^*_{m,n}$

and

$ (\cap_{m+1,n}  \two_{m,n+1}  X_{m+1,n}X_{m,n+1})
\cdot (\bar{H}^*_{m,n}  \bar{X}_{m,n+1}X_{m+1,n}X_{m,n+1})$

$\cdot (\cap_{m,n+1}\bar{X}_{m+1,n}  \bar{S}_{2;m,n})
\cdot (\cap_{m,n+1}  \bar{\two}^*_{m+1,n}  X_{m,n+1}\bar{X}_{m+1,n})
\cdot (W^*_{m,n+1}  \bar{X}_{m+1,n})$

$= (W^*_{m+1,n}  X_{m,n+1}) \cdot
\bar{H}^*_{m,n}$

\item An elliptic confluence of cusps: $T^*_{m,n} T_{m,n} = 1$.

\item A hyperbolic confluence of cusps: $T_{m,n} T^*_{m,n} = 1$.

\item A swallowtail on the fold lines (the swallowtail coherence
law):

$(\cap_{m,n}  T^\dagger_{m+1,n})
\cdot (N^*_{\cap_{m,n},\cap_{m+2,n}}  \cup_{m+1,n+1})
\cdot (\cap_{m,n} T_{m,n+1})
=1_{\cap_{m,n}}$

\item Removing redundant double points crossing the fold lines:
$H_{m,n} H^*_{m,n}$ and \hfill \break
$H^*_{m,n} H_{m,n}$ are identity 2-morphisms.

\item A branch point passes through a cusp:

$(\cap_{m+1,n}  H^{\dagger}_{m,n})
\cdot (\bar{W}^*_{m+1,n}  \cup_{m,n+1})
\cdot T^{\dagger *}_{m,n}
= (\bar{H}^*_{m,n}  \cup_{m+1,n})
\cdot (\cap_{m,n+1}  W^{\dagger}_{m+1,n})
\cdot T_{m,n}$

and

$(\cap_{m,n+1}  \bar{H}^{\dagger *}_{m,n})
\cdot (\bar{W}^*_{m,n+1}  \cup_{m+1,n})
\cdot T_{m,n}
= (H_{m,n}  \cup_{m,n+1})
\cdot (\cap_{m+1,n}  W^{\dagger}_{m,n+1})
\cdot T^{\dagger *}_{m,n}$

\item A double arc passes over a fold line near a cusp:

$(\cap_{m,n+2} H^\dagger_{m+1,n})
\cdot (N^*_{\cap_{m,n+2},\bar{X}_{m+2,n}} \cup_{m+1,n+1}) 
\cdot (\bar{X}_{m,n} T_{m,n+1})$

$= (H_{m,n+1} \cup_{m+2,n})
\cdot (\cap_{m+1,n+1} N^*_{\bar{X}_{m,n+2}, \cup_{m+2, n}})
\cdot (T_{m+1,n} \bar{X}_{m,n})$
 
\item A triple point near a fold line:

$(\bar{H}_{m+1,n}  A_{i;m,n+2}B_{i;m+1,n+1})
\cdot (\cap_{m+2,n}  Z_{i;m,n+1})
\cdot (N_{B_{i;m,n},\cap_{m+2,n}}  A_{i;m+1,n+1}X_{m,n+2})$

$\cdot (B_{i;m,n}  J_{A;m+1,n+1}  X_{m,n+2})
\cdot (B_{i;m,n}\cap_{m+1,n+1}  N_{X_{m,n+2},\bar{A}_{i;m+2,n}})$

$= (\cap_{m+1,n+1}  N_{A_{i;m,n+2},\bar{X}_{m+2,n}}  B_{i;m+1,n+1})
\cdot (J_{A;m,n+1}  \bar{X}_{m+2,n}B_{i;m+1,n+1})$

$\cdot (\cap_{m,n+2}  \tilde{Z}_{i;m+1,n})\cdot
(N^*_{\cap_{m,n+2},B_{i;m+2,n}}
 \bar{X}_{m+1,n+1}  \bar{A}_{i;m+2,n})
\cdot (B_{i;m,n}   \bar{H}_{m,n+1}  \bar{A}_{i;m+2,n})$

where $i = 0,1,2$, and
$A_{0;j,k}=X_{j,k}$, $B_{0;j,k}=X_{j,k}$, $Z_{0;j,k} = S^*_{0;j,k}$,
$\tilde{Z}_{0;j,k} = \bar{S}_{1;j,k}$;
$A_{1;j,k}=X_{j,k}$, $B_{1;j,k}=\bar{X}_{j,k}$,
$Z_{1;j,k} = \bar{S}^*_{2;j,k}$,  $\tilde{Z}_{1;j,k} = \bar{S}_{0;j,k}$;
and $A_{2;j,k}=\bar{X}_{j,k}$, $B_{2;j,k}=\bar{X}_{j,k}$,
$Z_{2;j,k} = \bar{S}^*_{1;j,k}$, $\tilde{Z}_{2;j,k} = S_{2;j,k}$.
Also, we set $J_{X;j,k} = \bar{H}^*_{j,k}$ and $J_{\bar{X};j,k} = H^*_{j,k}$.

\item $N_{Y_{m,n},Z_{i,j}}$ is unitary for $Y,Z\in
\{X,\bar{X},\cap ,\cup \}$, where $m,n,i,j$ are chosen so that
$N_{Y_{m,n},Z_{i,j}}$ is defined.

\item
$(N^*_{Y_{m,n},Y'_{i,j}} Y''_{k,l})
\cdot (Y'_{i',j}  N^*_{Y_{m,n'},Y''_{k,l}})
\cdot (N^*_{Y'_{i',j},Y''_{k',l}}  Y_{m,n''})$

$= (Y_{m,n}  N^*_{Y'_{i,j},Y''_{k,l}})
\cdot (N^*_{Y_{m,n},Y''_{k'',l}}  Y'_{i,j'})
\cdot (Y''_{k''',l}  N^*_{Y_{m,n'''},Y'_{i,j'}})$

where, here and in the relations below, $Y,Y',Y'' \in \{X,\bar{X},\cap
,\cup \}$ and the subscripts are chosen so that all the above
2-morphisms and composites are defined.

\item $(N^*_{Y_{m,n}, A_{j,k+1}}  B_{j+1,k} C_{j,k+1})
\cdot (A_{j',k+1} N^*_{Y_{m,n}, B_{j+1,k}}  C_{j,k+1})$

$\cdot (A_{j',k+1}B_{j'+1, k} N^*_{Y_{m,n},C_{j,k+1}})
\cdot (Z_{A,B,C;j',k}  Y_{m,n} )$

$= (Y_{m,n}  Z_{A,B,C;j,k})
\cdot (N^*_{Y_{m,n}, C_{j+1,k}}  B_{j,k+1}A_{j+1,k})$

$\cdot (C_{j'+1,k}  N^*_{Y_{m,n},B_{j,k+1}}  A_{j+1,k})
\cdot (C_{j'+1, k}B_{j',k+1} N^*_{Y_{m,n},A_{j+1,k}})$

where $n>k+2$, $A,B,C \in \{X,\bar{X}\}$ satisfy $A_{j,k+1} B_{j+1,k}
C_{j,k+1}= \source (Z_{A,B,C;j,k})$ for $Z_{A,B,C;j,k} = S_{i;j,k}$ or
$\bar{S}_{i;j,k}$, and $j'$ equals $j$ or $j\pm 2$ depending on $Y$.
Also similar relations where $m>j+2$ and $N^*_{Y_{m,n},
\chi_{j,k}}$ is replaced by $N_{\chi_{j,k},Y_{m,n}}$
for $\chi = A,B,C$.


\item $(N^*_{Y_{m,n},\cap_{j,k+1}}  \cup_{j+1,k})
\cdot(\cap_{j',k+1} N^*_{Y_{m,n+2},\cup_{j+1,k}})
\cdot (T_{j',k} Y_{m,n})
= Y_{m,n} T_{j,k}$

where $n>k$ and $j'$ equals $j$ or $j\pm 2$ depending on $Y$.  Also:

$(N_{\cap_{j,k'+1},Y_{m+2,n}} \cup_{j+1,k})
\cdot(\cap_{j,k'+1} N_{\cup_{j+1,k'}, Y_{m,n}})
\cdot(T_{j,k'} Y_{m,n})
= Y_{m,n} T_{j,k}$

where $m>j$ and $k'$ equals $k$ or $k\pm 2$ depending on $Y$.

\item $(N_{Y_{m,n},\cap_{j,k}})
\cdot (Y_{m,n} W_{j,k})
=(W_{j',k} Y_{m,n+2})
\cdot (\cap_{j',k} N_{Y_{m,n+2},X_{j,k}})
\cdot (N_{Y_{m,n},\cap_{j,k}} X_{j,k})$


where $n \geq k$ and $j'$ equals $j$ or $j\pm 2$ depending on $Y$.  Also:

$(N^*_{\cap_{j,k},Y_{m+2,n}})
\cdot (Y_{m,n} W_{j,k'})
=(W_{j,k} Y_{m+2,n})
\cdot (\cap_{j,k} N^*_{X_{j,k},Y_{m+2,n}})
\cdot (N^*_{\cap_{j,k},Y_{m+2,n}} X_{j,k'})$

where $m \geq j$ and $k'$ equals $k$ or $k\pm 2$ depending on $Y$.

\item $(N^*_{Y_{m,n},\cap_{j+1,k}}  \bar{X}_{j,k+1})
\cdot (\cap_{j'+1,k}  N^*_{Y_{m,n+2},\bar{X}_{j,k+1}})
\cdot ({H}^*_{j',k}  Y_{m,n+2})$

$= (Y_{m,n} {H}^*_{j,k})
\cdot (N^*_{Y_{m,n},\cap_{j,k+1}}  X_{j+1,k})
\cdot (\cap_{j',k+1}  N^*_{Y_{m,n+2},X_{j+1,k}})$

where $n>k$, and $j'$ equals $j$ or $j\pm 2$ depending on $Y$.  Also a
similar relation where $m>j+2$ and $N^*_{Y, \chi}$ is
replaced by $N_{\chi,Y}$ for $\chi = \cap, X, \bar{X}$.

\item A double point arc becomes tangent to the plane of projection:

$(\cap_{m+1,n+1} N^*_{X_{m,n+2},\bar{X}_{m+2,n}}\cup_{m+1,n+1})
\cdot (\cap_{m+1,n+1}\bar{X}_{m+2,n}  H^{\dagger}_{m,n+1})$

$\cdot (\bar{H}_{m+1,n} \bar{X}_{m+1,n+1}\cup_{m,n+2})
\cdot (\cap_{m+2,n} \two^*_{m+1,n+1} \cup_{m,n+2})
\cdot (N_{\cup_{m,n},\cap_{m,n}})$

$=(\bar{H}^*_{m,n+1}  \bar{X}_{m+2,n}\cup_{m+1,n+1})
\cdot (\cap_{m,n+2}\bar{X}_{m+1,n+1}  H^{\dagger *}_{m+1,n})$

$\cdot (\cap_{m,n+2} \bar{\two}^*_{m+1,n+1} \cup_{m+2,n})
\cdot (N^*_{\cap_{m,n+2},\cup_{m+2,n}})$

\item $(Y_{m,n}  Z_{i,j})
\cdot (N^*_{Y_{m,n},\chi_{i,j}} \chi^*_{i,j})
= (Z_{i',j}  Y_{m,n})
\cdot (\chi_{i',j}  N_{Y_{m,n'}, \chi^*_{i,j}})$

where $i \geq m$, $Z\in \{\two, I,E^*\}$, $\target(Z) = \chi
\chi^*$, and the subscripts are chosen to be compatible with the
restrictions for defining $N$ and the compositions.   Also a
similar condition where $j \geq n$.


\item A cusp on the double point set:
$(X_{m,n}  \bar{\two}_{m,n})\cdot (\two^*_{m,n}  X_{m,n}) =
1_{X_{m,n}}$

Also, a cusp on the set of fold-lines:
$(I_{m,n}  \cap_{m,n})\cdot (\cap_{m,n}  E_{m,n}) =
1_{\cap_{m,n}}$

\item A horizontal cusp:

$(\cap_{m+1,n}  E^*_{m,n+1})
\cdot (T^{\dagger *}_{m,n}  \cap_{m,n+1})
 = (T^*_{m,n}  \cap_{m+1,n})
\cdot (\cap_{m,n+1}  E_{m+1,n})$

\item A triple point passing through a maximum on the double point set:

$(\two_{m,n+1}  A_{m+1,n}B_{m,n+1})
\cdot (X_{m,n+1}  \bar{S}_{i;m,n})$

$= (A_{m+1,n}B_{m,n+1}  \two_{m+1,n})
\cdot (Z_{i;m,n}  \bar{X}_{m+1,n})$

where $A,B \in \{X, \bar{X}\}$, $i$ satisfies $\source (\bar{S}_{i;m,n})
= \bar{X}_{m,n+1}A_{m+1,n}B_{m,n+1}$, and $Z_{i;m,n} = S^*_{j;m,n}$ or
$\bar{S}^*_{j;m,n}$ satisfies $\source (Z_{i;m,n}) =
A_{m+1,n}B_{m,n+1}X_{m+1,n}$.

\item A maximum point of the double point set being pushed through a
branch point:

$(\cap_{m,n} \two_{m,n})
\cdot (W^*_{m,n}  \bar{X}_{m,n})
= (\bar{W}_{m,n})$

\item A branch point passes over a maximum point of the surface:

$I_{m,n}
\cdot (W_{m,n}  \cup_{m,n})
= I_{m,n}
\cdot (\cap_{m,n}  \bar{W}^{\dagger *}_{m,n})$

\item A double point arc passes over a fold line near a maximum point:

$I_{m,n+1}
\cdot (\cap_{m,n+1}  \two_{m+1,n}   \cup_{m,n+1})
\cdot (H_{m,n}  \bar{X}_{m+1,n}\cup_{m,n+1})$

$= I_{m+1,n}
\cdot (\cap_{m+1,n}  \bar{\two}_{m,n+1}  \cup_{m+1,n})
\cdot (\cap_{m+1,n}\bar{X}_{m,n+1}  H^\dagger_{m,n})$

\item A branch point passes over a saddle point of the surface:

$(\cup_{m,n}  W_{m,n})
\cdot (E_{m,n}  X_{m,n})
=
(\bar{W}^{\dagger *}_{m,n}  \cap_{m,n})
\cdot (X_{m,n}  E_{m,n})$

\item A double point arc passes over a fold line near a saddle point:

$(H^{\dagger *}_{m,n}  \cap_{m+1,n} \bar{X}_{m,n+1})
\cdot (X_{m,n+1}  E_{m+1,n}  \bar{X}_{m,n+1})
\cdot \two^*_{m,n+1}$

$=(\bar{X}_{m+1,n}\cup_{m,n+1}  H^*_{m,n})
\cdot (\bar{X}_{m+1,n}  E_{m,n+1}  X_{m+1,n})
\cdot (\bar{\two}^*_{m+1,n})$

\end{enumerate}

\noindent In addition to each of the above relations $\alpha = \beta$,
we include the analogous relations:

\begin{itemize}

\item  $\bar{\alpha } = \bar{\beta}$, where we impose the relations
$\bar{I}_{m,n} = I_{m,n}$, $\bar{E}_{m,n} = E_{m,n}$, $\bar{T}_{m,n} =
T_{m,n}$, $\bar{N}_{Y_{m,n},Z_{i,j}} = N_{\bar{Y}_{m,n},
\bar{Z}_{i,j}}$, $\bar{1}_f = 1_{\bar{f}}$, and we define $\bar{\alpha}$
for 2-morphisms $\alpha$ other than the generating 2-morphisms listed in
items 1 through 9 using the relations $\bar{\bar{\alpha}} = \alpha$,
$\overline{\alpha\cdot\beta} = \bar{\alpha}\cdot \bar{\beta}$, and
$\overline{\alpha \beta} = \bar{\alpha} \bar{\beta}$,
$\overline{\alpha^\dagger} = \bar{\alpha}^\dagger$, and
$\overline{\alpha^*} = \bar{\alpha}^*$.

\item $\alpha^\dagger = \beta^\dagger$, where we impose the relations
$I^\dagger = I^{*}$, $E^\dagger = E^{*}$, $\two^\dagger = \two^{*}$,
$S^\dagger_{0,m,n} = \bar{S}^{*}_{0,m,n}$, $S^\dagger_{1,m,n} = S^{*}_{2,m,n}$,
$N^\dagger_{f,g} = N_{f^*,g^*}$ and $1_f^\dagger= 1_{f^*}$, and
define $\alpha^\dagger$ for 2-morphisms other than the generating 2-morphisms
listed in items 1 through 10 using the relations
$\alpha^{\dagger \dagger} = \alpha$, $(\alpha\cdot \beta)^\dagger =
\beta^\dagger \cdot \alpha^\dagger$, $(\alpha \beta)^\dagger =
\beta^\dagger \alpha^\dagger$, $\bar{\alpha}^\dagger =
\overline{\alpha^\dagger}$, and $\alpha^{* \dagger} = \alpha^{\dagger *}$.

\item  $\alpha^* = \beta^*$, where we impose the relation $1_f^* = 1_f$,
and define $\alpha^*$ for 2-morphisms other
than the generating 2-morphisms listed in items 1 through 11 using the
relations $\alpha^{**} = \alpha$, $(\alpha \cdot \beta)^* = \beta^*\cdot
\alpha^*$, and $(\alpha \beta)^* = \alpha^* \beta^*$, $\bar{\alpha}^* =
\overline{\alpha^*}$, and $\alpha^{\dagger *} = \alpha^{* \dagger}$.

\end{itemize}

\subsection {$\T$ and $\C$ are Isomorphic}

We now show that $\T$ and $\C$ are isomorphic as 2-categories.  Recall
that two 2-categories $\cal A, \cal B$ are said to be `isomorphic' if
there exist 2-functors $F \from \cal A \to \cal B$ and $G \from
\cal B \to \cal A$ that are inverses in the strictest possible sense:
$FG = 1_{\cal A}$ and $GF = 1_{\cal B}$.

\begin{thm}\et $\T$ and $\C$ are isomorphic 2-categories.
\end{thm}

It suffices to construct a 2-functor $F\from \T \to \C$ and show that it
is bijective on objects, 1-morphisms and 2-morphisms.   Each object $A \in
\T$ is determined by the number of points in a representative of $A$. We
define $F(A)=A_n$, where $n$ is the number of points in a representative
of $A$. Clearly, $F$ is well defined and bijective on objects, since
objects in both $\T$ and $\C$ are uniquely determined by natural
numbers.

Each 1-morphism $f_{\T}$ in $\T$ is represented by a generic tangle.  The
planar diagram of this tangle has crossings and extrema at distinct
heights, and thus determines a 1-morphism $f_{\C}$ in $\C$ given as a
composite of the generating 1-morphisms.  Let $F(f_{\T}) = f_{\C}$.
Due to the level-preserving property of equivalence isotopies, two
generic tangles representing the same 1-morphism in $\T$ cannot differ
by Reidemeister moves or by changing the order of heights of crossings
or extrema, so $F$ is well defined on 1-morphisms.   Using standard
techniques one can construct an equivalence isotopy between any pair of
tangles whose equivalence classes are mapped by $F$ to the same
1-morphism in $\C$, so $F$ is injective on 1-morphisms.   Finally, every
1-morphism in $\C$ can be realized as the image of a 1-morphism in $\T$.
We conclude that $F$ is bijective on 1-morphisms.  One can also check that
$F$ as defined on objects and 1-morphisms is in fact a functor from the
underlying category of $\T$ to the underlying category of $\C$.

Let $\D$ be the 2-category with the same underlying category as
$\C$, but with 2-morphisms freely generated as horizontal and vertical
composites of the generating 2-morphisms listed in Section
\ref{2-morphisms}.   Each 2-morphism $\alpha_{\T}$ in $\T$ is represented
by some generic 2-tangle $S$.  By Theorem 3.5.4 of Carter, Rieger and
Saito \cite{CRS}, each singularity of the projection of this 2-tangle to
the square $I_3\times I_4$ corresponds to a generating 2-morphism for
$\C$.    Furthermore, the proof of this theorem yields a procedure for
assigning to $S$ a unique 2-morphism $\alpha_{\D}$ in $\D$.

Now suppose that $S'$ is another generic 2-tangle representing the
2-morphism $\alpha_{\T}$, and let $\alpha'_{\D}$ be the corresponding
2-morphism in $\D$.  Then there is an equivalence isotopy carrying $S$
to $S'$.  By Theorem 3.5.5 of Carter, Rieger and Saito, it follows that
one can go from $\alpha_{\D}$ to  $\alpha'_{\D}$ using the relations
listed in Section \ref{2-morphisms}, together with the equational laws
in the definition of a 2-category.   It follows that $\alpha_{\D}$ and
$\alpha'_{\D}$ map to the same 2-morphism $\alpha_{\C}$ under the
canonical 2-functor from $\D$ to $\C$.  Thus we may define $F$ on
2-morphisms by $F(\alpha_{\T}) = \alpha_{\C}$.   Theorem 3.5.4 of
Carter, Rieger and Saito implies that $F$ is surjective on 2-morphisms,
while their Theorem 3.5.5 implies that it is injective.  One can easily
check that $F \maps \T \to \C$ is a 2-functor.  \qed

Since $\T$ is a braided monoidal 2-category with duals, we can use the
isomorphism $F\maps \T \to \C$ to give $\C$ the structure of a braided
monoidal 2-category with duals.  We then have the following monoidal,
braiding and duality structures on $\C$:

\begin{enumerate}

\item $I = A_0$.

\item $A_m\tensor A_n = A_{m+n}$.

\item $A_n\tensor Y_{i,j} = Y_{n+i,j}$ and $Y_{i,j} \tensor A_n =
Y_{i,j+n}$ for $Y=X,\bar{X}, \cap$ or $\cup$.  $A_n\tensor 1_m =
1_n\tensor A_m = 1_{n+m}$.  The tensor products of objects with  other
1-morphisms are determined by the relations
\[  A \tensor (fg) = (A \tensor f)(A \tensor g) \]
and
\[ (fg) \tensor A = (f \tensor A)(g \tensor A) .\]

\item  $\btensor_{Y_{i,j},Z_{m,n}} = N_{Y_{i,j + z},Z_{y +
m,n}}$, $\btensor{}_{Y_{i,j},1_n} = 1_{Y_{i,j+n}}$ and
$\btensor_{1_n,Y_{i,j}} = 1_{Y_{n+i,j}}$ for $Y, Z \in \{ X, \bar{X},
\cup, \cap\}$, where $\source (Z_{m,n}) = A_z$ and $\target (Y_{i,j}) =
A_y$.  The tensorator is determined for other 1-morphisms by the
relations
\[{\btensor}_{f,gg'} =
((A\tensor g) {\btensor}_{f,g'})\cdot ({\btensor}_{f,g} (A'\tensor g')) \]
for $f\colon A\to A'$, $g\colon B\to B'$ and $g'\colon
B'\to B''$, and
\[{\btensor}_{ff',g}= ({\btensor}_{f,g}(f'\ten B'))\cdot ((f\ten B)
{\btensor}_{f',g}) \]
for $f\colon A\to A'$, $f'\colon A'\to A''$ and
$g\colon B\to B'$.

\item $A_n\tensor Y_{i,j} = Y_{n+i,j}$ and $Y_{i,j} \tensor A_n =
Y_{i,j+n}$ for $Y = I, E, W, \two, S_{i,}, H, T$ and the $\,\bar{}\,$,
$*$, or $\dagger$ of these 2-morphisms; $A_n\tensor N_{Y_{i,j},Z_{k,l}}
= N_{Y_{n + i,j},Z_{n + k,l}}$ and $N_{Y_{i,j},Z_{k,l}}\tensor A_n =
N_{Y_{i,j+n},Z_{k,l+n}}$; $A_n\tensor 1_f = 1_{A_n\tensor f}$ and $1_f
\tensor A_n = 1_{f\tensor A_n}$.  The tensor products of objects with
other 2-morphisms are determined by the relations
\[ A \tensor (\alpha \beta) = (A \tensor \alpha)(A \tensor \beta)\]
and
\[ (\alpha \beta) \tensor A = (\alpha \tensor A)(\beta \tensor A) .\]

\item $R_{A_0,A_n} = R_{A_n,A_0} = 1_n$, $R_{A_1,A_1} = X_{0,0}$,
and $R_{A_n,A_m}$ is determined for other $n,m$ by the relations
\[ R_{A_1,A_{m + n}} = (R_{A_1, A_m}\tensor A_n)(A_m\tensor R_{A_1,A_n})\]
and
\[ R_{A_{m + n},A_i} = (A_m\tensor R_{A_n,A_i})(R_{A_m,A_i}\tensor A_n). \]

\item $\tilde{R}_{(A,B|C)} = 1_{(A\tensor R_{B,C})(R_{A,C}\tensor B)}$.

\item $\tilde{R}_{(A_1|B,C)} = 1_{(R_{A_1,B}\tensor C)(B\tensor
R_{A_1,C})}$, while $\tilde{R}_{(A_n|B,C)}$ is determined for other $n$
using condition $((\bullet \tensor \bullet) \tensor (\bullet \tensor
\bullet))$ in the definition of braided monoidal 2-category.

\item $R_{A_1,X_{0,0}} = S^*_{0;0,0}$, $R_{X_{0,0},A_1} = S_{0;0,0}$,
$R_{A_1,\bar{X}_{0,0}} = S^*_{1;0,0}$, $R_{\bar{X}_{0,0},A_1} =
\bar{S}_{2;0,0}$,

$R_{A_1,\cap_{0,0}} = (\bar{H}^*_{0,0} X_{1,0})\cdot
(\cap_{0,1}\bar{\two}^*_{1,0})$,
$R_{\cap_{0,0},A_1} = (H_{0,0}X_{0,1})\cdot
(\cap_{1,0}\bar{\two}^*_{0,1})$,

$R_{A_1,\cup_{0,0}}
= (\two_{0,1}\cup_{1,0})(X_{0,1} \bar{H}^\dagger_{0,0})$,
$R_{\cup_{0,0},A_1} =
(\two_{1,0}\cup_{0,1})(X_{1,0}H^{\dagger *}_{0,0})$,

$R_{A_1,1_n} = 1_{R_{A_1,A_n}}$ and $R_{1_n,A_1} = 1_{R_{A_n,A_1}}$,
while $R_{A,f}$ and $R_{f,A}$ for other 1-morphisms $f$ and objects $A$
are determined using conditions $((\bull \ten\bull)\ten\rarr)$,
$(\rarr\ten (\bull\ten\bull))$, $((\rarr\ten\bull)\ten\bull)$,
$((\bull\ten\rarr)\ten\bull)$, $(\bull\ten(\rarr\ten\bull))$, and
$(\bull\ten(\bull\ten\rarr))$ in the definition of braided monoidal
2-category.

\item $A_n^* = A_n$, and the duals of 1-morphisms and 2-morphisms
are given as in the definition of $\C$.

\item $i_{A_0} = 1_{A_0}$ and $i_{A_1} = \cap_{0,0}$, while the unit
is determined for other objects by the relation
\[  i_{A\tensor B} = i_A (A\tensor i_B \tensor A^*).  \]

\item $e_{A_0} = 1_{A_0}$ and $e_{A_1} = \cup_{0,0}$, while the counit
is determined for other objects by the relation
\[  e_{A\tensor B} = (A^*\tensor e_B \tensor A) e_A .\]

\item $i_{\cap_{m,n}} = I_{m,n}$, $i_{\cup_{m,n}} = E^*_{m,n}$,
$i_{X_{m,n}} = \two_{m,n}$, $i_{\bar{X}_{m,n}} = \bar{\two}_{m,n}$,
$i_{1_n} = 1_{1_n}$.   For other 1-morphisms, the unit is determined
using the relation
\[   i_{fg} = i_f\cdot (f i_g f^*) .\]

\item $e_{\cap_{m,n}} = E_{m,n}$, $e_{\cup_{m,n}} = I^*_{m,n}$,
$e_{X_{m,n}} = \bar{\two}^*_{m,n}$, $e_{\bar{X}_{m,n}} = \two^*_{m,n}$,
$e_{1_n} = 1_{1_n}$.   For other 1-morphisms, the unit is determined
using the relation
\[   e_{fg} = (g^* e_f g)\cdot e_g .\]

\item $T_{A_0} = 1_{1_0}$,  $T_{A_1} = T_{0,0}$, and $T_{A_n}$ is
determined for other $n$ using the relation
\[T_{A\tensor B} = ((i_A\ten A\ten B)(A\ten {\btensor}^{-1}_{i_B,
e_A}\ten B)(A\ten B\ten e_B))\cdot ((T_A\ten B)(A\ten T_B)).\]

\end{enumerate}

\noindent Moreover, the definition of the $\dagger$ operation on
2-morphisms of $\C$ agrees with the general definition of `adjoint' for
2-morphisms in a monoidal 2-category with duals.  For example, for the
generator $H_{m,n}$, this follows from relations 23 and 28 in Section
\ref{2-morphisms}.

In what follows we sometimes apply the notation developed for $\C$ to
the 2-category $\T$, using the isomorphism $F$.  Note in particular that
$A_1 \in \C$ is an unframed self-dual object with writhing $W_{0,0}$.

\section{The Universal Property of 2-Tangles}

In this section, we use the isomorphism between $\T$ and $\C$ to
characterize $\T$ in terms of a universal property.  Namely, we show
that $\T$ is the free braided monoidal 2-category with duals on the
unframed self-dual object $Z$ corresponding to a single point in the
unit square.  To do this, we first define what it means for a braided
monoidal 2-category with duals to be `generated' by an unframed
self-dual object, and show that $\T$ is generated by $Z$ in this sense.
To describe the sense in which $\T$ is `freely' generated by $Z$, we
define what it means for a strict monoidal 2-functor to `semistrictly
preserve braiding and duals on an unframed self-dual generating object'.
Then we show that for any braided monoidal 2-category with duals $\A$
containing an unframed self-dual object $A$, there exists a unique
strict monoidal 2-functor $F\maps \T \to \A$ with this property that
maps $Z$ to $A$.   This universal property characterizes $\T$ up to
isomorphism.

\begin{defn}\label{generated}\et We say a braided monoidal 2-category is
{\rm generated} by an unframed self-dual object $A$ if:

\begin{enumerate}
{\rm
\item Every object can be obtained by tensoring from:
\begin{alphalist}
\item $I$,
\item $A$.
\end{alphalist}
\item Every 1-morphism can be obtained by composition from:
\begin{alphalist}
\item $1_A$,
\item $i_A$,
\item $R_{A,A}$,
\item tensor products of the above 1-morphisms with arbitrary objects,
\item duals of the above 1-morphisms.
\end{alphalist}
\item Every 2-morphism can be obtained by horizontal and vertical
composition from:
\begin{alphalist}
\item 2-morphisms $1_f$ for arbitrary 1-morphisms $f$,
\item 2-morphisms $\bigotimes_{f,g}$ for arbitrary 1-morphisms $f$ and $g$,
\item 2-morphisms $R_{A,f}$ and $R_{f,A}$ for the 1-morphisms $f$ listed
in a) - e) above,
\item 2-morphisms $i_f$ for arbitrary 1-morphisms $f$,
\item $T_A$,
\item $W_A$,
\item tensor products of arbitrary objects with the above 2-morphisms,
\item duals of the above 2-morphisms.
\end{alphalist}
}
\end{enumerate}
\end{defn}

\begin{thm}\et $\T$ is a braided monoidal 2-category generated by the
unframed self-dual object $Z$. \end{thm}

Proof - Using the isomorphism between $\T$ and $\C$, the definition of
$\C$ immediately implies that every object in $\T$ is either $I$ or a
tensor product of copies of $Z$.  It also implies that every 1-morphism
is a composite of the 1-morphisms $1_Z, i_Z, R_{Z,Z}$, tensor products
of these 1-morphisms with objects of $\T$, and duals thereof.

Similarly, we immediately see that any 2-morphism in $\T$ can be
obtained by horizontal and vertical composition from the 2-morphism
generators listed in Section \ref{2-morphisms}.    Thus it suffices to
describe these generators as horizontal and vertical composites of the
2-morphisms listed in clauses 3a) - 3h) of Definition \ref{generated}.
We can simplify this task with the help of a few observations.   First,
the case of identity 2-morphisms is trivial.  Second, for any 2-morphism
generator $\alpha$, the corresponding 2-morphisms $\alpha^*$ and
$\alpha^\dagger$ as given in Section \ref{2-morphisms} are the same as
the dual and adjoint in the 2-category sense, so we do not need to
describe these variants of $\alpha$.  Similarly, we do not need to
describe the variant $\bar{\alpha}$ when it equals $\alpha$, and we do
need to describe $\bar{N}_{Y,Z}$, since $\bar{N}_{Y,Z} =
N_{\bar{Y},\bar{Z}}$.  Finally, since $\alpha_{m,n} = A_m\ten
\alpha_{0,0}\ten A_n$, it suffices to describe the following
2-morphisms:

\begin{enumerate}
\item $I_{0,0} = i_{i_Z}$

\item $E_{0,0} = e_{i_Z}$

\item $W_{0,0} = W_Z$,
$\bar{W}_{0,0} = (i_Z i_{R_{Z,Z}})\cdot (W^*_Z R^*_{Z,Z})$

\item $\two_{0,0} = i_{R_{Z,Z}}$, $\bar{\two}_{0,0} = i_{R^*_{Z,Z}}$

\item $S_{0;0,0} = R_{R_{Z,Z},Z}$,
$S_{1;0,0} = R^*_{Z,R^*_{Z,Z}}$,
$S_{2;0,0} = R^\dagger_{Z,R^*_{Z,Z}}$

$\bar{S}_{0;0,0} = R^\dagger_{Z,R_{Z,Z}}$,
$\bar{S}_{1;0,0} = R^{\dagger *}_{R^*_{Z,Z},Z}$,
$\bar{S}_{2;0,0} = R_{R^*_{Z,Z},Z}$

\item $H_{0,0} = ((i_Z\tensor Z)(Z\tensor R_{Z,Z})(i_{R_{Z,Z}}\tensor Z))
\cdot (R_{i_Z,Z}(R^*_{Z,Z}\tensor Z))$

$\bar{H}_{0,0} = ((i_Z\ten Z)(Z\ten R^*_{Z,Z})(i_{R^*_{Z,Z}}\ten Z))
\cdot (R^\dagger_{Z,i^{*}_Z}(R_{Z,Z}\ten Z))$

\item $T_{0,0} = T_Z$

\item $N_{Y_{m,n}, Z_{i,j}} = \btensor_{Y_{m,n'}, Z_{0,j}}$ for
some $n'$

\end{enumerate}

All the above are composites of the 2-morphisms listed in clauses 3a) -
3h) of Definition \ref{generated}, so the braided monoidal 2-category
$\T$ is generated by the unframed self-dual object $Z$.  \qed

Now we describe the sense in which $\T$ is `freely' generated by the
unframed self-dual object $Z$.  Naively, one might hope that for any
braided monoidal 2-category $\B$ containing an unframed self-dual object
$B$, there would exist a unique 2-functor $F \maps \T \to \B$ sending
$Z$ to $B$ and strictly preserving all the structure in sight, at least
on the generating object: the monoidal structure, braiding and duals.
However, this is too much to ask.  Uniqueness is no problem, but such a
2-functor might not exist, because $\tilde{R}_{(Z|Z,Z)}$ and
$\tilde{R}_{(Z,Z|Z)}$ are identity 2-morphisms, while this might not be
true of $\tilde{R}_{(B|B,B)}$ and $\tilde{R}_{(B,B|B)}$.  One way to
deal with this is to include conditions ensuring this in the definition
of `unframed self-dual object'.  This is basically the approach we took
in our earlier short paper \cite{BL}.

While this approach is consistent, it seems unnecessarily restrictive to
impose these extra conditions on the object $B$.  Probably there is a
strictification theorem saying that these conditions represent no
essential loss of generality.  However, such a theorem has not yet been
proved.  The generalized center construction of HDA1 shows that one can
safely assume either that $\tilde{R}_{(\cdot|\cdot,\cdot)}$ or
$\tilde{R}_{(\cdot,\cdot|\cdot)}$ is the identity.   Unfortunately, we
do not see how to use it to simultaneously set {\it both} these braiding
coherence 2-morphisms equal to the identity, even for a single object.

The approach we take here is thus to weaken our insistence that $F \maps
\T \to \B$ strictly preserve all the structure in sight.   This allows
for the possibility that $\tilde{R}_{(B|B,B)}$ and $\tilde{R}_{(B,B|B)}$
are not identity 2-morphisms.  As a result, these 2-morphisms appear as
`padding' in some of the equations in Definition \ref{semistrictly}.

\begin{defn}\et For monoidal 2-categories $\A$ and $\B$, a {\rm strict
monoidal 2-functor} $F\maps \A \to \B$ is a 2-functor such that $F(I) =
I$, $F(A\tensor X) = F(A)\tensor F(X)$ and $F(X\tensor A) = F(X)\tensor
F(A)$ for any object $A$ and object, 1-morphism or 2-morphism $X$, and
$F(\btensor_{f,g}) = \btensor_{F(f),F(g)}$ for any 1-morphisms $f$ and
$g$.   \end{defn}

\begin{defn}\label{semistrictly}\et Suppose that $\A,\B$ are braided
monoidal 2-categories with duals and that $\A$ is generated by an
unframed self-dual object $A$ with $\tilde{R}_{(A|A,A)} = 1$,
$\tilde{R}_{(A,A|A)} = 1$.  We say a strict monoidal 2-functor $F \maps
\A \to \B$ mapping $A$ to an unframed self-dual object $B \in \B$ {\rm
preserves braiding and duals semistrictly on the generator} if:

{\rm
\begin{enumerate}
\item $F(X^*) = F(X)^*$ for every object, morphism, or 2-morphism $X$,
\item $F(i_f) = i_{F(f)}$ for every morphism $f$,
\item $F(i_A) = i_B$,
\item $F(T_A) = T_B$,
\item $F(W_A) = W_B$,
\item $F(R_{A,A}) = R_{B,B}$,
\item
$F(R_{A,f}) = (B\tensor F(f)) \tilde{R}_{(B|B,B)})
\cdot R_{B,F(f)} \cdot (\tilde{R}^{-1}_{(B|B,B)}(F(f)\tensor B))$

and

$F(R_{f,A}) = ((F(f)\tensor B) \tilde{R}_{(B,B|B)})
\cdot R_{F(f),B} \cdot (\tilde{R}^{-1}_{(B,B|B)}(B\tensor F(f)))$

for $f = R_{A,A}, R^*_{A,A}$.

\item
$F(R_{A,i_A}) = ((B\tensor i_B)\tilde{R}_{(B|B,B)})\cdot R_{B,i_B}$

and

$F(R_{i_A,A})= ((i_B\tensor B) \tilde{R}_{(B,B|B)})\cdot R_{i_B,B}$

\item
$F(R_{A,e_A}) = R_{B,e_B}\cdot (\tilde{R}^{-1}_{(B|B,B)}(e_B \tensor B))$

and

$F(R_{e_A,A}) = R_{e_B,B} \cdot (\tilde{R}^{-1}_{(B,B|B)}(B\tensor e_B))$

\end{enumerate}
}
\end{defn}

\noindent  Note that in the above definition we still assume that
$\tilde{R}_{(A|A,A)} = 1$ and $\tilde{R}_{(A,A|A)} = 1$ for the unframed
self-dual object $A$ generating the 2-category $\A$.  We could drop
this assumption at the expense of still more padding in conditions 7
-- 9, but we do not need this extra generality.

The following theorem is the main result of this paper:

\begin{thm}\et For any braided monoidal 2-category with duals $\B$ and
unframed self-dual object $B \in \B$, there exists a unique strict monoidal
2-functor $F\maps \T \to \B$ with $F(Z) = B$ that preserves braiding
and duals semistrictly on the generator. \end{thm}

Proof - Uniqueness follows straightforwardly from the fact that $\T$ is
generated by $Z$, together with the fact that $F$ is a strict monoidal
2-functor preserving the braiding and duals semistrictly on the
generator.  Together these suffice to determine $F$ on any object,
morphism, or 2-morphism of $\T$.

For existence we use the isomorphism $\T \iso \C$ to describe $\T$ using
generators and relations, and show that all the relations are mapped by
$F$ to equations that actually hold in $\B$.  For objects and
1-morphisms there are no nontrivial relations.  For 2-morphisms we need
to check that $F(\alpha) = F(\beta)$ for every equation $\alpha = \beta$
in the list of 30 relations given in Section \ref{2-morphisms}.
In addition we need to show that $F(\alpha^*) = F(\beta^*)$,
$F(\alpha^\dagger) = F(\beta^\dagger)$, and $F(\bar{\alpha}) =
F(\bar{\beta})$.  The first two follow automatically from $F(\alpha) =
F(\beta)$, since
\[       F(\alpha^*) = F(\alpha)^* = F(\beta)^*  = F(\beta^*)  \]
and
\[  F(\alpha^\dagger) = F(\alpha)^\dagger = F(\beta)^\dagger =
F(\beta^\dagger),  \]
using the fact that if $\alpha \maps f \tto g$,
\ban
F(\alpha^\dagger)
&=&F((g^*i_f) \cdot (g^*\alpha f^*) \cdot (e_g f^*)) \\
&=&(F(g)^* i_{F(f)}) \cdot (F(g)^* F(\alpha) F(f)^*) \cdot (e_{F(g)} F(f)^*) \\
&=&F(\alpha)^\dagger .
\ean
Thus we only need to check the `barred' version $F(\bar{\alpha}) =
F(\bar{\beta})$.  We skip this in cases where $\bar{\alpha} = \alpha$
and $\bar{\beta} = \beta$.

In what follows we sketch how to show $F(\alpha) = F(\beta)$ and
$F(\bar{\alpha}) = F(\bar{\beta})$ for each equation $\alpha = \beta$ in
the list of relations in Section \ref{2-morphisms}.  We only use the
definition of a braided monoidal 2-category, the definition of an
unframed self-dual object, and the properties of $F$.  The full
arguments are lengthy and complicated, so we just indicate the main ideas.

1,2.  These follow from the unitarity of the writhing $W_B$.

3,4.  These follow from the unitarity of $i_{R_{B,B}}$.

5.  These and their barred versions follow from clause 7
in Definition \ref{semistrictly}, together with the fact that
$R_{R_{B,B},B}$, $R_{B,R_{B,B}}$, $\tilde{R}_{(B|B,B)}$,
and $\tilde{R}_{(B,B|B)}$ are unitary.    For example,
\ban
F(S_{0;0,0}) &=& F(R_{R_{Z,Z},Z})  \\
&=& ((R_{B,B} \ten B) \tilde{R}_{(B,B|B)}) \cdot
    R_{R_{B,B},B} \cdot (\tilde{R}^{-1}_{(B,B|B)} (B \ten R_{B,B}).  \ean
Being a composite of unitary 2-morphisms, $F(S_{0;0,0})$ is unitary.
Using the relations between tensoring and duality, it follows that
$F(S_{0;n,m})$ is unitary for all $n,m$.

6.  Let us write $\alpha = \beta$ for any of these relations holding in
$\T$.  These are all variants of the Zamolodchikov tetrahedron
equation.  Because of clause 7 in Definition \ref{semistrictly}, the
corresponding equations $F(\alpha) = F(\beta)$ contain extra
`padding' built from the 2-isomorphisms $\tilde{R}_{(B|B,B)}$,
$\tilde{R}_{(B,B|B)}$, and their inverses.  However, one can cancel
out this padding, reducing $F(\alpha) = F(\beta)$ to the corresponding
version of the Zamolodchikov tetrahedron equation for the object $B$ in 
$\B$.  This, in turn, follows from the definition of a braided monoidal
2-category (see Kapranov and Voevodsky \cite{KV} and the comments in
HDA1).  The barred versions $F(\bar{\alpha}) = F(\bar{\beta})$
follow from the fact that $\bar{R}_{f,A} = R^\dagger_{A,\bar{f^{*}}}$
and $\bar{R}_{A,f} = R^\dagger_{\bar{f^{*}},A}$.

7.   These follow from conditions $(\bullet \tensor \Downarrow)$  and
$(\Downarrow \tensor \bullet)$ in the definition of braided monoidal
2-category, applied to the 2-morphisms $W_B$ and $\bar{W}_B$.  Again,
the equation in $\B$ contains padding built from $\tilde{R}_{(B|B,B)}$,
$\tilde{R}_{(B,B|B)}$, and their inverses, but this padding cancels.
And again, the barred versions follow from the fact that $\bar{R}_{f,A}
= R^\dagger_{A,\bar{f}}$ and $\bar{R}_{A,f} = R^\dagger_{\bar{f},A}$.

8,9.  These follow from the unitarity of the triangulator $T_B$.

10.  This follows from the swallowtail coherence law.

11.  This and its barred version follow from clause 6 in Theorem 
19, together with the fact that $i_{R_{B,B}}$,
$i_{R^*_{B,B}}$, $R_{i_B,B}$, $R_{B,i_B}$, $\tilde{R}_{(B,B|B)}$, and
$\tilde{R}_{(B|B,B)}$ are unitary.

12.  This follows from the coherence law satisfied by the writhing
$W_B$.

13.  This follows from condition $(\Downarrow \tensor \bullet)$ in the
definition of braided monoidal 2-category, applied to the 2-morphism
$T$.   The barred version works similarly.


14.  Half of these and their barred versions follow from condition $(\to
\ten \to)$ in the definition of braided monoidal 2-category,
applied to the 1-morphisms $R_{B,B}$ or its dual and $i_B$ or its
dual.  The rest follow from $(\Downarrow \ten \bullet)$ and
$(\bullet \ten \Downarrow)$ applied to the 2-morphisms $R_{Z,i_Z}$
and $R^\dagger_{e_Z,Z}$, together with extensive use of the
monoidal 2-category axioms and equation 25 below.  The equations
to be proved contain padding built from $\tilde{R}_{(B|B,B)}$,
$\tilde{R}_{(B,B|B)}$, and their inverses, but this padding cancels.

15.  These follow from the unitarity of the tensorator
$\bigotimes_{f,g}$ for generating 1-morphisms $f,g$.

16.  This is essentially the Yang-Baxter equation for the tensorator,
which follows from conditions (vii) and (viii) in the definition of a
monoidal 2-category, using an argument analogous to the usual proof of
the Yang-Baxter equation for the braiding in a braided monoidal
category.  (Indeed, the latter is a special case of the former, since a
monoidal 2-category with one object is a braided monoidal category.)

17. These follow from conditions (vi) and (vii) in the definition of a
monoidal 2-category, applied to generating 1-morphisms and the
2-morphisms $R_{R_{B,B},B}$, $R^*_{B,R^*_{B,B}}$, $R^\dagger_{B,R^*_{B,B}}$,
$R^\dagger_{Z,R_{Z,Z}}$, $R^{\dagger *}_{R^*_{Z,Z},Z}$, and
$R_{R^*_{Z,Z},Z}$.  Again, the equations to be proved contain padding
built from $\tilde{R}_{(B|B,B)}$ and $\tilde{R}_{(B,B|B)}$, which cancels.

18.  These follow from conditions (vi) and (vii) in the definition of a
monoidal 2-category, applied to generating 1-morphisms and the
2-morphism $T_B$.

19.  These and their barred versions follow from conditions (vi) and
(vii) in the definition of a monoidal 2-category, applied to
generating 1-morphisms and the 2-morphisms $W_B$ and $\bar{W}_B$.

20.  These and their barred versions follow from conditions (vi)  and
(vii) in the definition of a monoidal 2-category, applied to  generating
1-morphisms and the 2-morphisms $H_{B,B}$ and  $\bar{H}_{B,B}$.  Again,
the equations to be proved contain padding built from
$\tilde{R}_{(B|B,B)}$ and $\tilde{R}_{(B,B|B)}$, which cancels.

21.  This and its barred version follow from condition $(\to \ten \to)$
in the definition of braided monoidal 2-category, applied to the
1-morphisms $i_B$ and $e_B$, together with extensive use of the
2-category axioms, the triangle equations, and the formulas
\[  R_{e_B,B} = (i_{R_{B \tensor B,B}} (e_B \tensor B)) \cdot
                (R_{B \tensor B,B} R^{\dagger *}_{i_B,B}),  \]
\[  R_{i_B,B} = (R^{\dagger *}_{e_B,B} R_{B \tensor B,B}) \cdot
                ((B \tensor i_B) e_{R_{B \tensor B,B}}) ,\]
together with similar formulas for $R_{B,e_B}$ and $R_{B,i_B}$.
Again, the equations to be proved contain padding built from
$\tilde{R}_{(B|B,B)}$ and $\tilde{R}_{(B,B|B)}$, which cancels.

22. These and their barred versions follow from conditions (vi) and
(vii) in the definition of a monoidal 2-category, applied to
generating 1-morphisms and the 2-morphisms $i_{R_{B,B}}, i_{R^*_{B,B}},
i_{i_B},$ and $i_{e_B}$.  Again, the equations to be proved contain
padding built from $\tilde{R}_{(B|B,B)}$ and $\tilde{R}_{(B,B|B)}$,
which cancels.

23.  These and their barred versions follow from the triangle equation
in the definition of a monoidal 2-category with duals, $(i_f f) \cdot
(f e_f) = 1_f$, applied to the 1-morphisms $R_{B,B}, R^*_{B,B},$
and $i_B$.

24.  This follows from the triangle equation $(f^* i_f) \cdot
(e_f f^*) = 1_{f^*}$ applied to the 1-morphism $e_B$, making
extensive use of the monoidal 2-category axioms.

25.  These and their barred versions follow from conditions $(\bullet
\tensor \Downarrow)$  and $(\Downarrow \tensor \bullet)$ in the
definition of braided monoidal 2-category, applied to the 2-morphisms
$i_{R_{B,B}}$ and $i_{R^*_{B,B}}$.    For some the proof is fairly
straightforward; for the rest one must cleverly exploit the monoidal
2-category axioms, the invertibility of $i_{R_{B,B}}$ and
$i_{R^*_{B,B}}$, and the definition of adjoint 2-morphisms.  In all
cases, the equations to be proved contain padding built from
$\tilde{R}_{(B|B,B)}$ and $\tilde{R}_{(B,B|B)}$, which cancels.

26.   This and its barred version follow from the coherence law
for the writhing $W_B$, with the help of Lemma \ref{writhing.lemma}.

27.  This and its barred version follows from the equation
$\alpha^{\dagger *} = \alpha^{* \dagger}$ in the definition of a
monoidal 2-category with duals applied to the 2-morphism $W_B$,
together with the coherence law for the writhing $W_B$, the
triangle equations, and extensive use of the 2-category axioms.

28.  This follows from condition $(\Downarrow \tensor \bullet)$ in the
definition of braided monoidal 2-category applied to the 2-morphism
$i_{i_B}$, together with extensive use of the 2-category axioms and the
formulas for $R_{e_B,B}$ and $R_{i_B,B}$ used in the proof of 21.
Again, the equation to be proved contains padding built from
$\tilde{R}_{(B|B,B)}$ and $\tilde{R}_{(B,B|B)}$, which cancels.  The
barred version works similarly.

29.  This equation is actually redundant, since it follows from 23 and 27
using only the 2-category axioms.

30.  This and its barred version follow from conditions $(\bullet
\tensor \Downarrow)$ and $(\Downarrow \tensor \bullet)$ in the
definition of braided monoidal 2-category, applied to the 2-morphism
$e_{i_B}$.  The equations to be proved contain padding built from
$\tilde{R}_{(B|B,B)}$ and $\tilde{R}_{(B,B|B)}$, which cancels.
 \hbox{\hskip 30em} \qed

\section{Conclusions}

We have shown that the 2-category of 2-tangles in 4 dimensions is the
free braided monoidal 2-category with duals on an unframed self-dual
object.   Given any self-dual unframed object in a braided monoidal
2-category with duals, our proof of this result gives a concrete recipe
for computing a 2-tangle invariant.  Of course, for this to be useful,
we need more examples of braided monoidal 2-categories with duals.
Obtaining these will require further work in higher-dimensional algebra.
There are a number of promising strategies.

First, there is plenty of evidence that a certain class of braided
monoidal 2-categories with duals, the braided monoidal 2-groupoids,
are essentially the same as homotopy 2-types of double loop spaces
\cite{BD,BD3}.   These should give 2-tangle invariants with a
`purely homotopy-theoretic' flavor.  In particular, for any compact
2-dimensional submanifold $\Sigma \subset \R^4$, these invariants
should depend only on the homotopy type of $\R^4 - \Sigma$.

Second, one can construct braided monoidal 2-categories as the `quantum
doubles' of monoidal 2-categories \cite{BN,C}.  It seems plausible
that applying this construction to a monoidal 2-category with duals
will give a braided monoidal 2-category with duals.  This reduces the
question to obtaining monoidal 2-categories with duals.  The 2-category
of unitary representations of a 2-groupoid should be a monoidal
2-category with duals, just as the category of unitary representations
of a groupoid is a monoidal category with duals \cite{B}.   Moreover,
the 2-category of representations of any Hopf category is a
monoidal 2-category \cite{Neuchl}, and when `unitary' representations
can be be defined, the 2-category of unitary representations should
be a monoidal 2-category with duals.   Some examples of Hopf categories
and related structures have been studied by Neuchl \cite{Neuchl} as well
as Crane and Yetter \cite{CY,CY2}.  Also, Crane and Frenkel have sketched an
interesting construction of a Hopf category from Kashiwara and Lusztig's
canonical basis of a quantum group \cite{CF}.

Third, just as one can construct braided monoidal categories
from solutions of the Yang-Baxter equation, one can construct braided
monoidal 2-categories from solutions of the Zamolodchikov tetrahedron
equations \cite{KV}.  Many such solutions are known \cite{CS}, so one
may hope that some give braided monoidal 2-categories with duals.

Finally, one expects `braided monoidal 3-Hilbert spaces' to be
interesting examples of braided monoidal 2-categories with duals
\cite{B}.  However, to obtain these we will probably need to use
some of the constructions sketched above.

Our result and its proof can probably be improved in various ways.
First, we expect similar algebraic characterizations of the 2-category
of framed and/or oriented 2-tangles in 4 dimensions, where we drop the
conditions that the object be unframed and/or self-dual.  Of course, our
current definition of `unframed' applies only to self-dual objects, so
we need to more clearly separate these concepts.  Moreover, our
definition of `braided monoidal 2-category with duals' may need some
extra conditions to handle framed tangles.   For example, the 2-morphism
corresponding to the framed Reidemeister I move exists under the current
definition of a braided monoidal 2-category with duals, without using
the writhing, but we are unable to use this definition to show that
this 2-morphism is unitary, even though topological considerations
say it should be.

Second, one should be able to compress the definition of `monoidal
2-category with duals' using more of the language of 2-category theory.
Doing so will shed more light on the still mysterious general notion of
`$n$-category with duals'.  It bodes well that the triangulator and
its swallowtail coherence law have already been observed by Street in
his study of adjunctions between 2-categories \cite{Street3}.  In
general, we expect a close relation between the theory of $n$-categories
with duals and the theory of adjunctions between $n$-categories.  This
has already been noted in work on 2-Hilbert spaces \cite{B,M}, and
the patterns found here should continue for higher $n$-Hilbert spaces.

Finally, and most importantly, there must be a way to state and prove
the tangle hypothesis for all values of $n$ and $k$ that does not
involve long lists of equations.  Our treatment of the case $n = k = 2$
resembles moving a house across the country by taking it apart, sending
the pieces by mail, and then rebuilding it on the other side.  First
Carter, Saito and Rieger deduced their list of movie moves  using a
classification of singularities.  Then we showed that their movie moves
are equivalent to a long form of the definition of `braided monoidal
2-category with duals generated by an unframed self-dual object'.  But
this long definition can presumably expressed much more tersely using
more sophisticated higher-dimensional algebra.  There should thus be a
more conceptual approach that proceeds at a higher level of abstraction.
While already desirable for $n = k = 2$, the advantages of such an
approach will be even greater for larger $n$ and $k$.

Finding a more conceptual approach to the tangle hypothesis poses many
interesting challenges in $n$-category theory.  In particular, it will
require a deeper understanding of the mysterious relationship between
$n$-categories and singularity theory.  We hope the present work
provides some useful clues.

\subsection*{Acknowledgments}
Many of the ideas underlying this work were developed with James Dolan.
We thank Marco Mackaay and Ross Street for useful discussions on
duality in monoidal 2-categories, and we are especially grateful to
Scott Carter and Masahico Saito for invaluable correspondence
regarding 2-tangles.  We also thank the referees for their careful
reading of this paper.  

\section*{Errata}

Our previous paper \cite{BL} has the following errors:

\begin{enumerate}
\item  The definition of an `unframed object' $A$ should include the
requirement that $A^\ast = A$.  All appearances of $A^*$ in this
definition can thus be replaced by $A$.

\item  The pictures in Figure 1 and Figure 3 should be switched.
\end{enumerate}

\end{document}